%
%
%


\hsize=14.5cm \vsize=20.5cm \hoffset=0.9cm \voffset=2cm
\baselineskip=13pt \footline={\hss{\vbox to
1cm{\vfil\hbox{\rm\folio}}}\hss}

\input amssym.def
\input amssym.tex

\font\csc=cmcsc10
\font\title=cmr17
\font\small=cmr8

\font\teneusm=eusm10
\font\seveneusm=eusm7
\font\fiveeusm=eusm5 \newfam\eusmfam
\def\eusm{\fam\eusmfam\teneusm} \textfont\eusmfam=\teneusm
\scriptfont\eusmfam=\seveneusm \scriptscriptfont\eusmfam=\fiveeusm

\font\teneufb=eufb10
\font\seveneufb=eufb7
\font\fiveeufb=eufb5 \newfam\eufbfam
\def\eufb{\fam\eufbfam\teneufb} \textfont\eufbfam=\teneufb
\scriptfont\eufbfam=\seveneufb \scriptscriptfont\eufbfam=\fiveeufb

\def\varGamma{{\mit \Gamma}}
\def\Re{{\rm Re}\,}
\def\Im{{\rm Im}\,}
\def\txt#1{{\textstyle{#1}}}
\def\scr#1{{\scriptstyle{#1}}}
\def\r#1{{\rm #1}}
\def\B#1{{\Bbb #1}}
\def\e#1{{\eusm #1}}
\def\f#1{{\eufb #1}}
\def\b#1{{\bf #1}}

\centerline{\title SUM FORMULA FOR KLOOSTERMAN SUMS}
\bigskip
\centerline{\title AND FOURTH MOMENT OF THE DEDEKIND}
\bigskip
\centerline{\title ZETA-FUNCTION OVER THE GAUSSIAN}
\bigskip
\centerline{\title NUMBER FIELD} \vskip 1cm \centerline{BY ROELOF W.
BRUGGEMAN AND YOICHI MOTOHASHI} \vskip 1cm
\noindent
{\csc Abstract}. We prove the Kloosterman--Spectral sum formula for
$\r{PSL}_2(\B{Z}[i])\backslash \r{PSL}_2(\B{C})$, and apply it to
derive an explicit spectral expansion for the fourth power moment of
the Dedekind zeta function of the Gaussian number field. Our sum
formula, Theorem 13.1, allows the extension of the spectral theory of
Kloosterman sums to all algebraic number fields.
\footnote{}{Mathematics Subject Classification: 11M41 11F72 11L05
22E30}
\vskip 1cm \centerline{\bf 1.
Introduction}
\smallskip
\noindent
Our principal aim is to establish an explicit formula for the fourth
moment of the Dedekind zeta-function $\zeta_\r{F}$ of the Gaussian
number field $\r{F}=\B{Q}(i)$:
$$
{\eusm Z}_2(g,\r{F})=\int_{-\infty}^\infty
|\zeta_\r{F}(\txt{1\over2}+it)|^4g(t)dt,\leqno(1.1)
$$
where the weight function $g$ is assumed, for the sake of simplicity,
to be entire and of rapid decay in any fixed horizontal strip. The
basic implement to be utilized is a sum formula for Kloosterman sums
over $\r{F}$
$$
S_\r{F}(m,n;\,c)=\sum_{\scr{d\bmod c}\atop\scr {(d,c)=1}}\exp\Big(2\pi
i\Re((md+n\tilde{d})/c))\Big) \leqno(1.2)
$$
with $c,d,m,n\in\B{Z}[i]$ and $d\tilde{d}\equiv1\bmod c$. It is a
genuine counterpart of Kuznetsov's sum formula for Kloosterman sums
over the rationals. As is to be detailed shortly, the existing
version of this sum formula concerns only $K$-trivial automorphic
forms over $\r{PSL}_2(\B{C})$, and is unsuitable to handle sums of
$S_\r{F}$ explicitly. We extend it to all $K$-types and invert the
Bessel transformation occurring in it. We stress that this inversion
has so far been obtained only for $\r{PSL}_2(\B{R})$ and
infinitesimally isomorphic groups.
\smallskip
The explicit formula in Theorem 14.1 expresses $\e{Z}_2(g,\r{F})$ as a
sum of a term $M_\r{F}(g)$, given by an integral transform of $g$,
and a term based on the spectral decomposition of
$L^2(\r{PSL}_2(\B{Z}[i])\backslash \r{PSL}_2(\B{C}))$. It is a
generalization to $\r{F}$ of Theorem 4.2 of [30], which gives for the
fourth moment $\e{ Z}_2(g,\B{Q})$ of the Riemann zeta-function a
similar explicit formula, based on spectral data for
$L^2(\r{PSL}_2(\B{Z})\backslash \r{PSL}_2(\B{R}))$; see also Section
15 below. In [6], we have generalized that theorem to real quadratic
number fields with class number one. In all the three cases, the
fourth moment $\e{Z}_2(g,\cdot)$ is linked to spectral data via
Kloosterman sums. For $\r{F}$ the step from the fourth moment to sums
of $S_\r{F}$ is carried out in Section 2, and follows the same lines
as in the rational case. In this respect, the real quadratic case in
[6] is harder. There we had to deal with infinitely many units. As
far as we see, it is essential for the method to assume that the
class number is one, for both real and imaginary quadratic fields.
Thus one may say that these three cases are built on an essentially
common structure. There are, however, notable differences among them
as well. In the rational case, the term $M_{\B Q}(g)$ is the main
term, overshadowing the other explicitly spectral terms. For the
quadratic cases, the same does not hold. We shall briefly discuss
this peculiar fact for the present case in the final section. In the
real quadratic case, there is even less reason to call the term
corresponding to $M_\r{F}(g)$ the main term. Yet the outward
similarity among these spectral expansions of the moments is highly
remarkable.
\smallskip
Various sums of Kloosterman sums can be related to Fourier
coefficients of automorphic forms --- by a Kloosterman-Spectral sum
formula. In the case of the Riemann zeta function, this is
Kuznetsov's sum formula in [21], [22]. There a sum of rational
Kloosterman sums is expressed in terms of a bilinear form in Fourier
coefficients of automorphic forms over the upper half-plane. The test
functions on both sides of this equality are related by an integral
transformation, given by a Bessel function. In applications of the
sum formula, it is important to have control over this integral
transformation, and, in particular, to be able to invert it.
Kuznetsov did this in [21], [22] in an ingenious way. The sum formula
for the upper half-plane has been discussed at many places. A self
contained treatment along classical lines can be found in the first
two chapters of [30]. For a spectral formulation of the sum formula,
the version in [1] has the advantage to stress that the spectral data
are tied not to automorphic forms but, in fact, to irreducible
subspaces of the right regular representation of $\r{PSL}_2(\B{R})$
in $L^2(\r{PSL}_2(\B{Z})\backslash \r{PSL}_2(\B{R}))$.
\par
Sums of Kloosterman sums over a number field ask for generalization of
the sum formula to $\r{SL}_2$
over the field. For the purpose of [6] we could appeal to [4], where
the case of totally real number fields is treated. The restriction
 there to totally positive number fields is
due to the fact that for products of copies of $\r{PSL}_2(\B{R})$ the
Bessel transformation in the sum formula is not essentially more
difficult than for one copy.
\par
For our present group $\r{PSL}_2(\B{C})$ such a reduction to smaller
groups does not hold; in fact it is the first step of an induction,
and we have to start essentially from scratch. It is true that
Miatello and Wallach have given in [26] a wide generalization of the
sum formula, to Lie groups of real rank one. They have, however, good
reasons to restrict themselves to irreducible representations with a
$K$-trivial vector which are comparable to automorphic forms of
weight zero. The relevant integral transform is described by a power
series expansion of its kernel function; hence its behavior is only
known near the origin. Moreover, the restriction to $K$-trivial
representations makes it unlikely that the same kernel function can
be used to describe the inverse transformation. However, we need, for
the purpose of the present paper, a sum formula that relates sums of
the form
$$
\sum_{c\in\B{Z}[i]\backslash\{0\}}S_\r{F}(m,n;c)f(c)\leqno(1.3)
$$
to spectral data for rather arbitrary test functions $f$ on
$\B{C}\backslash\{0\}$. That is, we are given a sum of Kloosterman
sums to begin with, but not spectral expressions as in [26]. This
requires a good control of the relevant integral kernel in much the
same manner as Kuznetsov's theory allows us to do for
$\r{PSL}_2(\B{R})$. We achieve this by deriving a complete ---
including all $K$-types --- sum formula of Kuznetsov type for
$\r{PSL}_2(\B{C})$, with the discrete subgroup $\r{PSL}_2(\B{Z}[i])$.
The major part of the present article is devoted to the development
of such a generalization of the Kuznetsov sum formula. As is
mentioned above, it is the first Lie group other than $\r{PSL}_2(\B
{R})$ and its covering groups for which this has been carried out.
The results are stated Theorems 10.1 and 13.1. The former theorem can
be used to get information on spectral data; and it is the basis of
the latter, in which sums of the type $(1.3)$ are spectrally
decomposed. The integral transform in these formulae has a product of
two Bessel functions as its kernel, see $(6.21)$ and $(7.21)$.
Theorem 11.1 gives the inversion of the integral transformation. The
integral representation in Theorem 12.1 allows us to bound the kernel
function in $(7.21)$ in a practical way for applications, especially
to treat $\e{Z}_2(g,\r{F})$.
\par
Once the sum formula in Theorem 13.1 is available, we can proceed with
the study of $\e{Z}_2(g,\r{F})$. The general approach is the same as
in the rational case, but the computations are by far more involved,
as can be expected.
\smallskip
We could try to deal with general imaginary quadratic number fields
with class number one, but have exploited arithmetical
simplifications offered by the specialization to $\r{F}=\B{Q}(i)$.
The derivation of the sum formula as given in Theorems 10.1 and 13.1
could be carried out for any imaginary quadratic number field. If the
class number is larger than one, the contribution of the continuous
spectrum is more complicated. The discrete subgroup can be any
congruence subgroup, provided we have a Weil type bound of the
corresponding generalized Kloosterman sums; actually, any non-trivial
estimate suffices. Without such a bound, we would run into additional
technical difficulties, that would require the method of [25] for
their resolution.
\smallskip
\noindent
{\csc Remark.} Main results of the present article have been announced
in our note [5].
\smallskip
\noindent
{\csc Convention.} Notations become available at their first
appearances and will continue to be effective throughout the sequel.
This applies to those in the above as well. We stress two points
especially: (i) The terms left/right invariance/irreducibility are,
respectively, abbreviations for the invariance/irreducibility of the
relevant function space with respect to the left/right translations
by the elements of the group under consideration. (ii) There are
mainly two kinds of summation variables, rational and Gaussian
integers. The distinction between them will easily be made from the
context. Also, group elements, operators, and spaces appear as
variables. They are explicitly indicated if there is any danger of
confusion.
\medskip
\centerline{\bf 2. A sum of Kloosterman sums}
\smallskip
\noindent
The aim of this section is to reduce $\e{Z}_2(g,\r{F})$ to a sum of
$S_\r{F}$ with variable arguments and modulus, indicating the core of
the problem that we are going to resolve. We shall partly follow a
discussion developed in [30] on the same subject.
\smallskip
Thus, let $g$ be as in $(1.1)$. Closely related to $\e{Z}_2(g,\r{F})$
is the function
$$
\e{I}(z_1,z_2,z_3,z_4;g) =\int_{-\infty}^\infty\zeta_\r{F}(z_1+it)
\zeta_\r{F}(z_2+it)\zeta_\r{F}(z_3-it)\zeta_\r{F}(z_4-it)g(t)dt,
\leqno(2.1)
$$
where all $\Re z_j$ are larger than $1$. Shifting the contour upward
appropriately, this can be continued meromorphically to the whole of
${\Bbb C}^4$. It is regular in a neighbourhood of the point
$\r{p}_{1\over2}=({1\over2},{1\over2},{1\over2},{1\over2})$, and
$$
{\eusm Z}_2(g,\r{F})=\e{I}(\r{p}_{1\over2};g)
+a_0g(\txt{1\over2}i)+b_0g(-\txt{1\over2}i)
+a_1g'(\txt{1\over2}i)+b_1g'(-\txt{1\over2}i)\leqno(2.2)
$$
with certain absolute constants $a_0$, $a_1$, $b_0$, $b_1$ which could
be made explicit. On the other hand, expanding the integrand and
integrating term by term, we get
$$
\e{I}(z_1,z_2,z_3,z_4;g)={1\over16}\sum_{kl\not=0}
{\sigma_{z_1-z_2}(k)
\sigma_{z_3-z_4}(l)\over{|k|^{2z_1}|l|^{2z_3}}}\hat{g}
\big(2\log|l/k|\big)\leqno(2.3)
$$
with $k,l\in\B{Z}[i]$. Here
$$
\sigma_\nu(n)=\sigma_\nu(n,0),\quad
\sigma_\nu(n,p)={1\over4}\sum_{d|n}(d/|d|)^{4p}|d|^{2\nu}\leqno(2.4)
$$
with $p\in\B{Z}$ and the divisibility inside $\B{Z}[i]$, and
$$
\hat{g}(x)=\int_{-\infty}^\infty g(t)e^{ixt}dt.\leqno(2.5)
$$
Classifying the summands according as $k=l$ and $k\not=l$, we have, in
the region of absolute convergence,
$$
\leqalignno{ \e{I}(z_1,z_2,z_3,z_4;g,\r{F})&=
{{\zeta_\r{F}(z_1+z_3)\zeta_\r{F}(z_1+z_4)\zeta_\r{F}(z_2+z_3)
\zeta_\r{F}(z_2+z_4)}\over4\zeta_\r{F}(z_1+z_2+z_3+z_4)}
\hat{g}(0)&(2.6)\cr &+{1\over16}\sum_{m\ne0} |m|^{-2z_1-2z_3}
B_m(z_1-z_2,z_3-z_4; g^*(\cdot\,;z_1,z_3)),\cr
}
$$
where
$$
B_m(\alpha,\beta;h)=\sum_{n(n+m)\ne0}
\sigma_\alpha(n)\sigma_\beta(n+m)h(n/m),\leqno(2.7)
$$
and
$$
g^*(u;\gamma,\delta)
={\hat{g}(2\log|1+1/u|)\over{|u|^{2\gamma}|1+u|^{2\delta}}}.
\leqno(2.8)
$$
The first term on the right of $(2.6)$ is due to $(14.21)$.
\smallskip
In order to exploit the relation $(2.2)$, experience in the rational
case suggests that we should continue analytically the identity
$(2.6)$ to a neighbourhood of the point $\r{p}_{1\over2}$, and that
such a continuation should be accomplished via spectrally decomposing
the function $B_m(\alpha,\beta;g^*(\cdot;\gamma,\delta))$ with a sum
formula of Kuznetsov's type. We shall see, in the final section, that
this is indeed the case. The long process to reach there begins with
the following fact on the complex Mellin transform of $g^*$:
\smallskip
\noindent
{\bf Lemma 2.1.}\quad{\it We put, for $q\in\B{Z}$, $s\in\B{C}$,
$$
\tilde{g}_q(s;\gamma,\delta)={1\over2\pi}\int_{\B{C}^\times}
g^*(u;\gamma,\delta)(u/|u|)^{-q}|u|^{2s}d^\times\!u,\leqno(2.9)
$$
where $\B{C}^\times=\B{C}\backslash\{0\}$, and
$d^\times\!u=|u|^{-2}d_+\!u$ with the Lebesgue measure $d_+\!u$ on
$\B{C}$. Then $\tilde{g}_q(s;\gamma,\delta)$ is regular in the domain
$$
\Re (s-\gamma-\delta)<0\leqno(2.10)
$$
as a function of three complex variables. More precisely, all of its
singularities are in the set $\{\gamma+\delta+{1\over2}|q|+l:\B{Z}\ni
l\ge0\}$, as is implied by the representation
$$
\leqalignno{ \tilde{g}_q(s;\gamma,\delta)=&{1\over2}(-1)^q
{\Gamma(\gamma+\delta-s+{1\over2}|q|)\over
\Gamma(1+s-\gamma-\delta+{1\over2}|q|)}&(2.11)\cr
&\times\int_{-\infty}^\infty{\Gamma(1-\delta+it)\over
\Gamma(\delta-it)}{ \Gamma(s-\gamma-it+{1\over2}|q|)\over
\Gamma(1+\gamma-s+it+{1\over2}|q|)}g(t)dt, }
$$
where the contour separates the poles of $\Gamma(1-\delta+it)$ and
those of $\Gamma(s-\gamma-it+{1\over2}|q|)$ to the left and the
right, respectively; and $s,\,\gamma,\,\delta$ are assumed to be such
that the contour can be drawn. Moreover, if $\gamma$, $\delta$, and
$\Re s$ are bounded, then we have, regardless of $(2.10)$,
$$
\tilde{g}_q(s;\gamma,\delta)\ll (1+|q|+|s|)^{-A}\leqno(2.12)
$$
with any fixed $A>0$, as $|q|+|s|$ tends to infinity. }
\smallskip
\noindent
Proof. The first assertion follows from the observations that
$g^*(u;\gamma,\delta)\ll|u|^{-2\Re(\gamma+\delta)}$ as $u\to\infty$,
and that $g^*(u;\gamma,\delta)$ is of rapid decay as $u\to0,\,-1$,
which is a consequence of respective upward and downward shifts of
the contour in $(2.5)$. To prove the second assertion we assume,
temporarily, that
$$
\Re\gamma<\Re s<\Re(\gamma+\delta)<\Re\gamma+1, \leqno(2.13)
$$
which is of course contained in $(2.10)$. Moving to polar coordinates,
we have
$$
\tilde{g}_q(s;\gamma,\delta) ={1\over2\pi}\int_{-\infty}^\infty
g(t)\int_0^\infty r^{2(s-\gamma-it)-1}\int_{-\pi}^\pi
{e^{i|q|\theta}\over|1+re^{i\theta}|^{2(\delta-it)}}
d\theta\,dr\,dt.\leqno(2.14)
$$
This triple integral is absolutely convergent. Note that we need to
deal with the part corresponding to $|r-1|<\varepsilon$,
$|\theta\pm\pi|<\varepsilon$ with a small $\varepsilon>0$ separately.
The innermost integral is equal to
$$
\leqalignno{
&{1\over\Gamma(\delta-it)}\int_0^\infty y^{\delta-it-1}e^{-y(1+r^2)}
\int_{-\pi}^\pi\exp(i|q|\theta-2ry\cos\theta)d\theta\,dy&(2.15)\cr
&={2\pi(-1)^q\over\Gamma(\delta-it)}\int_0^\infty
y^{\delta-it-1}e^{-y(1+r^2)}I_{|q|}(2ry)dy. }
$$
Thus we have
$$
\leqalignno{ \tilde{g}_q(s;\gamma,\delta)&={(-1)^q}
\int_{-\infty}^\infty
{2^{-\delta+it}g(t)\over\Gamma(\delta-it)}\int_0^\infty
r^{2(s-\gamma)-\delta-it-1}&(2.16)\cr
&\hskip1cm\times\int_0^\infty
y^{\delta-it-1}e^{-{1\over2}y(r+1/r)}I_{|q|}(y)dy\, dr\,dt\cr
&={(-1)^q}\int_{-\infty}^\infty
{2^{1-\delta+it}g(t)\over\Gamma(\delta-it)}\cr
&\hskip1cm\times\int_0^\infty
y^{\delta-it-1}K_{2(s-\gamma)-\delta-it}(y)I_{|q|}(y)dy\,dt, }
$$
where the necessary absolute convergence follows from asymptotic
expansions of these Bessel functions. The last integral can be
evaluated as a limiting case of formula $(1)$ on p.\ 410 of [42]
coupled with the relation $i^{-|q|}J_{|q|}(iy)=I_{|q|}(y)$. It is
equal to
$$
\leqalignno{
&{\Gamma(s-\gamma-it+{1\over2}|q|)
\Gamma(\gamma+\delta-s+{1\over2}|q|)
\over2^{2-\delta+it}\Gamma(|q|+1)}&(2.17)\cr
&\hskip 1cm\times{}_2F_1( s-\gamma-it+\txt{1\over2}|q|,
\gamma+\delta-s+\txt{1\over2}|q|;|q|+1;1)\cr
=&2^{\delta-it-2}{\Gamma(1-\delta+it) \Gamma(s-\gamma-it+{1\over2}|q|)
\Gamma(\gamma+\delta-s+{1\over2}|q|)\over
\Gamma(1+\gamma-s+it+{1\over2}|q|)
\Gamma(1-\gamma-\delta+s+{1\over2}|q|)}, }
$$
where we have used Gauss' formula for the value of the hypergeometric
function at the point $1$. {}From these we get $(2.11)$ with the
contour $\Im t=0$, provided $(2.13)$. Having obtained this, we use
analytic continuation to have the representation $(2.11)$ for those
$s,\gamma,\delta$ with which the above separation of poles is
possible. Then, to find the location of singularities, we move the
contour in $(2.11)$ to $\Im t=L$ with a sufficiently large $L>0$. We
may encounter poles at $t=i(1-\delta)+il$ with integers $l\ge0$, but
the relevant residual contributions are easily seen to be entire over
$\B{C}^3$, whence the above assertion. As to the bound $(2.12)$, we
note that, when $|s|+|q|$ tends to infinity while $\gamma$, $\delta$,
and $\Re s$ are bounded, the contour in $(2.11)$ can be drawn. Then
we need only to push the contour down appropriately. The new integral
is readily estimated by Stirling's formula. We may encounter poles if
$q$ is bounded, but then resulting residues do not disturb $(2.12)$
because of the assumption on $g$. This ends the proof.
\smallskip
Now, returning to $(2.7)$, we note the Ramanujan identity over
$\r{F}$: We have, for any $n\in\B{Z}[i]$, $p\in2\B{Z}$,
$$
\sum_{c\not=0}(c/|c|)^{2p}S_\r{F}(n,0;c)
|c|^{-2s}={4\over\zeta_\r{F}(s,p/2)}\cdot \cases{\sigma_{1-s}(n,p/2)&
if $n\ne0$, $\Re s>1$,
\cr \zeta_\r{F}(s-1,p/2)& if $n=0$, $\Re s>2$,\cr}\leqno(2.18)
$$
where
$$
\zeta_\r{F}(s,p)={1\over4}\sum_{n\ne0}(n/|n|)^{4p}|n|^{-2s}\leqno(2.19)
$$
is the Hecke $L$-function of $\r{F}$ associated with the
Gr\"ossencharakter $(n/|n|)^{4p}$. Applying this with $p=0$ to the
factor $\sigma_\beta(n+m)$ in $(2.7)$, we see that if
$$
1+\max(0,\Re\alpha)
<\Re(\gamma+\delta),\quad \Re\beta<-1,\leqno(2.20)
$$
then we have the absolutely convergent expression
$$
\leqalignno{ B_m(\alpha,\beta;&g^*(\cdot\,;\gamma,\delta))=
{1\over4}\zeta_\r{F}(1-\beta)\sum_{c\ne0} |c|^{2(\beta-1)}&(2.21) \cr
&\times\sum_{\scr{d\,\bmod\,c}\atop\scr{(d,c)=1}} \exp(2\pi
i\Re(dm/c)) D_m(\alpha,d/c;g^*(\cdot\,;\gamma,\delta)) }
$$
with
$$
D_m(\alpha,d/c;g^*(\cdot\,;\gamma,\delta))
=\sum_{n\ne0}\sigma_\alpha(n) \exp(2\pi
i\Re(dn/c))g^*(n/m;\gamma,\delta). \leqno(2.22)
$$
Note that we have used $g^*(-1;\gamma,\delta)=0$ in $(2.21)$. We
expand $g^*(u;\gamma,\delta)$ into a Fourier series, and apply Mellin
inversion to each Fourier coefficient, so that in view of the last
lemma we have, for any $u\in\B{C}^\times$ and
$\tau<\Re(\gamma+\delta)$,
$$
g^*(u;\gamma,\delta)={1\over\pi i}\sum_{q\in\B{Z}}(u/|u|)^{q}
\int_{(\tau)}\tilde{g}_q(s;\gamma,\delta)|u|^{-2s}ds,\leqno(2.23)
$$
where $(\tau)$ is the vertical line $\Re s=\tau$. Thus we have
$$
D_m(\alpha,d/c;g^*(\cdot\,;\gamma,\delta))
={1\over{\pi{i}}}\sum_{q\in\B{Z}}(m/|m|)^{-q} \int_{(\tau)}
X_q(s,\alpha;d/c)|m|^{2s}\tilde{g}_q(s;\gamma,\delta)ds, \leqno(2.24)
$$
provided $1+\max(0,\Re\alpha)<\tau<\Re(\gamma+\delta)$, where
$$
X_q(s,\alpha;d/c)=\sum_{n\ne0}\sigma_\alpha(n) \exp(2\pi
i\Re(dn/c))(n/|n|)^q|n|^{-2s}\leqno(2.25)
$$
with $(d,c)=1$.
\par
Then we invoke (see [33]): If $q\ne0$, the function
$X_q(s,\alpha;d/c)$ of $s$ is regular for any $\alpha$. If $q=0$ and
$\alpha\ne0$, it is regular except for the simple poles at $s=1$ and
$s=1+\alpha$ with the residues
$\pi|c|^{2(\alpha-1)}\zeta_\r{F}(1-\alpha)$ and
$\pi|c|^{-2(\alpha+1)}\zeta_\r{F}(1+\alpha)$, respectively. Moreover,
we have, for any combination of parameters,
$$
\leqalignno{ X_q(s,\alpha;&d/c)
=(c/|c|)^{2q}(\pi/|c|)^{4s-2\alpha-2}&(2.26)\cr
&\times{\Gamma(1-s+{1\over2}|q|) \Gamma(1+\alpha-s+{1\over2}|q|)\over
\Gamma(s+{1\over2}|q|) \Gamma(s-\alpha+{1\over2}|q|)}
X_{-q}(1-s,-\alpha;\tilde{d}/c)\cr
}
$$
with $\tilde{d}d\equiv 1\bmod c$. By the convexity argument of
Phragm\'en and Lindel\"of we deduce from this that $X_q$ is of
polynomial growth with respect to all involved parameters as far as
$s$ remains in an arbitrary but fixed vertical strip.
\par
The last lemma allows us to shift the contour in $(2.24)$ to the left
as we like. The functional equation $(2.26)$ yields
\smallskip
\noindent
{\bf Lemma 2.2.}\quad{\it If
$$
1+\max(0,\Re\alpha)<\Re(\gamma+\delta),\quad
|\Re\alpha|+\Re\beta<-2,\leqno(2.27)
$$
then we have, for any non-zero $m\in\B{Z}[i]$,
$$
B_m(\alpha,\beta;g^*(\cdot\,;\gamma,\delta))=\big[B_m^{(0)}
+B_m^{(1)}\big](\alpha,\beta;g^*(\cdot\,;\gamma,\delta)).
\leqno(2.28)
$$
Here
$$
\leqalignno{
B_m^{(0)}&(\alpha,\beta;g^*(\cdot\,;\gamma,\delta))&(2.29)\cr
&=2\pi|m|^2 \sigma_{\alpha+\beta-1}(m){{\zeta_\r{F}(1-\alpha)
\zeta_\r{F}(1-\beta)}\over\zeta_\r{F}(2-\alpha-\beta)}
\tilde{g}_0(1;\gamma,\delta)\cr
&+2\pi|m|^{2(\alpha+1)}\sigma_{\beta-\alpha-1}(m)
{{\zeta_\r{F}(1+\alpha)\zeta_\r{F}(1-\beta)}
\over\zeta_\r{F}(2+\alpha-\beta)} \tilde{g}_0(1+\alpha;\gamma,\delta)
\cr
}
$$
and
$$
\leqalignno{
&B_m^{(1)}(\alpha,\beta;g^*(\cdot\,;\gamma,\delta))&(2.30)\cr
&=-{1\over4}i\pi^{2\beta-1}\zeta_\r{F}(1-\beta) |m|^{\alpha+\beta+1}
\sum_{n\ne0}\sigma_{-\alpha}(n)|n|^{
\alpha+\beta-1}S_{m,n}(\alpha,\beta,\gamma,\delta;g),\cr
}
$$
where
$$
\leqalignno{
&\qquad S_{m,n}(\alpha,\beta,\gamma,\delta;g)=\sum_{c\ne0} {1\over
{|c|^2}}S_\r{F}(m,n\,;c)\,[g]\!\Big({2\pi\over c}
\sqrt{mn};\alpha,\beta,\gamma,\delta\Big);&(2.31)\cr
&\qquad[g](u;\alpha,\beta,\gamma,\delta)
=(|u|/2)^{-2(1+\alpha+\beta)}\cr
&\qquad\times\sum_{q\in\B{Z}} (u/|u|)^{-2q}
\int_{(\eta)}{\Gamma(1-s+{1\over2}|q|)
\Gamma(1+\alpha-s+{1\over2}|q|)\over \Gamma(s+{1\over2}|q|)
\Gamma(s-\alpha+{1\over2}|q|)}
\tilde{g}_q(s;\gamma,\delta)(|u|/2)^{4s}ds&(2.32)\cr
}
$$
with $\eta<1+\min(0,\Re\alpha)$. All members on the left sides of
$(2.28)$--$(2.32)$ are regular functions of the four complex
variables in the domain $(2.27)$. }
\smallskip
\noindent
Proof. This is analogous to the rational case, which is developed in
Section 4.3 of [30] (see also [27], [28]). The first condition in
$(2.27)$ comes from $(2.20)$. The second condition there allows us to
shift the contour in $(2.24)$ to $(\xi)$ with $1+{1\over2}
\Re(\alpha+\beta)<\xi<\min(0,\Re\alpha)$. Then the right side of
$(2.29)$, which depends on $(2.18)$, is the contribution of the poles
at $s=1,\,1+\alpha$ which occur only when $q=0$. With this choice of
the contour, in place of $(\eta)$ in $(2.32)$, the absolute
convergence throughout $(2.30)$--$(2.32)$ follows solely from
$(2.12)$. This gives the regularity assertion. To finish the proof we
move the contour from $(\xi)$ to $(\eta)$ as is specified above.
\smallskip
Note that $[g](u;\alpha,\beta,\gamma,\delta)$ is even with respect to
$u$ so that the choice of square root makes no difference in
$(2.31)$. A relatively closed expression for the transform
$g\mapsto[g]$ is available, though it is not much relevant to our
present purpose; see Remark at the end of Section 12. The use of the
Weil bound $(8.14)$ for $S_\r{F}$ gives the refinement of the second
condition in $(2.27)$ to $|\Re\alpha|+\Re\beta<-{3\over2}$. It
should, however, be observed that the domain $(2.27)$ in $\B{C}^4$,
even with this refinement, does not contain the critical point
$(0,0,{1\over2},{1\over2})$ that corresponds to $\r{p}_{1\over2}$ in
$(2.2)$. In other words, the estimate $(8.14)$ of individual
Kloosterman sums does not suffice; we need a massive cancellation
among Kloosterman sums. We shall demonstrate, in the final section,
that the sum formula for $S_\r{F}$ of the Kuznetsov type serves this
purpose. It can be regarded as a device to separate the variables
$m,\, n$ in $(2.31)$. Taking the result of the separation into
$(2.30)$, a sum of products of Hecke series emerges. The fact that
these functions are {\it entire} and of polynomial growth gives rise
to the desired analytic, or more precisely, meromorphic continuation
of $(2.6)$ to $\B{C}^4$.
\smallskip
\noindent
{\csc Remark.} The right side of $(2.7)$ is called the complex binary
additive divisor sum; for its rational counterpart see [28]. The
dissection leading to $(2.6)$ is crucial in our argument. We stress
that it is not precisely the analogue of the Atkinson dissection in
the rational case (see Section 4.2 of [30]). Observe that on the
right side of $(2.6)$ three times the first term is hidden in the
second term; and thus the latter does not represent the {\it
non-diagonal} part of the sum $(2.3)$ in the traditional sense. A
direct extension of Atkinson's device might be to classify the
summands in $(2.3)$ according either to the norms of $k$ and $l$ or
to the ideals generated by them. Then it would, however, be difficult
to relate the result of dissection with any Kloosterman sums or like.
In our argument we exploit a geometric feature specific to the ring
$\B{Z}[i]$; that is, our dissection is based on the {\it lattice},
rather than arithmetic, structure of $\B{Z}[i]$. If one tries to
 consider the fourth moment of the Dedekind zeta-function of any
imaginary quadratic number field with class number larger than one,
the dissection argument will become an issue. We note that our
argument has, nevertheless, a certain generality as well; it extends
to
$$
\int_{-\infty}^\infty|\zeta_\r{F}(\txt{1\over2}+it,a)
\zeta_\r{F}(\txt{1\over2}+it,b)|^2g(t)dt\leqno(2.33)
$$
for any $a,\,b\in\B{Z}$.
\medskip
\centerline{\bf 3. The group $\r{PSL}_2(\B{C})$}
\smallskip
\noindent
The $S_{m,n}$ in $(2.31)$ is a sum of Kloosterman sums. The spectral
decomposition of it requires a considerable dose of the
representation theory of the Lie group $\r{PSL}_2(\B{C})$.
\smallskip
With this aim in mind, we shall work with functions on
$\r{SL}_2(\B{C})$ which are left-invariant over $\r{SL}_2(\B{Z}[i])$.
Since it is implied that they are even, i.e., $f(-\r{g})=f(\r{g})$,
we are actually dealing with
$$
G=\r{PSL}_2(\B{C}),\quad \varGamma=\r{PSL}_2(\B{Z}[i]).\leqno(3.1)
$$
Denoting by $\left[{a\atop c}{b\atop d}\right]$ the projective image
of the elements $\pm\left({a\atop c}{b\atop d}\right)$ of
$\r{SL}_2(\B{C})$, we put
$$
\r{n}[z]=\left[\matrix{1&z\cr&1}\right],\quad
\r{h}[u]=\left[\matrix{u&\cr&1/u}\right],\quad
\r{k}[\alpha,\beta]=\left[\matrix{\alpha&\beta\cr
-\bar\beta&\bar\alpha}\right]\leqno(3.2)
$$
for $z,u,\alpha,\beta\in\B{C}$, $u\not=0$, $|\alpha|^2+|\beta|^2=1$;
and also
$$
N=\{\r{n}[z]: z\in\B{C}\},\quad A=\{\r{a}[r]: r>0\},\quad
K=\r{PSU}(2)=\{\r{k}[\alpha,\beta]: \alpha,\beta\in\B{C}\}
\leqno(3.3)
$$
with $\r{a}[r]=\r{h}[\sqrt{r}]$. We have, for the Euler angles $
\varphi,\theta,\psi\in\B{R}$,
$$
\r{k}[\alpha,\beta]=\r{h}[e^{i\varphi/2}]
\r{k}[\cos(\txt{1\over2}\theta),
i\sin(\txt{1\over2}\theta)]\r{h}[e^{i\psi/2}]. \leqno(3.4)
$$
We have the Iwasawa decomposition
$$
G=NAK,\leqno(3.5)
$$
which we write, e.g.,
$\r{g}=\r{n}\r{a}\r{k}=\r{n}[z]\r{a}[r]\r{k}[\alpha,\beta]$, and
understand as a coordinate system on $G$. With it, Haar measures on
respective groups are given by
$$
d\r{n}=d_+\!z,\quad d\r{a}={1\over r}dr,\quad
d\r{k}={1\over8\pi^2}\sin\theta\,d\varphi\, d\theta\,d\psi,
\leqno(3.6)
$$
and
$$
d\r{g}={1\over r^2}d\r{n}\,d\r{a}\,d\r{k}.\leqno(3.7)
$$
We have, in particular,
$$
\int_Kd\r{k}=1,\qquad \int_{\varGamma\backslash G}d\r{g}=
{2\over\pi^2}\zeta_\r{F}(2),\leqno(3.8)
$$
where, with an obvious abuse of notation,
$$
\varGamma\backslash G=\varGamma\backslash G/K\cdot K\cong
\varGamma\backslash\B{H}^3\cdot K.\leqno(3.9)
$$
Here $\B{H}^3$ is the hyperbolic upper half-space, and
$\varGamma\backslash\B{H}^3$ is represented by the set
$$
\left\{(z,r)\in\B{H}^3:\,
|\Re{z}|\le\txt{1\over2},\,0\le\Im{z}\le\txt{1\over2},\,
|z|^2+r^2\ge1\right\},\leqno(3.10)
$$
which is the classical fundamental domain of the Picard group.
\smallskip
Next, the Lie algebra $\f{g}$ of $G$ has a basis consisting of the six
elements
$$
\leqalignno{ \b{H}_1&=\txt{1\over2}\pmatrix{1&\cr&-1\cr},\,
\b{V}_1=\txt{1\over2}\pmatrix{&1\cr1&\cr},\,
\b{W}_1=\txt{1\over2}\pmatrix{&1\cr-1&\cr},&(3.11)\cr
\b{H}_2&=\txt{1\over2}\pmatrix{i&\cr&-i\cr},\,
\b{V}_2=\txt{1\over2}\pmatrix{&i\cr -i&\cr},\,
\b{W}_2=\txt{1\over2}\pmatrix{&i\cr i&\cr}.
\cr
}
$$
These generate the universal enveloping algebra ${\cal U}(\f{g})$. We
identify them with right differentiations on $G$; that is, e.g.,
$$
(\b{H}_2f)(\r{g})=\lim_{\B{R}\ni t\to 0}{d\over
dt}f(\r{g}\exp(t\b{H}_2))=(\partial_\psi f)(\r{g}). \leqno(3.12)
$$
Then ${\cal U}(\f{g})$ is the set of all left-invariant differential
operators on $G$. The center ${\cal Z}(\f{g})$ of ${\cal U}(\f{g})$
is the polynomial ring $\B{C}[\Omega_+,\Omega_-]$ with the two
Casimir elements
$$
\Omega_\pm={1\over8}\left(({\b{H}_1\mp i\b{H}_2})^2+(\b{V}_1\mp
i\b{W}_2)^2 -(\b{W}_1\mp i\b{V}_2)^2\right), \leqno(3.13)
$$
where the factor $i$ means the complexification of respective
elements. In terms of the Iwasawa coordinates we have
$$
\leqalignno{ \Omega_+=&{1\over2}r^2\partial_z\partial_{\bar{z}}
+{1\over2}re^{i\varphi}\cot\theta\,\partial_z
\partial_{\varphi}
-{1\over2}ire^{i\varphi}\partial_z\partial_{\theta}
-{re^{i\varphi}\over2\sin\theta}\partial_z\partial_\psi
&(3.14)\cr
&+{1\over8}r^2\partial^2_r
-{1\over4}ir\partial_r\partial_\varphi
-{1\over8}\partial^2_\varphi
-{1\over8}r\partial_r
+{1\over4}i\partial_\varphi;
\cr
}
$$
and
$$
\leqalignno{ \Omega_-=&{1\over2}r^2\partial_z\partial_{\bar{z}}
+{1\over2}re^{-i\varphi}\cot\theta\,\partial_{\bar{z}}
\partial_{\varphi}
+{1\over2}ire^{-i\varphi}\partial_{\bar{z}}\partial_{\theta}
-{re^{-i\varphi}\over2\sin\theta}\partial_{\bar{z}}
\partial_\psi
&(3.15)\cr
&+{1\over8}r^2\partial^2_r
+{1\over4}ir\partial_r\partial_\varphi
-{1\over8}\partial^2_\varphi
-{1\over8}r\partial_r
-{1\over4}i\partial_\varphi.
\cr
}
$$
\smallskip
Restricting ourselves to the maximal compact subgroup $K$, we note
that its Lie algebra $\f{k}$, and thus its universal enveloping
algebra ${\cal U}(\f{k})$ are generated by $\b{H}_2$, $\b{W}_1$, and
$\b{W}_2$. The center ${\cal Z}(\f{k})$ of ${\cal U}(\f{k})$ is the
polynomial ring $\B{C}[\Omega_\f{k}]$ with
$$
\Omega_\f{k}={1\over2}(\b{H}^2_2+\b{W}_1^2+ \b{W}_2^2).\leqno(3.16)
$$
In terms of the Iwasawa coordinates we have
$$
\Omega_\f{k}={1\over2\sin^2\theta}\left(\partial^2_\varphi
+\sin^2\theta\,\partial_\theta^2
+\partial^2_\psi
-2\cos\theta\,\partial_\varphi\partial_\psi
+\sin\theta\cos\theta\,\partial_\theta\right).
\leqno(3.17)
$$
Let $L^2(K)$ be the Hilbert space of all functions on $K$ which are
square-integrable over $K$ with respect to the Haar measure $d\r{k}$.
To describe the structure of $L^2(K)$, and hence the unitary
representations of the compact group $K$, we put, for $|q|\le l$,
$$
(\alpha x-\bar{\beta})^{l-q}(\beta x+\bar{\alpha})^{l+q}
=\sum_{p=-l}^l\Phi^l_{p,q}(\r{k}[\alpha,\beta]) x^{l-p}.\leqno(3.18)
$$
We have
$$
\Phi_{-p,-q}^l=(-1)^{p+q}\overline{\Phi_{p,q}^l}\,.\leqno(3.19)
$$
Also we have, with $(3.4)$,
$$
\Phi^l_{p,q}(\r{k}[\alpha,\beta]) =e^{-ip\varphi-iq\psi}\Phi_{p,q}^l(
\r{k}[\cos(\txt{1\over2}\theta),
i\sin(\txt{1\over2}\theta)]),\leqno(3.20)
$$
from which the relations
$$
\Omega_\f{k}\Phi^l_{p,q}=-\txt{1\over2}(l^2+l)\Phi^l_{p,q}\,,\quad
\b{H}_2\Phi_{p,q}^l=-iq\Phi^l_{p,q}\, \leqno(3.21)
$$
follow with the convention $\Phi^l_{p,q}\equiv0$ if $|p|,\,|q| \le l$
is violated. The set
$$
\left\{\Phi_{p,q}^l:\,l\ge0,\,|p|,|q|\le l\right\}\leqno(3.22)
$$
is an orthogonal basis of $L^2(K)$ with norms
$$
\Vert{\Phi_{p,q}^l}\Vert_K=\Big[\int_K|\Phi_{p,q}^l(\r{k})|^2
d\r{k}\Big]^{1\over2}={1\over\sqrt{ l+{1\over2}}}{2l\choose
l-p}^{1\over2}{2l\choose l-q}^{-{1\over2}}. \leqno(3.23)
$$
The matrices $\Phi_l=(\Phi_{p,q}^l)$ realize all unitary
representations of $K$; in particular, we have
$$
\Phi_l(\r{k}_1\r{k}_2)=\Phi_l(\r{k}_1)\Phi_l(\r{k}_2),
\quad \r{k}_1,\,\r{k}_2\in K.\leqno(3.24)
$$
We arrange $(3.22)$ as
$$
L^2(K)=\overline{\bigoplus_{\scr{l,\,q}\atop \scr{|q|\le
l}}L^2(K;\,l,\,q)},
\quad L^2(K;\,l,\,q)=\bigoplus_{|p|\le l}\B{C}\Phi_{p,q}^l.
\leqno(3.25)
$$
We have
$$
L^2(K;\,l,\,q)=\left\{f\in L^2(K):\,\Omega_\f{k}
f=-\txt{1\over2}(l^2+l)f,\, \b{H}_2f=-iqf\right\}.\leqno(3.26)
$$
We call $L^2(K;\,l,\,q)$ the subspace of $L^2(K)$ of $K$-type $(l,q)$.
\smallskip
More generally, in a space in which $K$ acts, we shall say an element
has $K$-type $(l,q)$ if it is a simultaneous eigenvector of
$\Omega_\f{k}$ and $\b{H}_2$ with eigenvalues $-{1\over2}(l^2+l)$ and
$-iq$, respectively. This concept corresponds to the weight in the
theory of modular forms on the upper half plane.
\smallskip
We shall be concerned with representation spaces of $\f g$ for which
we can use the principal series representations as a model space. The
space of $K$-finite vectors in the principal series is
$$
H(\nu,p)=\{\,\hbox{\rm finite linear combinations of
$\varphi_{l,q}(\nu,p)$}\,\}\leqno(3.27)
$$
with
$$
\varphi_{l,q}(\nu,p)(\r{n}\r{a}[r]\r{k})=
r^{1+\nu}\Phi_{p,q}^l(\r{k})\quad (\nu\in\B{C}).\leqno(3.28)
$$
Formulas $(3.14)$--$(3.15)$ and $(3.20)$ imply that $H(\nu,p)$ is a
simultaneous eigenspace of $\Omega_\pm$:
$$
\Omega_\pm\varphi_{l,q}(\nu,p)={1\over8}((\nu\mp p)^2-1)
\varphi_{l,q}(\nu,p).\leqno(3.29)
$$
The space $H(\nu,p)$ is not $G$-invariant, but known to be
$\f{g}$-invariant and irreducible for values of $(\nu,p)$ that are of
interest for our purpose. Restricting functions in $H(\nu,p)$ to $K$,
we have a scalar product on $H(\nu,p)$. If $\Re\nu=0$, the group $G$
acts unitarily in the resulting Hilbert space -- unitary principal
series. In Section 8 we encounter these unitary representations as
models for irreducible subspaces of $L^2(\varGamma\backslash G)$ to
be defined there. The irreducibility can in fact be confirmed by
computing the actions of the six elements in $(3.11)$ over each
$\varphi_{l,q}(\nu,p)$, though we skip it. In what follows we shall
see that in most cases $H(\nu,p)$ is the space we are actually
dealing with via maps commuting with the action of $\f{g}$ or
equivalently of ${\cal U}(\f{g})$.
\smallskip
\noindent
{\csc Remark.} The general theory as well as specific treatments of
unitary representations of Lie groups can be found in [20], [38],
[40], and [41]. The fundamental region $(3.10)$ is due to Picard
[35], and its volume, given in $(3.8)$, to Humbert [15]. Formulas
$(3.14)$, $(3.15)$, and $(3.17)$ are obtained by first interpreting
$(3.13)$ and $(3.16)$ in terms of the local coordinates
$\r{g}=\r{n}[z_1]\r{h}[z_2]\left[{\atop1}{-1\atop}\right] \r{n}[z_3]$
with $(z_1,z_2,z_3)\in\B{C}^3$, over the big cell of the Bruhat
decomposition of $G$, and by changing variables according to the
Iwasawa coordinates. Our choice of basis elements $(3.22)$ is
somewhat different from common practice as is indicated by $(3.23)$.
This is for the sake of convenience for our later discussion. The
maps which commute with the action of $\f{g}$ are usually called
intertwining operators.
\medskip
\centerline{\bf 4. Automorphic forms}
\smallskip
\noindent
Let $C^\infty(\varGamma\backslash G)$ be the space of all smooth left
$\varGamma$-invariant or $\varGamma$-automorphic functions on $G$. We
consider subspaces composed of simultaneous eigenfunctions of
$\Omega_\pm$, $\Omega_\f{k}$, and $\b{H}_2$: Let $\chi$ be a
character on ${\cal Z}(\f{g})$. We put
$$
A_{l,q}(\chi)=\Big\{f\in C^\infty(\varGamma\backslash G):\, \Omega_\pm
f=\chi(\Omega_\pm)f,\,\hbox{\rm and of $K$-type
$(l,q)$}\Big\}.\leqno(4.1)
$$
Elements of $A_{l,q}(\chi)$ are called left $\varGamma$-automorphic
forms on $G$ of $K$-type $(l,q)$ with character $\chi$. Obviously
they are counterparts of $\r{PSL}_2(\B{Z})$-automorphic forms on
$\r{PSL}_2(\B{R})$.
\smallskip
As being eigenvalues of differential operators, $\chi(\Omega_\pm)$
cannot be arbitrary:
\smallskip
\noindent
{\bf Lemma 4.1.}\quad{\it If $A_{l,q}(\chi)\not=\{0\}$ then
$\chi=\chi_{\nu,p}$. Here $\chi_{\nu,p}$ is the character of ${\cal
Z}(\f{g})$ defined by
$$
\chi_{\nu,p}(\Omega_\pm)={1\over8}((\nu\mp p)^2-1) \leqno(4.2)
$$
with certain $\nu\in\B{C}$ and $p\in\B{Z}$, $|p|\le l$, which are
uniquely determined modulo $(\nu,p)\mapsto (-\nu,-p)$. }
\smallskip
\noindent
This assertion is a consequence of a study of Fourier coefficients of
automorphic forms, which we are going to develop. Thus, for any $f\in
C^\infty(\varGamma\backslash G)$ we have the Fourier expansion
$$
f(\r{g})= \sum_{\omega\in\B{Z}[i]} F_\omega f(\r{g}),\leqno(4.3)
$$
where
$$
F_\omega f(\r{g})= \int_{\varGamma_N\backslash
N}\psi_\omega(\r{n})^{-1}f(\r{n} \r{g}) d\r{n}\leqno(4.4)
$$
with
$$
\varGamma_N=\varGamma\cap N,\quad\psi_\omega(\r{n}[z])=\exp(2\pi
i\Re(\omega z)). \leqno(4.5)
$$
Obviously the operator $F_\omega$ commutes with every element of
${\cal U}(\f{g})$, implying that, if $f\in A_{l,q}(\chi)$, then
$F_\omega f$ is in the space
$$
\leqalignno{ W_{l,q}(\chi,\omega) =\big\{h\in C^\infty&(G):\,
h(\r{n}\r{g}) =\psi_\omega(\r{n})h(\r{g}),&(4.6)\cr
&\hbox{\rm of $K$-type $(l,q)$ with character $\chi$} \big\}.\cr
}
$$
Thus the above lemma is a corollary of
\smallskip
\noindent
{\bf Lemma 4.2.}\quad{\it If $W_{l,q}(\chi,\omega)\not=\{0\}$, then
there exist $\nu\in\B{C}$ and $p\in\B{Z}$, $|p|\le l$, such that
$\chi=\chi_{\nu,p}$. }
\smallskip
\noindent
Proof. Let $h\in W_{l,q}(\chi,\omega)$. We note that for any fixed
$\r{g}\in G$ the function $h(\r{g}\r{k})$ of $\r{k}\in K$ belongs to
$ L^2(K;\,l,\,q)$. In particular we have
$$
h(\r{g})=\sum_{|p|\le l} h_p(\r{n}\r{a})
\Phi_{p,q}^l(\r{k}).\leqno(4.7)
$$
The formulas $(3.14)$--$(3.15)$ and $(3.20)$ imply that the condition
$\Omega_\pm h=\chi(\Omega_\pm)h$ is equivalent to
$$
\leqalignno{
\qquad\chi(\Omega_+&)h_p={1\over2}(l-p)r\partial_zh_{p+1}
+{1\over8}(r^2\partial_r^2-(1+2p)r\partial_r
+4r^2\partial_z\partial_{\bar{z}}+p(p+2))h_p,&(4.8)\cr
\chi(\Omega_-&)h_p=-{1\over2}(l+p)r\partial_{\bar{z}}
h_{p-1}+{1\over8}(r^2\partial_r^2+(2p-1)r\partial_r
+4r^2\partial_z\partial_{\bar{z}}+p(p-2))h_p,&(4.9)\cr
}
$$
where it is supposed that $h_p\equiv0$ if $|p|>l$. We shall first
consider the case $\omega=0$. Then $(4.8)$--$(4.9)$ can be written as
$$
\leqalignno{
&r^2h_p''-rh_p'+(p^2+1)h_p=\txt{1\over2}(a_++a_-)h_p,&(4.10)\cr
& prh_p'-ph_p=\txt{1\over4}(a_--a_+)h_p
\cr
}
$$
with $\chi(\Omega_+)={1\over8}(a_+-1)$,
$\chi(\Omega_-)={1\over8}(a_--1)$. If $p\not=0$ the second equation
has a solution space spanned by $r^{1+\nu}$ with $\nu=(a_--a_+)/4p$;
and the first equation gives $a_\pm= (\nu\mp p)^2$. If $p=0$, we have
$a_+=a_-$. If $a_+\not=0$, then $r^{1+\nu}$ and $r^{1-\nu}$ with
$\nu^2=a_+$ span the solutions of the first equation. If $a_+=0$ then
$r$ and $r\log r$ are the corresponding solutions; and we have
$\nu=0$. This settles the case $\omega=0$. We next move to the case
$\omega\not=0$. For each $t\in\B{C}\setminus\{0\}$, let $\ell_t$ be
the left translation
$$
\ell_t f(\r{g})=f(\r{h}[t]\r{g}). \leqno(4.11)
$$
We have
$$
\ell_tW_{l,q}(\chi,\omega)=W_{l,q}(\chi,t^2\omega),\leqno(4.12)
$$
which reduces the problem to the case $\omega=1$. Any $h\in
W_{\ell,q}(\chi,1)$ has the form
$$
h(\r{n}[z]\r{a}[r]\r{k}) =\exp(\pi i(z+\bar{z})) \sum_{|m|\le
l}h_m(r)\Phi^l_{m,q}(\r{k}). \leqno(4.13)
$$
Again by $(4.8)$--$(4.9)$ we have
$$
\leqalignno{ r^2 h_m''&-(2m+1)rh_m'+ (m^2+2m-4\pi^2
r^2-8\chi(\Omega_+)) h_m &(4.14)\cr
&= - 4\pi i(l-m)rh_{m+1},\cr
r^2 h_m''& + (2m-1)rh_m'(r) +
(m^2-2m-4\pi^2r^2-8\chi(\Omega_-))h_m&(4.15)\cr
&= 4\pi i(l+m)rh_{m-1}.\cr}
$$
We write $\chi(\Omega_\pm)={1\over8}(\mu_\pm^2-1)$. We may assume,
without loss of generality, that
$$
0\le\Re\mu_+\le\Re\mu_-\,.\leqno(4.16)
$$
It is immediate that there exist constants $c_\pm,\,d_\pm$ such that
$$
h_{\pm l}(r)=c_\pm r^{l+1}K_{\mu_\pm}(2\pi r)+d_\pm
r^{l+1}I_{\mu_\pm}(2\pi r).\leqno(4.17)
$$
We consider first the case $\mu_\pm\not\in\B{Z}$. Applying inductively
the equation $(4.15)$ to $h_l(r)$ we see that in the expansion of
$h_{-l}(r)$ all terms are multiples of either $r^{\mu_+-l+1+2m}$ or
$r^{-\mu_+-l+1+2n}$ with integers $m,\,n\ge0$. On the other hand, if
$c_-\not=0$ in $(4.17)$, then $h_{-l}(r)$ has a term equal to a
multiple of $r^{-\mu_-+l+1}$. Thus we have either
$\mu_+-l+1+2m=-\mu_-+l+1$ or $-\mu_+-l+1+2n=-\mu_-+l+1$. The first
identity gives $\mu_++\mu_-=2(l-m)$, whence $0\le m\le l$ and
$\mu_+=\nu+(l-m)$, $\mu_-=-\nu+(l-m)$ with a $\nu\in\B{C}$. The
second gives $\mu_+=\mu_--2(l-n)$; that is, $\mu_+=\nu-(l-n)$,
$\mu_-=\nu+(l-n)$ with a $\nu\in\B{C}$, where we have $0\le n\le l$
because of $(4.16)$. This settles the case $c_-\not=0$. In other
case, we should have $h_{-l}(r)= d_-r^{l+1}I_{\mu_-}(2\pi r)$.
Applying inductively the equation $(4.14)$ to $h_{-l}$, we proceed in
much the same way, and obtain the assertion of the lemma. Next, if
$\mu_-\in\B{Z}$, then the above procedure yields that all terms of
$h_{l}(r)$ are multiples of either $r^{\mu_--l+1+2m}\log r$ or
$r^{\mu_--l+1+2n}$ with integers $m,n\ge0$. According as either
$c_+\not=0$ or $=0$, we have $\mu_--l+1+2m=\mu_++l+1$ or
$\mu_--l+1+2n=\mu_++l+1$, respectively. Thus we are again led to the
same conclusion. Finally, we observe that $\mu_+\in\B{Z}$ implies
$\mu_-\in\B{Z}$; and we end the proof.
\smallskip
It should be remarked that in the above it is proved that if
$(\nu,p)\not=(0,0)$, then
$$
W_{l,q}(\chi_{\nu,p},0)=\B{C}\,\varphi_{l,q}(\nu,p)\oplus
\B{C}\,\varphi_{l,q}(-\nu,-p);\leqno(4.18)
$$
otherwise
$$
W_{l,q}(\chi_{0,0},0)=\B{C}\,\varphi_{l,q}(0,p)
\oplus\B{C}\,\partial_\nu\varphi_{l,q}(\nu,p)|_{\nu=0}\,.
\leqno(4.19)
$$
Thus $\dim{W_{l,q}(\chi_{\nu,p},0)}=2$, but for $\omega\not=0$ we know
only $\dim{W_{l,q}(\chi_{\nu,p},\omega)}\le2$ at this stage.
\par
If a function $f$ on $G$ satisfies the bound
$$
f(\r{n}\r{a}[r]\r{k})=O(r^b)\leqno(4.20)
$$
as $r\uparrow \infty$ with a certain real constant $b$, then we say
that $f$ is of polynomial growth. The dependency of the bound on the
set where $\r{n}$ and $\r{k}$ move around is to be mentioned in our
discussion. Automorphic forms with polynomial growth are the most
interesting ones. The Fourier terms inherit this growth property, and
we put
$$
W_{l,q}^{\rm pol}(\chi_{\nu,p},\omega)= \Big\{h\in
W_{l,q}(\chi_{\nu,p},\omega): \hbox{\rm of polynomial growth,
uniformly over $K$}\Big\}.\leqno(4.21)
$$
In this way we get rid of the $I$-Bessel term in $(4.17)$, which is of
exponential growth. On noting basic properties of the $K$-Bessel
function, we have readily
\smallskip
\noindent
{\bf Lemma 4.3.}\quad {\it Let $\omega\not=0$. If $W_{l,q}^{\rm
pol}(\chi_{\nu,p},\omega)$ is non-zero, then it has dimension one.
Any generator $h$ satisfies
$$
h(\r{n}\r{a}[r]\r{k}) = O\left( |\omega r|^b
e^{-2\pi|\omega|r}\right),\leqno(4.22)
$$
as $r\uparrow\infty$, with a certain real $b$. }
\smallskip
\noindent
We shall prove in the next section that we have actually $\dim
W_{l,q}^{\rm pol}(\chi_{\nu,p},\omega)=1$ for any $\omega\not=0$.
Moreover, we shall later show that $\dim
W_{l,q}(\chi_{\nu,p},\omega)=2$ always (see $(6.17)$).
\smallskip
We next introduce the notion of cusp forms: Let
$$
\leqalignno{ A^{\rm pol}_{l,q}(\chi_{\nu,p})=\big\{f\in
A_{l,q}(\chi_{\nu,p}):&\,\hbox{\rm of polynomial growth,}&(4.23)\cr
&\qquad\hbox{\rm uniformly over $N$ and $K$}\big\};
\cr
}
$$
and put
$$
{}^0\!A_{l,q}(\chi_{\nu,p})=\{f\in A^{\rm pol}_{l,q}(\chi_{\nu,p}):\,
F_0f=0\}.\leqno(4.24)
$$
This is the space of cusp forms of $K$-type $(l,q)$ with character
$\chi_{\nu,p}$. The description of the Fourier terms in $(4.17)$, and
$(3.14)$--$(3.15)$ imply:
\smallskip
\noindent
{\bf Lemma 4.4.}\quad{\it All cusp forms $f$ are real-analytic and of
exponential decay:
$$
f(\r{n}\r{a}[r]\r{k})=O(e^{-\pi r})\leqno(4.25)
$$
uniformly over $N$, $K$ as $r$ tends to infinity. }
\smallskip
\noindent
{\csc Remark.} Automorphic forms can be defined on much more general
Lie groups, see, e.g., Harish Chandra's lecture notes [13].
\medskip
\centerline{\bf 5. Jacquet operator}
\smallskip
\noindent
In order to give explicitly an element in the space $W_{l,q}^{\rm
pol}(\chi_{\nu,p},\omega)$, $\omega\not=0$, we shall appeal to the
Jacquet integral. This device turns up in the computation of the
Fourier expansion of Poincar\'e series.
\smallskip
Thus, let $f$ be a function on $\varGamma_N\backslash G$, with which
we generate the Poincar\'e series
$$
P_f(\r{g})={1\over2}
\sum_{\gamma\in\varGamma_N\backslash\varGamma}f(\gamma\r{g}).
\leqno(5.1)
$$
We shall ignore the convergence issue temporarily. Via the Bruhat
decomposition we have
$$
P_f(\r{g})={1\over2}\left(f(\r{g})+f(\r{h}[i]\r{g})\right)+
{1\over4}\sum_{c\not=0}\,\sum_{\scr{d\bmod
c}\atop\scr{(d,c)=1}}\sum_{\omega}
f(\r{n}[\tilde{d}/c]\r{h}[1/c]\r{w}\r{n}[d/c+\omega]\r{g}),
\leqno(5.2)
$$
where $\r{w}=\left[{\atop1}{-1\atop}\right]$;
$d\tilde{d}\equiv1\bmod\, c$. The innermost sum is, by the Poisson
sum formula, equal to
$$
\sum_{\omega}\exp(2\pi i\Re(d\omega/c)) \int_N
\psi_\omega(\r{n})^{-1}f(\r{n}[\tilde{d}/c]\r{h}[1/c]\r{w}
\r{n}\r{g})d\r{n}\leqno(5.3)
$$
with $\psi_\omega$ as in $(4.5)$. If we suppose further that $f$ is
such that
$$
f(\r{n}\r{g})=\psi_{\omega'}(\r{n})f(\r{g})\leqno(5.4)
$$
with an $\omega'\in\B{Z}[i]$, then we have
$$
\leqalignno{
P_f(\r{g})=&{1\over2}\left(f(\r{g})+f(\r{h}[i]\r{g})\right)&(5.5)\cr
&+{1\over4}\sum_\omega\sum_{c\not=0}S_\r{F}(\omega,\omega';c)
\int_N\psi_\omega(\r{n})^{-1}f(\r{h}[1/c]\r{w} \r{n}\r{g})d\r{n}.\cr
}
$$
Hence we have
$$
F_\omega P_f={1\over2}(\delta_{\omega,\omega'}f+
\delta_{\omega,-\omega'}\ell_if)
+{1\over4}\sum_{c\not=0}S_\r{F}(\omega,\omega';c)
\e{A}_\omega\ell_{1/c}f,\leqno(5.6)
$$
where $\delta$ is the Kronecker delta, $\ell$ is as in $(4.11)$, and
$$
\e{A}_\xi f(\r{g})=\int_N\psi_\xi(\r{n})^{-1}f(\r{w}\r{n}\r{g})d\r{n}
\leqno(5.7)
$$
is the Jacquet integral. A property of $\e{A}_\xi$ is
$$
\ell_t\e{A}_\xi\ell_t=|t|^4\e{A}_{t^2\xi}\leqno(5.8)
$$
for any $t\not=0$; thus only $\e{A}_0$, $\e{A}_1$ matter actually.
Obviously $\e{A}_\xi$ commutes with any element of ${\cal U}(\f{g})$.
\smallskip
Now, let us compute $\e{A}_\omega\varphi_{l,q}(\nu,p)$, which is in
$W_{l,q}(\chi_{\nu,p},\omega)$. We remark that
$$
\r{w}\r{n}[z]\r{a}[r] =\r{n}\left[{-\bar{z}\over
r^2+|z|^2}\right]\r{a}\left[{r\over r^2+|z|^2}\right]
\r{k}\bigg[{\bar{z}\over\sqrt{r^2+|z|^2}},\,
{-r\over\sqrt{r^2+|z|^2}}\bigg];\leqno(5.9)
$$
and thus
$$
\leqalignno{
&\e{A}_\omega\varphi_{l,q}(\nu,p;\r{n}\r{a}[r]\r{k}
[\alpha,\beta])&(5.10)\cr
=&\psi_\omega(\r{n})r^{1-\nu} \int_\B{C} {e^{-2\pi i\Re(\omega
rz)}\over (1+|z|^2)^{\nu+1}}
\Phi_{p,q}^l\left(\r{k}\bigg[{\bar{z}\over \sqrt{1+|z|^2}},\,
{-1\over\sqrt{1+|z|^2}}\bigg] \r{k}[\alpha,\beta]\right)d_+\!z.
\cr
}
$$
This shows that for $\Re\nu>0$
$$
\e{A}_\omega\varphi_{l,q}(\nu,p)\in W_{l,q}^{\rm
pol}(\chi_{\nu,p}).\leqno(5.11)
$$
We have, by $(3.24)$,
$$
\e{A}_\omega\varphi_{l,q}(\nu,p;\r{n}\r{a}[r]\r{k})
=\psi_\omega(\r{n})\sum_{|m|\le
l}v_m^l(r)\Phi_{m,q}^l(\r{k}),\leqno(5.12)
$$
where
$$
v_m^l(r) =r^{1-\nu} \int_\B{C} {e^{-2\pi i\Re(\omega rz)}\over
(1+|z|^2)^{\nu+1}} \Phi_{p,m}^l\left(\r{k}\bigg[{\bar{z}\over
\sqrt{1+|z|^2}},\, {-1\over\sqrt{1+|z|^2}}\bigg]\right)d_+\!z.
\leqno(5.13)
$$
The relation
$$
\r{k}[e^{-i\phi}\alpha,\beta]=\r{h}[e^{-i\phi/2}]
\r{k}[\alpha,\beta]\r{h}[e^{-i\phi/2}]\leqno(5.14)
$$
and $(3.20)$ imply that, after the change of variables $z=ue^{i\phi}$,
the last integral becomes
$$
\leqalignno{ \int_0^\infty{u\over
(1+u^2)^{\nu+1}}&\Phi_{p,m}^l\left(\r{k}\bigg[{u\over
\sqrt{1+u^2}},\, {-1\over\sqrt{1+u^2}}\bigg]\right)&(5.15)\cr
&\times\int_{-\pi}^\pi\exp\left((p+m)i\phi-2\pi i\Re(\omega
rue^{i\phi})\right)d\phi\,du. }
$$
Thus we see that if $\omega=0$, then
$$
v_m^l(r) =2\pi\delta_{m,-p}r^{1-\nu} \int_0^\infty {u\over
(1+u^2)^{\nu+1}} \Phi_{p,-p}^l\left(\r{k}\bigg[{u\over
\sqrt{1+u^2}},\, {-1\over\sqrt{1+u^2}}\bigg]\right)du\,;\leqno(5.16)
$$
and if $\omega\not=0$, then
$$
\leqalignno{ v_m^l(r)
&=2\pi r^{1-\nu}(i\omega/|\omega|)^{-p-m}
&(5.17)\cr
&\times\int_0^\infty u{J_{p+m}(2\pi |\omega| ru)\over (1+u^2)^{\nu+1}}
\Phi_{p,m}^l\left(\r{k}\bigg[{u\over \sqrt{1+u^2}},\,
{-1\over\sqrt{1+u^2}}\bigg]\right)du.
\cr
}
$$
The integral in $(5.16)$ is, by the definition of $\Phi_{p,-p}^l$,
equal to
$$
\leqalignno{
&(-1)^{l-p}\sum_{a=0}^{l-|p|}(-1)^a{l+|p|\choose a} {l-|p|\choose
a}\int_0^\infty{u^{2a+1} \over(1+u^2)^{\nu+l+1}}du&(5.18)\cr
=&{1\over2}{\Gamma(l+|p|+1)\Gamma(\nu+|p|)
\over\Gamma(\nu+l+1)\Gamma(2|p|+1)}
\sum_{a=0}^{l-|p|}(-1)^a{l-|p|\choose a}{(\nu+|p|)_a
\over(2|p|+1)_a}\cr
=&{1\over2}{\Gamma(l+1-\nu) \over\Gamma(l+1+\nu)}{\Gamma(|p|+\nu)
\over\Gamma(|p|+1-\nu)}.
\cr
}
$$
The last line depends on the identity
$$
\sum_{j=0}^k(-1)^j{k\choose j}{(\alpha)_j\over(\beta)_j}
={(\beta-\alpha)_k\over(\beta)_k},
\quad(\alpha)_j=\alpha(\alpha+1)\cdots(\alpha+j-1),\leqno(5.19)
$$
which can be shown by induction. On the other hand, we observe that if
$p+m<0$ in the integrand in $(5.17)$ then we may replace $(p,\,m)$ by
$(-p,\,-m)$ without affecting the value of the integral, since we
have $(3.19)$ and $J_{-a}=(-1)^aJ_a$ for $a\in\B{Z}$. Thus, according
as $\r{sgn}(p+m)=\pm$, the integral is equal to
$$
\leqalignno{
&(-1)^{l-p}\sum_{a=0}^{\min\{l\mp m,l\mp p\}} (-1)^a{l\mp m\choose
a}{l\pm m\choose l\mp p-a}&(5.20)\cr
&\hskip 3cm\times\int_0^\infty
{u^{|m+p|+1+2a}J_{|m+p|}(2\pi|\omega|ru) \over(1+u^2)^{\nu+l+1}}du
\cr
=&(-1)^{l-p}\sum_{a=0}^{\min\{l\mp m,l\mp p\}} \sum_{b=0}^a
(-1)^{b}{l\mp m\choose a}{l\pm m\choose l\mp p-a} {a\choose b}\cr
&\hskip 3cm\times\int_0^\infty {u^{|m+p|+1}J_{|m+p|}(2\pi|\omega|ru)
\over(1+u^2)^{\nu+l+1-b}}du.\cr
}
$$
Exchanging the order of summation, the sum over $a$ taken inside is
equal to
$$
\leqalignno{&\xi_p^l(m,b) ={b!(2l-b)!\over(l-p)!(l+p)!}&(5.21)\cr
&\times{l-{1\over2}(|m+p|+|m-p|)\choose b}
{l-{1\over2}(|m+p|-|m-p|)\choose b}.
\cr
}
$$
In fact it is
$$
\leqalignno{
&\sum_{a=b}^{\min\{l\mp m,l\mp p\}} {l\mp m\choose a}{l\pm m\choose
l\mp p-a} {a\choose b}&(5.22)\cr
&={(l-m)!(l+m)!\over b!}\sum_{a=b}^{\min\{A,B\}}{1\over
(A-a)!(B-a)!(a-b)!(a+c)!} }
$$
with $A=l\mp m$, $B=l\mp p$, $c=|m+p|$; and on the assumption $A\le B$
$$
\leqalignno{ \sum_{a=b}^{\min\{A,B\}}&={1\over(A-b)!(B-b)!(c+b)!}
\sum_{d=0}^{A-b}(-1)^d{A-b\choose d}{(b-B)_d\over(c+b+1)_d}&(5.23)\cr
&={(A+B+c-b)!\over(A-b)!(B-b)!(A+c)!(B+c)!} }
$$
because of $(5.19)$. Hence we get $(5.21)$.

 Collecting these and invoking the formula
$$
\int_0^\infty{u^{\tau+1}J_\tau(tu)\over(1+u^2)^{\eta+1}}du
={(t/2)^\eta\over\Gamma(\eta+1)} K_{\tau-\eta}(t)\quad
(\,t>0\,),\leqno(5.24)
$$
which holds for $-1<\Re\tau<2\Re\eta+{3\over2}$ (formula (2) on p.\
434 of [42]), we obtain
\smallskip
\noindent
{\bf Lemma 5.1.}\quad{\it We have that for $\omega=0$
$$
\e{A}_0\varphi_{l,q}(\nu,p)=\pi{\Gamma(l+1-\nu)\over
\Gamma(l+1+\nu)}{\Gamma(|p|+\nu)\over\Gamma(|p|+1-\nu)}
\varphi_{l,q}(-\nu,-p);\leqno(5.25)
$$
and for $\omega\not=0$
$$
\leqalignno{
&\e{A}_\omega \varphi_{l,q}(\nu,p)(\r{n}\r{a}[r]\r{k})&(5.26)\cr
&=2(-1)^{l-p}\pi^\nu|\omega|^{\nu-1}\psi_\omega(\r{n}) \sum_{|m|\le l}
(i\omega/|\omega|)^{-m-p}\alpha_m^l
(\nu,p;|\omega|r)\Phi_{m,q}^l(\r{k}) }
$$
with
$$
\leqalignno{
&\alpha_m^l(\nu,p;r)&(5.27)\cr &= \sum_{j=0}^{l-|m+p|/2-|m-p|/2}
(-1)^j\xi_p^l(m,j) {(\pi r)^{l+1-j}\over\Gamma(\nu+l+1-j)}
K_{\nu+l-|m+p|-j}(2\pi r).\cr
}
$$
}
\par
We see that with respect to $\nu$ the function
$\e{A}_0\varphi_{l,q}(\nu,p)$ is meromorphic, and for $\omega\not=0$
the function $\e{A}_\omega\varphi_{l,q}(\nu,p)$ is entire. Thus,
taking into account analytic continuation, we may extend
$\e{A}_\omega$ so that $\e{A}_\omega\varphi_{l,q}(\nu,p)$ is given by
the right side members of $(5.25)$--$(5.26)$, as far as they are
regular. In this way we define the {\it Jacquet operator\/}:
$$
\e{A}_\omega:\, H(\nu,p)\to W^{\rm
pol}(\chi_{\nu,p},\omega)\leqno(5.28)
$$
where the right side is the space spanned by all $W_{l,q}^{\rm
pol}(\chi_{\nu,p},\omega)$, $|p|\le l$, $|q|\le l$. The function
$\e{A}_\omega\varphi_{l,q}(\nu,p)$ spans the space $W_{l,q}^{\rm
pol}(\chi_{\nu,p},\omega)$, $\omega\not=0$, for all values of
$\nu\in\B{C}$, $p\in\B{Z}$. In particular, since the space
$W_{l,q}^{\rm pol}(\chi_{-\nu,-p},\omega)$ is identical to
$W_{l,q}^{\rm pol}(\chi_{\nu,p},\omega)$, the function
$\e{A}_\omega\varphi_{l,q}(-\nu,-p)$ is a multiple of
$\e{A}_\omega\varphi_{l,q}(\nu,p)$. Checking the coefficients of
$\Phi_{l,q}^l(\r{k})$ in these functions we find the functional
equation
$$
\leqalignno{ (\pi|\omega|)^{-\nu}
(i\omega/|\omega|)^p&\Gamma(l+1+\nu)\e{A}_\omega
\varphi(\nu,p)&(5.29)\cr
&=(\pi|\omega|)^{\nu}
(i\omega/|\omega|)^{-p}\Gamma(l+1-\nu)\e{A}_\omega
\varphi(-\nu,-p).\cr
}
$$
\smallskip
Note that the term Jacquet operator is limited to its application to
the space $H(\nu,p)$, whereas we use the term Jacquet integral
wherever it applies. This abuse of terminology should not cause
confusion in our later discussion.
\smallskip
Now, the most important example of automorphic forms that are not
cuspidal but of polynomial growth is offered by the Eisenstein series
of $K$-type $(l,q)$:
$$
e_{l,q}(\nu,p;\r{g})={1\over2}\sum_{\gamma\in
\varGamma_N\backslash\varGamma}\varphi_{l,q}
(\nu,p)(\gamma\r{g})\qquad (\,\Re\nu>1\,)\leqno(5.30)
$$
with $p\in2\B{Z}$. We need this condition on $p$ to have a non-trivial
sum; note that $(3.20)$ implies
$\varphi_{l,q}\left(\nu,p;\r{h}[i]\r{g}\right)=(-1)^p
\varphi_{l,q}(\nu,p;\r{g})$. The series converges absolutely in the
indicated domain of $\nu$, and is regular there, which is the same as
in the $K$-trivial case (see Section 3.2 of [8]). The Fourier
expansion of $e_{l,q}(\nu,p)$ can be obtained as an application of
the foregoing discussion. Obviously we have
$$
\e{A}_\omega\ell_{1/c}\varphi_{l,q}(\nu,p)(\r{a}[r]\r{k})
=(c/|c|)^{2p}|c|^{-2(\nu+1)}
\e{A}_\omega\varphi_{l,q}(\nu,p)(\r{a}[r]\r{k}).\leqno(5.31)
$$
Thus, by $(2.18)$, we have
\smallskip
\noindent
{\bf Lemma 5.2.}\quad{\it The Eisenstein series $e_{l,q}(\nu,p)$,
$p\in2\B{Z}$, is meromorphic over $\B{C}$ with respect to $\nu$. When
it is regular, we have the Fourier expansion
$$
\leqalignno{ e_{l,q}(\nu,p)&=\varphi_{l,q}(\nu,p)+ \pi^{2\nu}
{\Gamma(l+1-\nu)\over\Gamma(l+1+\nu)}
{\zeta_\r{F}(1-\nu,p/2)\over\zeta_\r{F}(1+\nu,p/2)}
\varphi_{l,q}(-\nu,-p)&(5.32)\cr
&+{1\over\zeta_\r{F}(1+\nu,p/2)}\sum_{\omega\not=0}
\sigma_{-\nu}(\omega,p/2) \e{A}_\omega\varphi_{l,q}(\nu,p).
\cr}
$$
We also have the functional equation}
$$
\leqalignno{ \pi^{-\nu}\Gamma(l+1+\nu)&\zeta_\r{F}(1+\nu,p/2)
e_{l,q}(\nu,p)&(5.33)\cr
&=\pi^\nu\Gamma(l+1-\nu)\zeta_\r{F}(1-\nu,p/2) e_{l,q}(-\nu,-p).
\cr}
$$
Proof. These assertions are consequences of the previous lemma, the
identity $(5.29)$, and the functional equation
$$
\pi^{-\nu}\Gamma(|p|+\nu)\zeta_\r{F}(\nu,p/2)
=\pi^{\nu-1}\Gamma(|p|+1-\nu)\zeta_\r{F}(1-\nu,p/2).\leqno(5.34)
$$
This ends the proof.
\smallskip
Note that in the present arithmetical situation we do not need to
establish Langlands' analytic continuation [23] of Eisenstein series.
We stress also that the above discussion implies that each cusp-form
$\psi\in{}^0\!A_{l,q}(\chi_{\nu,p})$ has the Fourier expansion
$$
\hbox{ $\displaystyle\psi=\sum_{\omega\ne0}c(\omega)
\e{A}_\omega\varphi_{l,q}(\nu,p)\quad$ or $\quad F_\omega\psi=
c(\omega)\e{A}_\omega\varphi_{l,q}(\nu,p)$}\leqno(5.35)
$$
with certain complex numbers $c(\omega)$. Because of this, instead of
considering individual automorphic forms, we study systems that
behave under the action of $\f{g}$ in the same way as the
$\varphi_{l,q}(\nu,p)$. Thus, automorphic representations move to the
focus of interest; that are linear maps from the model space
$H(\nu,p)$ to $C^\infty(\varGamma\backslash G)$ that commute with the
action of $\f{g}$. Specifically we have, for any $\b{X}\in\f{g}$,
$$
F_\omega \b{X}\psi=c(\omega)\e{A}_\omega\b{X}
\varphi_{l,q}(\nu,p).\leqno(5.36)
$$
This means that the result of a right differentiation applied to a
cusp-form is a sum of a finite linear combination of cusp-forms,
since $H(\nu,p)$ is $\f{g}$-invariant. Moreover, we see that the set
of automorphic functions $\f{g}\psi$ or ${\cal U}(\f{g})\psi$ share
the Fourier coefficients $\{c(\omega)\}$ in the sense expressed by
the identity $(5.36)$. In passing, we note that in $(5.35)$ we have
$$
c(-\omega)=(-1)^pc(\omega).\leqno(5.37)
$$
This is because $\r{h}[i]\in\varGamma$ and
$\ell_i\e{A}_\omega=\e{A}_{-\omega}\ell_i^{-1}$.
\smallskip
\noindent
{\csc Remark.} The operator $\e{A}_\omega$ has been studied by Jacquet
[17] for more general groups than $\r{PSL}_2(\B{C})$. For
$\r{PSL}_2(\B{R})$, one obtains an expression in terms of Whittaker
functions. Basic properties, like $(5.34)$, of Hecke $L$-functions
associated with Gr\"ossencharakters can be found in [14].
\medskip
\centerline{\bf 6. Goodman--Wallach operator}
\smallskip
\noindent
The Jacquet integral has given a solution to the system
$(4.14)$--$(4.15)$, which is at most of polynomial growth in the
sense of $(4.20)$; and it has fixed the operator $\e{A}_\omega$. The
formula $(4.17)$ suggests, however, the existence of a solution of
exponential growth. To construct such a solution we shall employ a
method due to Goodman and Wallach [12]. We shall have a map
$$
\e{B}_\omega:\, H(\nu,p)\to W(\chi_{\nu,p},\omega),\leqno(6.1)
$$
where the right side is spanned by all $W_{l,q}(\chi_{\nu,p},\omega)$,
$|p|\le l$, $|q|\le l$.
\smallskip
Thus, let $\varphi\in H(\nu,p)$ be arbitrary. We shall find a vector
$\{a(m,n):m,n\ge 0\}$, which depends only on $\nu,p,\omega$, so that
$$
\e{B}_\omega\varphi(\r{g})=
\sum_{m,n\ge0}a(m,n)\partial_z^m\,\partial_{\bar{z}}^n
\,\varphi(\r{w}\r{n}[z]\r{w}^{-1}\r{g})\vert_{z=0}\leqno(6.2)
$$
satisfies
$$
\e{B}_\omega\varphi(\r{n}\r{g})
=\psi_\omega(\r{n})\e{B}_\omega\varphi(\r{g}),\leqno(6.3)
$$
or
$$
\partial_t\e{B}_\omega\varphi(\r{n}[t]\r{g})\vert_{t=0}=\pi
i\omega\e{B}_\omega\varphi(\r{g}),\quad
\partial_{\bar{t}}\e{B}_\omega\varphi(\r{n}[t]\r{g})
\vert_{t=0}=\pi i\bar{\omega}\e{B}_\omega\varphi(\r{g}).\leqno(6.4)
$$
We note that
$$
\r{w}\r{n}[z]\r{w}^{-1}\pmatrix{0&t\cr0&0\cr}
(\r{w}\r{n}[z]\r{w}^{-1})^{-1}=t\pmatrix{0&1\cr0&0\cr}
+tz\pmatrix{1&0\cr0&-1\cr}
-tz^2\pmatrix{0&0\cr1&0\cr}.\leqno(6.5)
$$
Considering the exponential of the right side in a vicinity of $t=0$,
we have
$$
\leqalignno{
\partial_t\varphi&(\r{w}\r{n}[z]\r{w}^{-1}\r{n}[t]\r{g})|_{t=0}
&(6.6)\cr
&=\partial_t\varphi(\r{h}[e^{zt}]\r{w}\r{n}[z]
\r{w}^{-1}\r{g})|_{t=0}
+\partial_t\varphi(\r{w}\r{n}[z^2t+z]\r{w}^{-1}
\r{g})|_{t=0}\cr
&=(1+\nu-p)z\varphi(\r{w}\r{n}[z]\r{w}^{-1}\r{g})
+z^2\partial_z\varphi(
\r{w}\r{n}[z]\r{w}^{-1}\r{g}),
\cr
}
$$
since for $\xi=\xi_1+i\xi_2$ $(\xi_1,\,\xi_2\in\B{R})$
$$
\varphi(\r{h}[e^\xi]\r{n}\r{a}[r]\r{k})
=\varphi(\r{a}[e^{2\xi_1}r]\r{h}[e^{i\xi_2}]\r{k}) =e^{2(1+\nu)\xi_1
-2p\xi_2i}\varphi(\r{a}[r]\r{k}).\leqno(6.7)
$$
The formula $(6.6)$ gives
$$
\leqalignno{
\partial_t\,\partial_z^m\,\partial_{\bar{z}}^n
&\varphi(\r{w}\r{n}[z]\r{w}^{-1}\r{n}[t]\r{g})|_{(z,t)=(0,0)}&(6.8)\cr
&=m(\nu-p+m)\partial_z^{m-1}\,\partial_{\bar{z}}^n
\varphi(\r{w}\r{n}[z]\r{w}^{-1}\r{g})|_{z=0}\,.
\cr
}
$$
In just the same way one may show that
$$
\leqalignno{
\partial_{\bar{t}}\,\partial_z^m\,\partial_{\bar{z}}^n
&\varphi(\r{w}\r{n}[z]\r{w}^{-1}\r{n}[t]\r{g})|_{(z,t)=(0,0)}&(6.9)\cr
&=n(\nu+p+n)\partial_z^m\,\partial_{\bar{z}}^{n-1}
\varphi(\r{w}\r{n}[z]\r{w}^{-1}\r{g})|_{z=0}\,.
\cr
}
$$
{}From these and $(6.4)$ we see that the coefficients $a(m,n)$ should
satisfy the recurrence relation
$$
\leqalignno{ \pi i\omega\, a(m,n)&=(m+1)(\nu-p+m+1)a(m+1,n),&(6.10)\cr
\pi i\bar{\omega}\, a(m,n)&=(n+1)(\nu+p+n+1)a(m,n+1).
\cr
}
$$
We set the side condition
$$
a(0,0)=\{\Gamma(\nu+1+p)\Gamma(\nu+1-p)\}^{-1}.\leqno(6.11)
$$
Then we are led to
$$
a(m,n)= {(\pi i\omega)^m(\pi i\bar{\omega})^n \over
m!n!\Gamma(\nu+1-p+m)\Gamma(\nu+1+p+n)}\,.\leqno(6.12)
$$
With this choice of the vector, the sum $(6.2)$ converges absolutely
for any element $\varphi\in H(\nu,p)$. Indeed, the analyticity of
$z\mapsto \varphi(\r{w}\r{n}[z]\r{w}^{-1}g)$ provides us with a
necessary bound of the derivatives. We stress that the sum is entire
with respect to $\nu$.
\smallskip
Obviously the operator $\e{B}_\omega$ commutes with all elements of
${\cal U}(\f{g})$; and it maps $\varphi_{l,q}(\nu,p)$ into
$W_{l,q}(\chi_{\nu,p},\omega)$. Thus there should be an expansion of
$\e{B}_\omega\varphi_{l,q}(\nu,p)$ in terms of $\Phi_{m,q}^l$,
$|m|\le l\/$:
\smallskip
\noindent
{\bf Lemma 6.1.}\quad{\it We have, for any $\omega\not=0$,
$$
\leqalignno{ \e{B}_\omega&
\varphi_{l,q}(\nu,p)(\r{n}\r{a}[r]\r{k})&(6.13)\cr
&= (\pi|\omega|)^{-\nu-1}\psi_\omega(\r{n})\sum_{|m|\le
l}(-i\omega/|\omega|)^{p-m}\beta_m^l(\nu,p;|\omega|r)
\Phi_{m,q}^l(\r{k}), }
$$
where
$$
\beta_m^l(\nu,p;r)= \sum_{j=0}^{l-|m+p|/2-|m-p|/2} \xi_p^l(m,j){(\pi
r)^{l+1-j}\over\Gamma(\nu+l+1-j)} I_{\nu+l-|m+p|-j}(2\pi
r).\leqno(6.14)
$$
We have also
$$
\leqalignno{ \pi^{-2}(\pi|\omega|)^{-\nu}& (-i\omega/|\omega|)^p
\Gamma(l+1+\nu)\e{A}_\omega\varphi_{l,q}(\nu,p)&(6.15)\cr
=&-{(\pi|\omega|)^\nu\over\sin\pi\nu} (i\omega/|\omega|)^{-p}
\Gamma(l+1+\nu)\e{B}_\omega\varphi_{l,q}(\nu,p)\cr
&+{(\pi|\omega|)^{-\nu}\over\sin\pi\nu} (i\omega/|\omega|)^p
\Gamma(l+1-\nu)\e{B}_\omega\varphi_{l,q}(-\nu,-p)\,, }
$$
which is a refinement of $(5.29)$.}
\smallskip
\noindent
Proof. Let us suppose that $\nu\not\in\B{Z}$. On the right side of
$(5.27)$ we replace the $K$-Bessel function by its defining
expression:
$$
K_\xi(u)={\pi\over2\sin\pi\xi}(I_{-\xi}(u)-I_\xi(u)). \leqno(6.16)
$$
Then the function $\alpha_m^l(\nu,p;r)$ is a difference of two parts;
one is $r^\nu$ times a power series in $r$, and the other $r^{-\nu}$
times another power series. Taking these into the system
$(4.14)$--$(4.15)$, we see that each part satisfies the system. The
first part is equal to a multiple of $\beta_m^l(\nu,p;r)$, whence the
right side of $(6.13)$ belongs to $W_{l,q}(\chi_{\nu,p},\omega)$. The
other part yields another member of $W_{l,q}(\chi_{\nu,p},\omega)$;
and these two are linearly independent. Since we have shown $\dim
W_{l,q}(\chi_{\nu,p},\omega)\le2$ already, we find that
$$
\dim W_{l,q}(\chi_{\nu,p},\omega)=2\leqno(6.17)
$$
under the present specification. On the other hand, it is easy to see
that there is a power series $P$ such that $\e{B}_1
\varphi_{l,q}(\nu,p;\r{a}[r])=a(0,0)r^{1+\nu}P(r)$ with $P(0)=1$.
Hence $\e{B}_1\varphi_{l,q}(\nu,p)$ should be a constant multiple of
the right side of $(6.13)$ with $\omega=1$. The constant is equal to
$1$, as can be seen by checking the term with $m=p$. Observing that
$$
\ell_t\e{B}_1\ell_t^{-1}=\e{B}_{t^2}\leqno(6.18)
$$
because of $\r{h}[t^{-1}]\r{w}
\r{n}[z]\r{w}^{-1}\r{h}[t]=\r{w}\r{n}[t^2z]\r{w}^{-1}$, we get
$(6.13)$ for general non-zero $\omega$. As to $(6.15)$ we note that
$\e{B}_\omega\varphi_{l,q}(\nu,p)$ and
$\e{B}_\omega\varphi_{l,q}(-\nu,-p)$ are linearly independent
elements of $W_{l,q}^{\rm pol}(\chi_{\nu,p},\omega)$; and thus
$\e{A}_\omega\varphi_{l,q}(\nu,p)$ is a linear combination of them.
Computing the coefficients of $\Phi_{l,q}^l(\r{k})$ in these three
elements we obtain $(6.15)$. The case $\nu\in\B{Z}$ is settled with
analytic continuation, since both sides of $(6.13)$ are entire in
$\nu$; and $(6.15)$ is similar. This ends the proof.
\smallskip
Now, we shall show that operators $\e{A}_{\omega_1}$ and
$\e{B}_{\omega_2}$ are related in a way which will turn out to be
important in our later discussions of the sum formula for Kloosterman
sums. We observe that by $(6.13)$ we have $\e{B}_{\omega_2}
\varphi_{l,q}(\nu,p;\r{n}\r{a}[r]\r{k}) =O(r^{\Re\nu+1})$ as
$r\downarrow0$ for any $\omega_2\not=0$. Hence the Jacquet integral
$\e{A}_{\omega_1}\e{B}_{\omega_2}
\varphi_{l,q}(\nu,p;\r{n}\r{a}[r]\r{k})$ converges for $\Re\nu>0$:
\smallskip
\noindent
{\bf Lemma 6.2.}\quad{\it Let $\omega_2\not=0$, $\Re\nu>0$. Then we
have that
$$
\e{A}_0\e{B}_{\omega_2}\varphi_{l,q}(\nu,p)=(-1)^p {\sin\pi\nu\over
\nu^2-p^2}{\Gamma(l+1-\nu)\over\Gamma(l+1+\nu)}
\varphi_{l,q}(-\nu,-p)\,;\leqno(6.19)
$$
and for $\omega_1\not=0$
$$
\e{A}_{\omega_1}\e{B}_{\omega_2}\varphi_{l,q}(\nu,p)=
(\pi^2|\omega_1\omega_2|)^{-\nu}
(\omega_1\omega_2/|\omega_1\omega_2|)^p\e{J}_{\nu,p}(
2\pi\sqrt{\omega_1\omega_2}\,)\e{A}_{\omega_1}
\varphi_{l,q}(\nu,p)\leqno(6.20)
$$
with
$$
\e{J}_{\nu,p}(u)=|u/2|^{2\nu} (u/|u|)^{-2p} J^*_{\nu-p}(u)
J^*_{\nu+p}(\bar{u}).\leqno(6.21)
$$
Here $J^*_\nu(x)$ is the entire function of $x$ which is equal to
$J_\nu(x)(x/2)^{-\nu}$ for $x>0$. }
\smallskip
\noindent
Proof. Let $\varphi\in H(\nu,p)$ with $\Re\nu>0$. We may take the
integral defining $\e{A}_{\omega_1}\e{B}_{\omega_2}\varphi$ inside
the sum for $\e{B}_{\omega_2}$. The $(m,n)$-th term is equal to
$$
\leqalignno{ a(m,&n)\int_\B{C}\exp(-2\pi i\Re(\omega_1
z_1))\partial_z^m\partial_{\bar{z}}^n
\,\varphi(\r{w}\r{n}[z+z_1]\r{g})|_{z=0}\,d_+\!z_1&(6.22)\cr
= a(m,&n)\int_\B{C}\exp(-2\pi i\Re(\omega_1
z))\partial_z^m\partial_{\bar{z}}^n
\,\varphi(\r{w}\r{n}[z]\r{g})\,d_+\!z\cr
= a(m,&n)(\pi i \omega_1)^m(\pi i
{\bar\omega}_1)^n\e{A}_{\omega_1}\varphi(\r{g}).\cr
}
$$
This and $(5.25)$, $(6.12)$ readily give $(6.19)$--$(6.20)$.
\smallskip
\noindent
{\csc Remark.} The operator $\e{B}_\omega$ is due to Goodman and
Wallach [12], but for a more general context than $\r{PSL}_2(\B{C})$.
Miatello and Wallach use it to express Fourier coefficients of
Poincar\'e series in terms of their $\tau$-function, which coincides
with our $\e{J}_{\nu,p}$ if specialized to $\r{PSL}_2(\B{C})$; see
Proposition 2.7 in [25]. An extension of [12] is given in [24].
\medskip
\centerline{\bf 7. Lebedev transform}
\smallskip
\noindent
The Lebedev or the $K$-Bessel transform
$$
f\mapsto\int_0^\infty f(r)K_\nu(r){dr\over r}\leqno(7.1)
$$
plays a significant r\^ole in the theory of sum formulas for rational
Kloosterman sums. Since the function $K_\nu$ appears in the Fourier
expansion of the classical Eisenstein series over $\B{H}^2$, it
appears natural to anticipate that the corresponding function, i.e.,
$\e{A}_\omega\varphi_{l,q}(\nu,p)$, in the Fourier expansion of
$e_{l,q}(\nu,p)$ should work analogously in the present context. We
shall show in later sections that this is indeed the case. Here we
shall carry out some preparatory work. In particular we shall prove
an extension of the one-sided inversion of the Lebedev transform:
$$
\eta(\nu)={i\over\pi^2}\int_0^\infty K_\nu(r)r^{-1}
\int_{(0)}\eta(\xi)K_\xi(r)\xi\sin(\pi\xi) d\xi\,dr,\leqno(7.2)
$$
where $\eta$ is to satisfy an appropriate regularity and decay
condition.
\smallskip
Thus, let $\omega\not=0$ and put
$$
\leqalignno{P_{l,q}(N\backslash G,\omega)=\Big\{&f\in C^\infty(G):
f(\r{n}\r{g})=\psi_\omega(\r{n})f(\r{g}), \hbox{ of $K$-type
$(l,q)$},&(7.3)\cr
&\hbox{$f(\r{a}[r]\r{k})=O(r^{1+\sigma_0})$ as $r\downarrow0$,
$=O(r^{1-\sigma_\infty})$ as $r\uparrow\infty$}\Big\},
\cr
}
$$
where constants $\sigma_0,\,\sigma_\infty>0$ may depend on each $f$.
We define an extension $\e{L}_{l,q}^\omega$ of the Lebedev transform
applied to an $f\in P_{l,q}(N\backslash G,\omega)$ by
$$
\e{L}^\omega_{l,q} f(\nu,p)=\Gamma(l+1-\nu)
{(\pi|\omega|)^\nu(-i\omega/|\omega|)^{-p}
\over\pi^2\Vert\Phi_{p,q}^l\Vert_K} \int_{N\backslash G}f(\r{g})
\overline{\e{A}_{\omega}\varphi_{l,q}(-\bar{\nu},p)
(\r{g})}\,d\dot\r{g},\leqno(7.4)
$$
where $d\dot\r{g}=r^{-3}dr\,d\r{k}$ for $\r{g}=\r{a}[r]\r{k} \in
N\backslash G$. {}From $(5.27)$ it follows that the integral
converges absolutely for $|\Re\nu|<\sigma_0$, and that
$\e{L}^\omega_{l,q}f(\nu,p)$ is a regular function there. If
$\Re\nu<0$, then we have an integral representation for
$\e{A}_{\omega}\varphi_{l,q}(-\bar{\nu},p)(\r{g})$, which being
inserted into $(7.4)$ yields an absolutely convergent double integral
over $N\backslash G\times N$. Hence for $-\sigma_0<\Re\nu<0$ the last
integral is equal to
$$
\leqalignno{ \int_G
f(\r{g})\overline{\varphi_{l,q}(-\bar{\nu},p;\r{w}\r{g})}d\r{g}
&=\int_{N\backslash G}\int_N
f(\r{w}\r{n}\r{g})\overline{\varphi_{l,q}(-\bar{\nu},p;\r{g})}
d\r{n}\, d\dot\r{g}&(7.5)\cr
&=\int_{N\backslash
G}\e{A}_0f(\r{g})\overline{\varphi_{l,q}(-\bar{\nu},p;\r{g})}
d\dot\r{g}.\cr
}
$$
We may write
$$
\e{A}_0f(\r{n}\r{a}[r]\r{k})= \sum_{|m|\le
l}u_m(r)\Phi_{m,q}^l(\r{k}).\leqno(7.6)
$$
Then, we find that for $-\min(\sigma_0,\sigma_\infty)<\Re\nu<\sigma_0$
$$
\e{L}^\omega_{l,q} f(\nu,p) =\pi^{-2}(\pi|\omega|)^\nu
(-i\omega/|\omega|)^{-p} \Vert{\Phi_{p,q}^l}\Vert_K\Gamma(l+1-\nu)
\int_0^\infty u_p(r)r^{-\nu-2}dr.\leqno(7.7)
$$
Using this relation we shall show that there exists a one-sided
inversion of $\e{L}^\omega_{l,q}$:
\smallskip
\noindent
{\bf Theorem 7.1.}\quad{\it Let us assume that the function $\eta$ is
defined over the set
$$
\{(\nu,p)\in\B{C}\times\B{Z}:\, |\Re\nu|\le\sigma,\,|p|\le
l\}\leqno(7.8)
$$
with a fixed $\sigma>1$, and satisfies the conditions:
\item{1.} $\eta(\nu,p)$ is holomorphic on a neighbourhood of the strip
$|\Re\nu|\le\sigma$,
\item{2.} $\eta(\nu,p)\ll e^{-\pi|\Im\nu|/2}(1+|\Im\nu|)^{-A}$ for any
$A>0$,
\item{3.} $\eta(\nu,p)=\eta(-\nu,-p)$.
\par
\noindent
We put
$$
\leqalignno{ \e{M}_{l,q}^\omega\eta(\r{g})={1\over2\pi^3i}
\sum_{|p|\le l}&{(-i\omega/|\omega|)^p\over \Vert\Phi_{p,q}^l\Vert_K}
\int_{(0)}\eta(\nu,p) (\pi|\omega|)^{-\nu}\Gamma(l+1+\nu)&(7.9)\cr
&\times\e{A}_\omega \varphi_{l,q}(\nu,p)
(\r{g})\nu^{\epsilon(p)}\sin\pi\nu d\nu }
$$
with $\epsilon(0)=1$ and $\epsilon(p)=-1$ for $p\not=0$. Then
$\e{M}_{l,q}^\omega\eta\in P_{l,q}(N\backslash G, \omega)$, and we
have}
$$
\e{L}_{l,q}^\omega\e{M}^\omega_{l,q}\eta(\nu,p)
={2\over\pi}{\nu^{\epsilon(p)+1}\over p^2-\nu^2}\prod_{1\le j\le
l}(j^2-\nu^2)\cdot \eta(\nu,p). \leqno(7.10)
$$
Proof. We shall consider the first assertion. To estimate
$\e{M}_{l,q}^\omega\eta(\r{n}\r{a}[r]\r{k})$ as $r\uparrow\infty$, we
note that the integral formula
$$
K_\xi(u)={ \Gamma(\xi+\txt{1\over2})\over2\sqrt\pi(u/2)^\xi}
\int_0^\infty{e^{iux}\over (1+x^2)^{\xi+{1\over2}}}dx,\leqno(7.11)
$$
which holds for $u>0,\,\Re\xi>-{1\over2}$, gives, after a multiple use
of partial integration,
$$
K_\xi(u)\ll e^{-\pi|\xi|/2}((1+|\xi|)/u)^{\Re\xi+k}\leqno(7.12)
$$
for each fixed $k\ge1$, uniformly for $|\Re\xi|<{1\over2}k$, $u>0$.
This and $(5.26)$--$(5.27)$ imply that
$\e{M}_{l,q}^\omega\eta(\r{n}\r{a}[r]\r{k})$ is of rapid decay with
respect to $r$ as $r\uparrow\infty$. Next, to treat the case where
$r\downarrow0$ we observe that by $(7.12)$ the contour in $(7.9)$ can
be shifted to $(\alpha)$ with $0<\alpha<\sigma$. Then $(3.23)$,
$(6.15)$, and the condition {\it 3\/} give
$$
\leqalignno{ \e{M}_{l,q}^\omega\eta(\r{g})={i\over2\pi} \sum_{|p|\le
l}&{(i\omega/|\omega|)^{-p}\over\Vert\Phi_{p,q}^l\Vert_K}
\left\{\int_{(\alpha)}+\int_{(-\alpha)}\right\}
\eta(\nu,p)(\pi|\omega|)^\nu&(7.13)\cr
&\times\Gamma(l+1+\nu)\e{B}_\omega \varphi_{l,q}(\nu,p)
(\r{g})\nu^{\epsilon(p)}d\nu. }
$$
We may shift the contour $(-\alpha)$ to $(\alpha)$, and have
$$
\leqalignno{ \e{M}_{l,q}^\omega\eta(\r{g})&={i\over\pi} \sum_{|p|\le
l}{(i\omega/|\omega|)^{-p}\over\Vert\Phi_{p,q}^l\Vert_K}
\int_{(\alpha)} \eta(\nu,p)(\pi|\omega|)^\nu&(7.14)\cr
&\hskip 1cm\times\Gamma(l+1+\nu)\e{B}_\omega \varphi_{l,q}(\nu,p)
(\r{g})\nu^{\epsilon(p)}d\nu\cr
&+l!\sum_{1\le|p|\le l}{(i\omega/|\omega|)^{-p}\over
\Vert\Phi_{p,q}^l\Vert_K}\eta(0,p)
\e{B}_\omega\varphi_{l,q}(0,p)(\r{g}). }
$$
The formulas $(6.13)$--$(6.14)$ imply that as $r\downarrow0$ the first
sum on the right is $O(r^{1+\alpha})$, and also
$$
\leqalignno{
&l!\sum_{1\le|p|\le l}{(i\omega/|\omega|)^{-p}\over
\Vert\Phi_{p,q}^l\Vert_K}\eta(0,p)
\e{B}_\omega\varphi_{l,q}(0,p)(\r{g})&(7.15)\cr
&=b(\eta)\psi_\omega(\r{n})r^2\Phi_{0,q}^l(\r{k})+O(r^3)\cr
&=b(\eta)\e{B}_\omega\varphi_{l,q}(1,0)(\r{g})+O(r^3),
\cr
}
$$
where $l\ge1$, and $b(\eta)=-\pi l\cdot
l!\eta(0,1)|\omega|\Vert\Phi_{1,q}^l\Vert_K^{-1} $. Collecting these
we have indeed $\e{M}_{l,q}^\omega\eta\in P_{l,q}(N\backslash G,
\omega)$. We note that the last line in $(7.15)$ is to play a r\^ole
in Section 9.
\par
Consequently we may use $(7.7)$ in computing
$\e{L}_{l,q}^\omega\e{M}_{l,q}^\omega \eta(\nu,p)$. We thus apply
$\e{A}_0$ to the identity $(7.14)$. It is easy to check the absolute
convergence that is necessary to exchange the order of integration,
and by $(6.19)$ we have
$$
\leqalignno{
&\e{A}_0\e{M}_{l,q}^\omega\eta(\r{g})={i\over\pi} \sum_{|p|\le
l}{(-i\omega/|\omega|)^p\Phi_{p,q}^l(\r{k})
\over\Vert\Phi_{p,q}^l\Vert_K}&(7.16)\cr
&\times\int_{(\alpha)} \eta(-\nu,p)(\pi|\omega|)^\nu
\Gamma(l+1-\nu){\nu^{\epsilon(p)}\sin\pi\nu\over
\nu^2-p^2}r^{1-\nu}d\nu.\cr
}
$$
This and $(7.7)$ yield, via the Mellin inversion,
$$
\e{L}_{l,q}^\omega\e{M}_{l,q}^\omega\eta(\nu,p) =-2\pi^{-2}
\Gamma(l+1+\nu) \Gamma(l+1-\nu){\nu^{\epsilon(p)}\sin\pi\nu\over
\nu^2-p^2}\eta(\nu,p),\leqno(7.17)
$$
which ends the proof.
\smallskip
The above discussion implies in particular that
$\e{M}_{l,q}^\omega\eta\in L^2(N\backslash G)$. Related to this we
shall show a Parseval property of the transform $\e{M}_{l,q}^\omega$:
\smallskip
\noindent
{\bf Lemma 7.1.}\quad{\it Let $\eta$ and $\theta$ satisfy the three
conditions in the last theorem. Then we have}
$$
\leqalignno{
&\int_{N\backslash G}\e{M}_{l,q}^\omega\eta(\r{g})
\overline{\e{M}_{l,q}^\omega\theta(\r{g})}\,d\dot\r{g}&(7.18)\cr
&={1\over\pi^2i}\sum_{|p|\le l}\int_{(0)}\eta(\nu,p)
\overline{\theta(\nu,p)}{\nu^{2\epsilon(p)+1}\sin\pi\nu \over
p^2-\nu^2}\prod_{j=1}^l(j^2-\nu^2)\,d\nu. }
$$
Proof. We replace $\e{M}_{l,q}^\omega\theta(\r{g})$ by its defining
expression. The resulting double integral over $N\backslash G\times
i\B{R}$ is easily seen to be absolutely convergent; and consequently
we have
$$
\int_{N\backslash G}\e{M}_{l,q}^\omega\eta(\r{g})
\overline{\e{M}_{l,q}^\omega\theta(\r{g})}\,d\dot\r{g}={1\over2\pi
i}\sum_{|p|\le l}\int_{(0)}\e{L}_{l,q}^\omega\e{M}_{l,q}^\omega
\eta(\nu,p)\overline{\theta(\nu,p)}\nu^{\epsilon(p)}
\sin\pi\nu\,d\nu,\quad\leqno(7.19)
$$
which gives the assertion.
\smallskip
Further, we shall show
\smallskip
\noindent
{\bf Lemma 7.2.}\quad{\it For any non-zero $\omega_1$, $\omega_2$,
$\tau$ we define the map
$$
\kappa(\omega_1,\omega_2,\tau):\,\eta\mapsto
\e{K}_{\nu,p}(2\pi\tau\sqrt{\omega_1\omega_2})\eta,\leqno(7.20)
$$
where
$$
\e{K}_{\nu,p}(\xi)={1\over\sin\pi\nu} \{\e{J}_{-\nu,-p}(\xi)-
\e{J}_{\nu,p}(\xi)\}\leqno(7.21)
$$
with $\e{J}_{\nu,p}$ defined in $(6.21)$. Then we have, for $\eta$ as
in the last theorem,}
$$
\e{A}_{\omega_1}\ell_\tau\e{M}_{l,q}^{\omega_2}\eta(\nu,p)
=|\pi\tau|^2\e{M}_{l,q}^{\omega_1}\kappa(\omega_1,\omega_2, \tau)
\eta(\nu,p).\leqno(7.22)
$$
Proof. It is trivial that $\kappa(\omega_1,\omega_2,\tau)\eta$ satisfy
the three conditions in the last theorem. Hence the right side of
$(7.22)$ is well-defined. To transform the left side we use $(7.14)$.
Formally we have
$$
\leqalignno{
\e{A}_{\omega_1}\ell_\tau&\e{M}_{l,q}^{\omega_2}\eta(\nu,p)(\r{g})
={i\over\pi}|\tau|^2 \sum_{|p|\le
l}{(i\tau^2\omega_2/|\tau^2\omega_2|)^{-p}
\over\Vert\Phi_{p,q}^l\Vert_K} \int_{(\alpha)}
\eta(\nu,p)(\pi|\tau^2\omega_2|)^\nu&(7.23)\cr
&\hskip 1cm\times\Gamma(l+1+\nu)\e{A}_{\omega_1}\e{B}_{\tau^2\omega_2}
\varphi_{l,q}(\nu,p) (\r{g})\nu^{\epsilon(p)}d\nu\cr
&+l!|\tau|^2\sum_{1\le|p|\le
l}{(i\tau^2\omega_2/|\tau^2\omega_2|)^{-p}\over
\Vert\Phi_{p,q}^l\Vert_K}\eta(0,p)
\e{A}_{\omega_1}\e{B}_{\tau^2\omega_2} \varphi_{l,q}(0,p)(\r{g}),\cr
}
$$
where we have used $(6.18)$. To verify the exchange of the order of
integrals implicit in $(7.23)$ we need only to invoke $(7.15)$; note
that it also allows us to use $(6.20)$ even for
$\e{A}_{\omega_1}\e{B}_{\tau^2\omega_2} \varphi_{l,q}(0,p)$,
$p\not=0$. Thus we have, by $(6.20)$,
$$
\leqalignno{
&\e{A}_{\omega_1}\ell_\tau\e{M}_{l,q}^{\omega_2}\eta(\nu,p)(\r{g})
={i\over\pi}|\tau|^2 \sum_{|p|\le
l}{(-i\omega_1/|\omega_1|)^{p}\over\Vert\Phi_{p,q}^l\Vert_K}
\int_{(\alpha)} \eta(\nu,p)(\pi|\omega_1|)^{-\nu}&(7.24)\cr
&\hskip 1cm\times \Gamma(l+1+\nu)\e{J}_{\nu,p}
(2\pi\tau\sqrt{\omega_1\omega_2})\e{A}_{\omega_1}
\varphi_{l,q}(\nu,p) (\r{g})\nu^{\epsilon(p)}d\nu\cr
&+l!|\tau|^2\sum_{1\le|p|\le l}{(-i\omega_1/|\omega_1|)^{p}\over
\Vert\Phi_{p,q}^l\Vert_K}\eta(0,p)\e{J}_{0,p}
(2\pi\tau\sqrt{\omega_1\omega_2}) \e{A}_{\omega_1}
\varphi_{l,q}(0,p)(\r{g}).\cr
}
$$
We shift the contour $(\alpha)$ of one half of the last integral to
$(-\alpha)$; then the last sum disappears. To the integrand over
$(-\alpha)$ we apply the functional equation $(5.29)$. After a
rearrangement we get $(7.22)$.
\smallskip
\noindent
{\csc Remark.} This section is a detailed work-out of the last chapter
of [41] in the case of $\r{PSL}_2(\B{C})$. It is in fact the harmonic
analysis of the space of $N$-equivariant functions. The
$\e{L}_{l,q}^\omega$ could be called a Whittaker transform, but we
regard it rather as an extension of the Lebedev transform, paying
respect to its origin. For $(7.1)$--$(7.2)$ see Section 2.6 of [30].
For an interpretation of $\e{K}_{\nu,p}$ see Section 15.
\medskip
\centerline{\bf 8. The space $L^2(\varGamma\backslash G)$}
\smallskip
\noindent
In the next section we shall treat inner-products of certain
Poincar\'e series, especially their spectral decompositions. Here we
shall briefly develop the relevant spectral theory of the space
$L^2(\varGamma\backslash G)$ composed of all left
$\varGamma$-automorphic functions on $G$ which are square integrable
over $\varGamma\backslash G$ with respect to the measure induced by
$d\r{g}$. To this end we shall employ the unitary representation of
$G$ realized over $L^2(\varGamma\backslash G)$ via right translations
by elements of $G$. We shall see that automorphic forms on $G$,
especially the basis elements for the Parseval formula over
$L^2(\varGamma\backslash G)$, do not occur singly but are
parametrized through maps of the model space $H(\nu,p)$ and live in
right-irreducible subspaces of $L^2(\varGamma\backslash G)$ sharing
Fourier coefficients, as is indicated at the end of Section 5. This
point of view will be essential in describing the sum formula for
$S_\r{F}$.
\smallskip
We first observe that the constant function and all cusp-forms over
$G$ belong to $L^2(\varGamma\backslash G)$, because of
$(3.8)$--$(3.10)$ and $(4.25)$. We have
$$
L^2(\varGamma\backslash G)=\B{C}\oplus {}^0\!L^2(\varGamma\backslash
G)\oplus {}^e\!L^2(\varGamma\backslash G).\leqno(8.1)
$$
Here ${}^0\!L^2(\varGamma\backslash G)$ is the subspace spanned by all
cusp-forms, and called the cuspidal subspace; the subspace
${}^e\!L^2(\varGamma\backslash G)$ is the orthogonal complement. The
space ${}^0\!L^2(\varGamma\backslash G)$ is $G$-invariant, and we
have a decomposition
$$
{}^0\!L^2(\varGamma\backslash G)=\overline{\bigoplus V}\leqno(8.2)
$$
into countably many subspaces $V$ irreducible with respect to the
action of $G$. Each $V$ has a dense subspace that is a common
eigenspace of the Casimir elements, and Lemma 4.1 implies that we
should have
$$
\Omega_\pm|_V=\chi_{\nu_V,p_V}(\Omega_\pm)\cdot 1.\leqno(8.3)
$$
 It is known
that for the group $\varGamma$ all $V$ are of unitary principal
series type, and we can suppose that
$$
\nu_V\in i[0,\infty).\leqno(8.4)
$$
Moreover, there exists a linear isomorphism
$$
T_V:\,H(\nu_V,p_V)\to V,\leqno(8.5)
$$
which has a dense image and commutes with the action of $\f{g}$. This
has the following immediate consequences: We have the decomposition
$$
V=\overline{\bigoplus_{|p_V|\le l,\,|q|\le l}V_{l,q}},\quad V_{l,q}=
\B{C}\,T_V\varphi_{l,q}(\nu_V,p_V),\leqno(8.6)
$$
where $V_{l,q}$ is the subspace spanned by all cusp-forms of $K$-type
$(l,q)$ in $V$; that is, $\dim V_{l,q}=1$. Besides, the unitary
structure of $H(\nu_V,p_V)$ mentioned at the end of Section 3 is
transferred by $T_V$ into $(8.6)$, and we have
$$
\Vert T_V\varphi_{l,q}(\nu_V,p_V)\Vert_{\varGamma\backslash G}=
\Vert\Phi_{p_V,q}^l\Vert_K\leqno(8.7)
$$
with the norm $\Vert\cdot\Vert_{\varGamma\backslash G}$ corresponding
to $(8.9)$. By $(5.35)$ we have the Fourier expansion
$$
T_V\varphi_{l,q}(\nu_V,p_V)=\sum_{\omega\ne0}c_V(\omega)
\e{A}_\omega\varphi_{l,q}(\nu_V,p_V).\leqno(8.8)
$$
The Fourier coefficients $c_V(\omega)$ depends only on $V$ and
$\omega$. This is because both $T_V$ and $\e{A}_\omega$ commute with
the action of $\f{g}$. The vector $\{c_V(\omega)\}$ is fixed by $V$
up-to an arbitrary multiplier of unit absolute value.
\smallskip
We now restrict the decomposition $(8.1)$ to the subspace
$L^2(\varGamma\backslash G)_{l,q}$ spanned by all square-integrable
left $\varGamma$-automorphic functions of $K$-type $(l,q)$. Then the
cuspidal part is well described by the above assertions. What remains
is the non-cuspidal part, and it is rendered in terms of Eisenstein
series of $K$-type $(l,q)$, as is embodied in the fundamental
\smallskip
\noindent
{\bf Theorem 8.1.}\quad{\it Let $f_1, f_2\in L^2(\varGamma\backslash
G)_{l,q}$, and denote their inner-product by
$$
\langle{f_1,f_2}\rangle_{\varGamma\backslash G}
=\int_{\varGamma\backslash G}f_1(\r{g})\overline{f_2(\r{g})}d\r{g}.
\leqno(8.9)
$$
Then we have the Parseval identity
$$
\leqalignno{
\qquad\langle{f_1,f_2}\rangle_{\varGamma\backslash G}&=\delta_{l,0}
{\pi^2\over2\zeta_\r{F}(2)}\langle{f_1,1}
\rangle_{\varGamma\backslash G}\langle{1,f_2}
\rangle_{\varGamma\backslash G}&(8.10)\cr
&+\sum_V {1\over\Vert\Phi_{p_V,q}^l\Vert_K^2}
\langle{f_1,T_V\varphi_{l,q}(\nu_V,p_V)} \rangle_{\varGamma\backslash
G} \langle{T_V\varphi_{l,q}(\nu_V,p_V),f_2}
\rangle_{\varGamma\backslash G}\cr
&+\sum_{\scr{|p|\le l}\atop\scr{p\in2\B{Z}}}{1\over2\pi
i\Vert\Phi_{p,q}^l \Vert_K^2}\int_{(0)}E_{l,q}[\nu,p;f_1]
\overline{E_{l,q}[\nu,p;f_2]}d\nu,\cr
}
$$
where the convergence is absolute throughout. Here $V$ runs over a
complete orthogonal system of right-irreducible cuspidal subspaces of
$L^2(\varGamma\backslash G)$ that intersect the space
$L^2(\varGamma\backslash G)_{l,q}$ non-trivially. Also
$$
E_{l,q}[\nu,p;f_j]= \int_{\varGamma\backslash
G}f_j(\r{g})\overline{e_{l,q} (-\bar{\nu},p;\r{g})}d\r{g}\leqno(8.11)
$$
in the sense of norm convergence. }
\smallskip
\noindent
Proof. This is a special case of a general result due to Langlands
[23] (see also [13]). We stress, however, that our particular
assertion could be established in a direct way.
\smallskip
Next, we shall take into consideration the action of Hecke operators:
We define the Hecke operator labeled with $n\in\B{Z}[i]$ by
$$
{\cal T}_n:\psi(\r{g})\mapsto{1\over4|n|}
\sum_{d|n}\;\sum_{b\bmod\,d}\psi(\r{n}[b/d]
\r{h}[\sqrt{n}/d]\r{g}),\leqno(8.12)
$$
where $\psi$ is to be left $\varGamma$-automorphic; the choice of the
square root is irrelevant. If $\psi\in{}^0\!A_{l,q}(\chi_{\nu,p})$,
then we have, from $(5.8)$ and $(5.35)$,
$$
{\cal T}_n\psi={(n/|n|)^p\over4|n|^\nu}\sum_{\omega\ne0}c(\omega)
\sum_{d|(\omega,n)}|d|^{2\nu}(d/|d|)^{-2p}\e{A}_{\omega n/d^2}
\varphi_{l,q}(\nu,p).\leqno(8.13)
$$
The commutativity of the algebra $\{{\cal T}_n\}$ and the metric
property of each ${\cal T}_n$ in $L^2(\varGamma\backslash G)$ are
analogous to the rational case. Since the right side of $(8.12)$
commutes with the right translation by elements of $G$, we may assume
that every $V$ in $(8.2)$ is an eigenspace of ${\cal T}_n$ for all
$n$ with the eigenvalue $t_V(n)\in{\Bbb R}$. The Weil bound
$$
S_\r{F}(\omega_1,\omega_2;c)\ll |(\omega_1,\omega_2,c)|
|c|\sigma_0(c,0)\leqno(8.14)
$$
yields the estimate
$$
t_V(n)\ll |n|^{1/2+\varepsilon}\leqno(8.15)
$$
for any fixed $\varepsilon>0$, where the implicit constant depends
only on $\varepsilon$; see Corollary 10.1 below. The equation
$(8.13)$ implies that in $(8.8)$
$$
t_V(n)c_V(\omega)={1\over4}|n|^{\nu_V}(n/|n|)^{-p_V}
\sum_{d|(\omega,n)} |d|^{-2\nu_V}(d/|d|)^{2p_V}c_V(\omega
n/d^2).\leqno(8.16)
$$
In particular we have, for all $n\in\B{Z}[i]$,
$$
c_V(n)=c_V(1)|n|^{-\nu_V}(n/|n|)^{p_V}t_V(n),\leqno(8.17)
$$
where we have used $(5.37)$. This and $(8.16)$ give
$$
t_V(m)t_V(n)={1\over4}\sum_{d|(m,n)}t_V(mn/d^2).\leqno(8.18)
$$
We have
$$
t_V(1)=1,\quad t_V(-n)=t_V(n),\quad
t_V(in)=\epsilon_Vt_V(n)\leqno(8.19)
$$
with $\epsilon_V=\pm1$.
\par
Further, let
$$
H_V(s)={1\over4}\sum_{n\ne0}t_V(n) |n|^{-2s}\leqno(8.20)
$$
be the Hecke series associated with the irreducible subspace $V$ under
the convention $(8.17)$; note that when $\epsilon_V=-1$ this vanishes
identically. The identity $(8.18)$ implies that in the region of
absolute convergence
$$
H_V(s_1)H_V(s_2) ={1\over4}\zeta_\r{F}(s_1+s_2)\sum_{n\ne0}
\sigma_{s_1-s_2}(n)t_V(n)|n|^{-2s_1}.\leqno(8.21)
$$
Properties of $H_V(s)$ as a function of $s$ can be read from
\smallskip
\noindent
{\bf Lemma 8.1.}\quad{\it Let $b\in\B{Z}$, and put
$$
H_V(s,b)={1\over4}\sum_{n\ne0}t_V(n)(n/|n|)^b|n|^{-2s}. \leqno(8.22)
$$
Take $b\in2\B{Z}$ and $i^b=\epsilon_V$ to have a non-trivial sum. Then
$H_V(s,b)$ is entire in $s$ and satisfies the functional equation}
$$
\leqalignno{
&\pi^{1-2s}\Gamma(s+\txt{1\over2}(|p_V+b|+\nu_V))
\Gamma(s+\txt{1\over2}(|p_V-b|-\nu_V))H_V(s,b)&(8.23)\cr
&=(-1)^{\max(|b|,|p_V|)}
\pi^{2s-1}\Gamma(1-s+\txt{1\over2}(|p_V-b|+\nu_V))\cr
&\hskip1.5cm\times\Gamma(1-s+\txt{1\over2}
(|p_V+b|-\nu_V))H_V(1-s,-b).\cr
}
$$
Proof. We consider the integral
$$
\leqalignno{
&\qquad Y_V(s,b;\r{k})&(8.24)
\cr
&=\int_{\B{C}^\times}|u|^{4s-2}
(u/|u|)^{-2b}T_V\varphi(\r{h}[u]\r{k})d^\times\!u
\cr
&=\int_{|u|\ge1}\Big(|u|^{4s-2}(u/|u|)^{-2b}
T_V\varphi(\r{h}[u]\r{k})+|u|^{2-4s}(u/|u|)^{2b}T_V
\varphi(\r{h}[u]\r{w}\r{k})\Big)d^\times\!u,\qquad }
$$
where $\varphi=\varphi_{l,l}(\nu_V,p_V)$ with $l=\max(|b|,|p_V|)$, and
$\r{w}$ is as in $(5.2)$. Obviously $Y_V(s,b;\r{k})$ is entire in
$s$, and $Y_V(1-s,-b;\r{w}\r{k})=Y_V(s,b;\r{k})$. In view of $(5.8)$,
$(8.8)$, and $(8.17)$, we have, for $\Re s$ sufficiently large,
$$
Y_V(s,b;\r{k})=4c_V(1)H_V(s,b)\int_{\B{C}^\times}
|u|^{4s-2}(u/|u|)^{-2b}\e{A}_1\varphi
(\r{h}[u]\r{k})d^\times\!u.\leqno(8.25)
$$
By $(5.26)$--$(5.27)$ this integral is equal to
$$
\leqalignno{
&2(-1)^{l-p_V}i^{b-p_V}\pi^{1+\nu_V}\Phi_{-b,l}^l(\r{k}) \int_0^\infty
r^{2s-2}\alpha_{-b}^l(\nu_V,p_V;r)dr&(8.26)\cr
=&2{(-1)^{l-p_V}i^{b-p_V} \pi^{l+2+\nu_V}\over\Gamma(l+1+\nu)}
\xi_{p_V}^l(-b,0)\Phi_{-b,l}^l (\r{k}) \int_0^\infty
r^{2s+l-1}K_{l+\nu-|p_V-b|}(2\pi r)dr.
\cr
}
$$
On noting that $\Phi_{b,l}^l(\r{w}\r{k})=(-1)^{l+b}
\Phi_{-b,l}^l(\r{k})$, we obtain $(8.23)$.
\smallskip
\noindent
{\csc Remark.} For the spectral theory of automorphic forms on
semisimple Lie groups see, e.g., [13]. The assertion $(8.4)$ depends
on the absence of exceptional eigenvalues for the non-Euclidean
Laplacian over $\varGamma\backslash\B{H}^3$ (see Proposition 6.2 in
Chapter 7 of [10]). That suffices, as the complementary series occurs
only for $p=0$. The bound $(8.14)$ is a special case of Theorem 10 of
[2], which applies to all number fields. It should be stressed that
for our purpose it is enough to have any non-trivial exponent in
place of $1+\varepsilon$, which is best possible. In the proof of the
last lemma we followed [18]; see the proof of Theorem 6.4 there.
\par
As to $(8.2)$ it may be worth mentioning the following {\it
multiplicity one} result: For given $\pm(\nu,p)\in i\B{R}\times\B{Z}$
and $\{t(n)\in\B{C}: n\in \B{Z}[i],\, n\neq0\}$, there is at most one
irreducible subspace $V$ of
$\displaystyle{{}^0\!L^2(\varGamma\backslash G)}$ with $(\nu_V,p_V) =
\pm (\nu,p)$, and $t_V(n)=t(n)$ for all $n$. Indeed, the Fourier
expansion $(5.35)$ and $(8.17)$ show that an automorphic form of a
given $K$-type is determined by $(\nu,p)$ and the $t(n)$ up to a
scalar factor. This shows that the decomposition $(8.2)$ is unique if
we impose the condition that the spaces $V$ are invariant under all
Hecke operators.
\medskip
\centerline{\bf 9. Preliminary sum formula}
\smallskip
\noindent
We now enter into the discussion of the sum formula for Kloosterman
sums $S_\r{F}$. The formula will be derived via spectral and
geometric computations of an inner-product of two particular
Poincar\'e series. This is analogous to the rational case. However,
the choice of these series gives rise to a discussion. A possible way
to take is to use an explicit function as a seed to generate the
Poincar\'e series, which extends Selberg's argument for the rational
case. This works well if we restrict ourselves to the $K$-trivial
case, but it does not seem to extend easily to the $K$-non-trivial
situation. On the other hand, a method that Miatello--Wallach [25],
[26] developed for a far more general situation offers us a flexible
way to choose the seed function. Here we shall follow their argument,
adopting it to our present specifications.
\smallskip
Thus we shall employ $\e{M}_{l,q}^\omega\eta$ to generate a Poincar\'e
series, where $\eta$ is to satisfy the three conditions given in
Theorem 7.1. We notice immediately that $(7.15)$ causes, in general,
a convergence problem. This reminds us a similar situation that Hecke
encountered in his investigation of holomorphic modular forms of
weight 2. He used analytic continuation to overcome the difficulty.
In much the same spirit we shall consider, in view of the last line
of $(7.15)$, the sum
$$
[\e{B}_\omega]\varphi_{l,q}(\nu,p)(\r{g})={1\over2}
\sum_{\gamma\in\varGamma_N\backslash\varGamma}
\e{B}_\omega\varphi_{l,q}(\nu,p)(\gamma\r{g}) \leqno(9.1)
$$
with non-zero $\omega\in\B{Z}[i]$, though we actually need only the
case $p=0$. By $(5.30)$ and $(6.13)$ we see that the sum converges
absolutely, for $\Re\nu>1$, to a left $\varGamma$-automorphic
function of $K$-type $(l,q)$ with character $\chi_{\nu,p}$. The
combination of $(5.5)$--$(5.7)$ and $(6.18)$--$(6.21)$ yields that
$$
\leqalignno{ [\e{B}_\omega]\varphi_{l,q}(\nu,p)(\r{g})&={1\over2}
(\e{B}_\omega+(-1)^p\e{B}_{-\omega})\varphi_{l,q}
(\nu,p)(\r{g})&(9.2)\cr
&+(-1)^p{\sin\pi\nu\over\nu^2-p^2}
{\Gamma(l+1-\nu)\over\Gamma(l+1+\nu)}
{\sigma_{-\nu}(\omega,p/2)\over\zeta_\r{F}(1+\nu,p/2)}
\varphi_{l,q}(-\nu,-p)(\r{g})\cr
&+\sum_{\omega'\ne0}{\cal J}_{\nu,p}(\omega,\omega')
\e{A}_{\omega'}\varphi_{l,q}(\nu,p)(\r{g}),
\cr
}
$$
where the second line appears only when $p\in2\B{Z}$; and
$$
{\cal J}_{\nu,p}(\omega,\omega')={1\over4\pi^2|\omega\omega'|^\nu}
\left({\omega\omega'\over |\omega\omega'|}\right)^p\sum_{c\ne0}
{1\over|c|^2}S_\r{F}(\omega,\omega';c) \e{J}_{\nu,p}\left({2\pi\over
c}\sqrt{\omega\omega'}\right).\leqno(9.3)
$$
The bound $(8.14)$ implies that for $\Re\nu>{1\over2}$ the function
${\cal J}_{\nu,p} (\omega,\omega')$ is regular and of polynomial
order in $\omega$, $\omega'$. Thus
$[\e{B}_\omega]\varphi_{l,q}(\nu,p)(\r{g})$ is regular for
$\Re\nu>{1\over2}$, and analytically continues to a left
$\varGamma$-automorphic function of $K$-type $(l,q)$ with character
$\chi_{\nu,p}$. It is, however, of exponential growth with respect to
$r$, $\r{g}=\r{n}\r{a}[r]\r{k}$; and thus it does not belong to
$L^2(\varGamma\backslash G)$. We then appeal to a common practice: we
attach a factor $\rho(\gamma\r{g})$ to each summand of $(9.1)$. Here
$\rho(\r{g})=\rho(r)$ with an abuse of notation, and $\rho(r)$ is
smooth, being equal to $1$ for $r\le r_0$ and to $0$ for $r>r_0+1$
with $r_0\ge2$. If $r>1$ then this affects actually only two terms in
$(9.1)$, which correspond to the cosets represented by $\gamma=1$ and
$\r{h}[i]$. We thus have, instead of $(9.1)$--$(9.2)$, that
$$
[\rho\e{B}_\omega]\varphi_{l,q}(\nu,p)(\r{g})={1\over2}
\sum_{\gamma\in\varGamma_N\backslash\varGamma}
\rho(\gamma\r{g})\e{B}_\omega\varphi_{l,q}(\nu,p)(\gamma\r{g})
\leqno(9.4)
$$
for $\Re\nu>1$, and that if $r>1$, $\Re\nu>{1\over2}$,
$$
\leqalignno{ [\rho\e{B}_\omega]\varphi_{l,q}(\nu,p)(\r{g})&=
[\e{B}_\omega]\varphi_{l,q}(\nu,p)(\r{g})&(9.5)\cr
&+{1\over2}(\rho(r)-1)(\e{B}_\omega+(-1)^p
\e{B}_{-\omega})\varphi_{l,q} (\nu,p)(\r{g}). }
$$
Note that as $r\uparrow\infty$
$$
[\rho\e{B}_\omega]\varphi_{l,q}(\nu,p)(\r{g}) \ll r^{1-\Re\nu}
\leqno(9.6)
$$
uniformly for $\Re\nu>{1\over2}$.
\smallskip
Returning to $(7.14)$--$(7.15)$ we define
$\e{M}_{l,q}^{\omega,\ast}\eta$ by
$$
\e{M}_{l,q}^\omega\eta=\e{M}_{l,q}^{\omega,\ast}
\eta+b(\eta)\rho\e{B}_\omega \varphi_{l,q}(1,0),\leqno(9.7)
$$
and $[\e{M}_{l,q}^\omega]\eta$ by
$$
[\e{M}_{l,q}^\omega]\eta=[\e{M}_{l,q}^{\omega,\ast}]
\eta+b(\eta)[\rho\e{B}_\omega] \varphi_{l,q}(1,0),\leqno(9.8)
$$
where
$$
[\e{M}_{l,q}^{\omega,\ast}]\eta={1\over2}\sum_{\gamma\in
\varGamma_N\backslash
\varGamma}\e{M}_{l,q}^{\omega,\ast}\eta(\gamma\r{g}).\leqno(9.9)
$$
By the construction we have that
$\e{M}_{l,q}^{\omega,\ast}\eta(\r{g})\ll r^{2+\varepsilon}$ with an
$\varepsilon>0$ as $r\downarrow0$, and $\ll r^{-A}$ for any $A>0$ as
$r\uparrow\infty$. Thus $[\e{M}_{l,q}^{\omega,\ast}] \eta(\r{g})\ll
r^{-\varepsilon}$ as $r\uparrow\infty$. This and $(9.6)$ give
$$
[\e{M}_{l,q}^\omega]\eta(\r{g})\ll 1.\leqno(9.10)
$$
\par
Now, let $\eta$, $\theta$ satisfy the three conditions given in
Theorem 7.1, and let us consider the inner-product
$\langle[\e{M}_{l,q}^{\omega_1}]\eta,
[\e{M}_{l,q}^{\omega_2}]\theta\rangle_{\varGamma\backslash G}$ with
$\omega_1,\omega_2\in\B{Z}[i]$, $\omega_1\omega_2\ne0$. We are going
to apply Theorem 8.1 to it. To this end we note first that the above
discussion implies
$$
\langle[\e{M}_{l,q}^{\omega}]\eta, f\rangle_{\varGamma\backslash G}
=\langle[\e{M}_{l,q}^{\omega,\ast}]\eta,
f\rangle_{\varGamma\backslash G}
+b(\eta)\lim_{\nu\to1+0}\langle[\rho\e{B}_{\omega}]
\varphi_{l,q}(\nu,0),f\rangle_{\varGamma\backslash G} \leqno(9.11)
$$
for any left $\varGamma$-automorphic $f$ which is integrable over
$\varGamma\backslash G$. In this we have, by the unfolding argument,
$$
\leqalignno{ \langle[\e{M}_{l,q}^{\omega,\ast}]\eta,
f\rangle_{\varGamma\backslash G}&={1\over2}\int_{N\backslash G}
\e{M}_{l,q}^{\omega,\ast}\eta(\r{g})
\overline{F_{\omega}f(\r{g})}d\dot\r{g},&(9.12)
\cr
\langle[\rho\e{B}_{\omega}]
\varphi_{l,q}(\nu,0),f\rangle_{\varGamma\backslash G}&=
{1\over2}\int_{N\backslash G}\rho(\r{g})\e{B}_{\omega}
\varphi_{l,q}(\nu,0)(\r{g}) \overline{F_{\omega}f(\r{g})}d\dot\r{g}.
\cr}
$$
Thus, assuming that
$$
\leqalignno{ \lim_{\nu\to1+0}\int_{N\backslash
G}&\rho(\r{g})\e{B}_{\omega} \varphi_{l,q}(\nu,0)(\r{g})
\overline{F_{\omega}f(\r{g})}d\dot\r{g}&(9.13)\cr
&=\int_{N\backslash G}\rho(\r{g})\e{B}_{\omega}
\varphi_{l,q}(1,0)(\r{g}) \overline{F_{\omega}f(\r{g})}d\dot\r{g},
\cr
}
$$
we have
$$
\langle[\e{M}_{l,q}^{\omega}]\eta, f\rangle_{\varGamma\backslash G}
={1\over2}\int_{N\backslash G} \e{M}_{l,q}^{\omega}\eta(\r{g})
\overline{F_{\omega}f(\r{g})}d\dot\r{g}.\leqno(9.14)
$$
The functions $T_V\varphi_{l,q}(\nu_V,p_V)$ and $e_{l,q}(\nu,p)$ with
$\Re\nu=0$ are integrable over $\varGamma\backslash G$ and satisfy
$(9.13)$. Hence we have, on noting $(7.4)$, $(8.8)$, and $(8.17)$,
$$
\leqalignno{
&\langle[\e{M}_{l,q}^{\omega}]\eta,
T_V\varphi_{l,q}(\nu_V,p_V)\rangle_{\varGamma\backslash G}&(9.15)\cr
&={1\over2}(-i)^{p_V}\pi^{2-\nu_V} \overline{c_V(1)}t_V(\omega)
{\Vert\Phi_{p_V,q}^l\Vert_K
\over\Gamma(l+1-\nu_V)}\e{L}_{l,q}^{\omega}
\e{M}_{l,q}^{\omega}\eta(\nu_V,p_V)\cr
&=(-i)^{p_V}\pi^{-\nu_V}\overline{c_V(1)}t_V(\omega)
\Vert\Phi_{p_V,q}^l\Vert_K\Gamma(l+1+\nu_V)
{\nu_V^{\epsilon(p_V)}\sin\pi\nu_V \over
p_V^2-\nu_V^2}\eta(\nu_V,p_V) }
$$
by virtue of Theorem 7.1 or rather $(7.17)$. Similarly we have, from
$(5.32)$,
$$
\leqalignno{ E_{l,q}[\nu,p;[\e{M}_{l,q}^{\omega}]\eta]
=&(-1)^{p/2}{(\pi|\omega|)^{-\nu}(\omega/|\omega|)^{p}
\sigma_{\nu}(\omega,-p/2)\over\zeta_\r{F}(1-\nu,-p/2)}&(9.16)\cr
&\times\Vert\Phi_{p,q}^l\Vert_K\Gamma(l+1+\nu)
{\nu^{\epsilon(p)}\sin\pi\nu \over p^2-\nu^2}\eta(\nu,p)\cr
}
$$
for $\Re\nu=0$, $p\in2\B{Z}$. Further, we have obviously
$\langle[\e{M}_{l,q}^{\omega}]\eta, 1\rangle_{\varGamma\backslash
G}=0$.
\smallskip
Collecting these we get, by Theorem 8.1,
$$
\leqalignno{
&\qquad\langle[\e{M}_{l,q}^{\omega_1}]\eta,
[\e{M}_{l,q}^{\omega_2}]\theta\rangle_{\varGamma\backslash G} =\sum_V
|c_V(1)|^2t_V(\omega_1)t_V(\omega_2) \lambda_l(\nu_V,p_V)
\eta(\nu_V,p_V)\overline{\theta(\nu_V,p_V)}&(9.17)\cr
&\qquad+\sum_{\scr{|p|\le l}\atop\scr{p\in2\B{Z}}}{(\omega_1\omega_2
/|\omega_1\omega_2|)^p\over2\pi i}
\int_{(0)}{\sigma_\nu(\omega_1,-p/2)\sigma_\nu(\omega_2,-p/2)
\over|\omega_1\omega_2|^\nu|\zeta_\r{F}(1+\nu,p/2)|^2}
\lambda_l(\nu,p)\eta(\nu,p)\overline{\theta(\nu,p)}\,d\nu,
\cr
}
$$
where $V\cap L^2(\varGamma\backslash G)_{l,q}\ne\{0\}$, and
$$
\lambda_l(\nu,p)=\Gamma(l+1+\nu)\Gamma(l+1-\nu)
\left({\nu^{\epsilon(p)}\sin\pi\nu\over\nu^2-p^2}\right)^2.
\leqno(9.18)
$$
\par
Next, we move to the geometric computation of the inner-product. To
this end we make a trivial observation that
$$
F_{\omega_1}[\e{M}_{l,q}^{\omega_2}]\theta=F_{\omega_1}
[\e{M}_{l,q}^{\omega_2,\ast}]\theta+b(\theta)\lim_{\nu\to1+0}
F_{\omega_1}[\rho\e{B}_{\omega_2}]\varphi_{l,q}(\nu,0). \leqno(9.19)
$$
Here we have, by $(5.6)$,
$$
\leqalignno{ F_{\omega_1}
[\e{M}_{l,q}^{\omega_2,\ast}]\theta&={1\over2}
\left(\delta_{\omega_1,\omega_2}
\e{M}_{l,q}^{\omega_2,\ast}\theta+\delta_{\omega_1,-\omega_2}
\ell_i\e{M}_{l,q}^{\omega_2,\ast}\theta\right)&(9.20)
\cr
&+{1\over4}\sum_{c\ne0}S_\r{F}(\omega_1,\omega_2;c)
\e{A}_{\omega_1}\ell_{1/c}\e{M}_{l,q}^{\omega_2,\ast}\theta }
$$
as well as
$$
\leqalignno{ F_{\omega_1}[\rho\e{B}_{\omega_2}]\varphi_{l,q}(\nu,0)&=
{1\over2} \left(\delta_{\omega_1,\omega_2}
\rho\e{B}_{\omega_2}\varphi_{l,q}(\nu,0)+ \delta_{\omega_1,-\omega_2}
\ell_i\rho\e{B}_{\omega_2}\varphi_{l,q}(\nu,0)\right)&(9.21)
\cr
&+{1\over4}\sum_{c\ne0}S_\r{F}(\omega_1,\omega_2;c)
\e{A}_{\omega_1}\ell_{1/c}\rho\e{B}_{\omega_2}\varphi_{l,q} (\nu,0).
}
$$
We have, by definition,
$$
\leqalignno{
&\e{A}_{\omega_1}\ell_{1/c}\rho\e{B}_{\omega_2}\varphi_{l,q}
(\nu,0)(\r{g})&(9.22)\cr
&=(c/|c|)^{2p}|c|^{-2(1+\nu)}\int_N
\psi_{\omega_1}(\r{n})^{-1}\rho(\r{h}[1/c]\r{w}\r{n}\r{g})
\e{B}_{\omega_2/c^2}\varphi_{l,q} (\nu,0)(\r{w}\r{n}\r{g})\,d\r{n}. }
$$
In view of $(5.9)$ this integral is bounded by a constant multiple of
$$
\int_{N^*}|\e{B}_{\omega_2/c^2}\varphi_{l,q}
(\nu,0)(\r{w}\r{n}[z]\r{g})|d_+z,\leqno(9.23)
$$
where $N^*=\{z: r\le2r_0|c|^2(r^2+|z|^2)\}$ with
$\r{g}=\r{a}[r]\r{k}$. By Lemma 6.1 we have, for $z\in N^*$,
$$
\e{B}_{\omega_2/c^2}\varphi_{l,q} (\nu,0)(\r{w}\r{n}[z]\r{g})\ll
\left({r\over|c|^2(r^2+|z|^2)}\right)^{1+\Re\nu},\leqno(9.24)
$$
where the implicit constant depends on $\nu$, $\omega_2$ but neither
on $c$ nor on $r$. In this way we get the estimate
$$
\e{A}_{\omega_1}\ell_{1/c}\rho\e{B}_{\omega_2}\varphi_{l,q}
(\nu,0)(\r{g})\ll r|c|^{-4-2\Re\nu}\leqno(9.25)
$$
uniformly for $r\ge0$, $\Re\nu>0$, and $\B{Z}[i]\ni c\ne0$. This means
that we may take the limit of $(9.19)$ inside the sum of $(9.21)$.
Then, by virtue of Lemma 7.2, we have
$$
\leqalignno{ F_{\omega_1} [\e{M}_{l,q}^{\omega_2}]\theta&={1\over2}
\left(\delta_{\omega_1,\omega_2}
\e{M}_{l,q}^{\omega_2}\theta+\delta_{\omega_1,-\omega_2}
\ell_i\e{M}_{l,q}^{\omega_2}\theta\right)&(9.26)
\cr
&+{\pi^2\over4}\sum_{c\ne0}{1\over|c|^2} S_\r{F}(\omega_1,\omega_2;c)
\e{M}_{l,q}^{\omega_1}\kappa(\omega_1,\omega_2,1/c)\theta. }
$$
The formulas $(7.22)$ and $(7.24)$ with appropriate changes of
notation imply that
$$
\sum_{c\ne0}{1\over|c|^2} |S_\r{F}(\omega_1,\omega_2;c)|
|\e{M}_{l,q}^{\omega_1}\kappa(\omega_1,\omega_2,1/c)\theta
(\nu,p)(\r{a}[r]\r{k})|\ll \cases{r^{-A}&as $r\uparrow\infty$,\cr
r^{1\over2}& as $r\downarrow0$,\cr
}\leqno(9.27)
$$
where $A>0$ is arbitrary. Thus $F_{\omega_1}
[\e{M}_{l,q}^{\omega_2}]\theta$ satisfies the condition $(9.13)$ with
$\omega=\omega_1$, and $(9.14)$ gives
$$
\langle[\e{M}_{l,q}^{\omega_1}]\eta, [\e{M}_{l,q}^{\omega_2}]\theta
\rangle_{\varGamma\backslash G} ={1\over2}\int_{N\backslash G}
\e{M}_{l,q}^{\omega_1}\eta(\r{g})
\overline{F_{\omega_1}[\e{M}_{l,q}^{\omega_2}]\theta(\r{g})}
\,d\dot\r{g}.\leqno(9.28)
$$
Moreover, we may insert $(9.26)$ into this and perform the integration
inside the infinite sum, getting
$$
\leqalignno{ \langle[&\e{M}_{l,q}^{\omega_1}]\eta,
[\e{M}_{l,q}^{\omega_2}]\theta \rangle_{\varGamma\backslash G}=
{1\over4}(\delta_{\omega_1,\omega_2}+
\delta_{\omega_1,-\omega_2})\int_{N\backslash G}
\e{M}_{l,q}^{\omega_1}\eta(\r{g})
\overline{\e{M}_{l,q}^{\omega_1}\theta (\r{g})}\,d\dot\r{g}&(9.29)\cr
&+{\pi^2\over8}\sum_{c\ne0}{1\over|c|^2}
S_\r{F}(\omega_1,\omega_2;c)\int_{N\backslash G}
\e{M}_{l,q}^{\omega_1}\eta(\r{g}) \overline{\e{M}_{l,q}^{\omega_1}
\kappa(\omega_1,\omega_2,1/c)\theta(\r{g})}\,d\dot\r{g}.
\cr}
$$
\par
Invoking Lemma 7.1, we have, from $(9.12)$ and $(9.29)$,
\smallskip
\noindent
{\bf Lemma 9.1.}\quad{\it Let $\eta$, $\theta$ satisfy the three
conditions given in Theorem 7.1. Then we have, for any non-zero
$\omega_1,\omega_2\in\B{Z}[i]$,
$$
\leqalignno{\qquad
&\sum_V|c_V(1)|^2t_V(\omega_1)t_V(\omega_2) \lambda_l(\nu_V,p_V)
\eta(\nu_V,p_V)\overline{\theta(\nu_V,p_V)}&(9.30)\cr
&+\sum_{\scr{|p|\le l}\atop\scr{p\in2\B{Z}}}{1\over2\pi
i}\left({\omega_1\omega_2 \over|\omega_1\omega_2|}\right)^p
\int_{(0)}{\sigma_\nu(\omega_1,-p/2)\sigma_\nu(\omega_2,-p/2)
\over|\omega_1\omega_2|^\nu|\zeta_\r{F}(1+\nu,p/2)|^2}
\lambda_l(\nu,p)\eta(\nu,p)\overline{\theta(\nu,p)}\,d\nu\cr
&={\delta_{\omega_1,\omega_2}+ \delta_{\omega_1,-\omega_2}\over4\pi^3
i}\sum_{|p|\le l}\int_{(0)}\lambda_l(\nu,p)\eta(\nu,p)
\overline{\theta(\nu,p)}(p^2-\nu^2)d\nu\cr
&+\sum_{c\ne0}{S_\r{F}(\omega_1,\omega_2;c) \over8\pi
i|c|^2}\sum_{|p|\le l}\int_{(0)}\e{K}_{\nu,p}\left({2\pi\over
c}\sqrt{\omega_1\omega_2}\right)\lambda_l(\nu,p)\eta(\nu,p)
\overline{\theta(\nu,p)}(p^2-\nu^2)d\nu,
\cr
}
$$
where $V\cap L^2(\varGamma\backslash G)_{l,q}\ne\{0\}$. }
\smallskip
\noindent
{\csc Remark.} For the idea of Hecke see Section 2.2 of [30]. The
inner product of two Poincar\'e series is the basis of almost all
proofs of the sum formula. In Kuznetsov's original proof in [22], the
seed function is explicit; and the same is in [32], where the
$K$-trivial case is treated. A more general class of seed functions
is used in [1] for $\r{PSL}_2(\B{R})$, and in [26] for the
$K$-trivial case on Lie groups of real rank one. Our discussion in
this section is different from that of Miatello and Wallach [26] in
that we positively exploit the arithmetical situation. Any
non-trivial estimate of Kloosterman sums suffices for the
continuation of $[\e{B}_\omega]\varphi_{l,q}(\nu,p)$ to a
neighbourhood of $\nu=1$, as has been indicated already. In the
general situation considered by Miatello and Wallach a spectral
decomposition is needed for analytic continuation. It gives in fact a
meromorphic continuation to $\B{C}$. In this respect our argument is
specific.
\par
It seems possible that the Poincar\'e series
$[\e{B}_\omega]\varphi_{l,q}(\nu,p)$ is relevant to the automorphic
resolvent of $\Omega_\pm$. In fact, the formula $(9.2)$ reminds us a
similar result for $\r{PSL}_2(\B{R})$, which is related to the
Fourier expansion of the automorphic resolvent for the Casimir
operator $\Omega$ (see [11]). We hope to return to this point
elsewhere.
\medskip
\centerline{\bf 10. Sum formula. I}
\smallskip
\noindent
Based on the above discussion, we shall establish the first version of
our sum formula, in which a given bilinear sum of Hecke eigenvalues,
or equivalently Fourier coefficients, of cuspidal irreducible
subspaces of $L^2(\varGamma\backslash G)$ is expressed in terms of
the arithmetic sums $S_\r{F}\;$:
\smallskip
\noindent
{\bf Theorem 10.1. (Spectral--Kloosterman sum formula)}\quad{\it Let
$h(\nu,p)$ be a function defined on a set $\{ \nu\in \B C\;:\;
|\Re\nu|\le {1\over2}+a\}\times\B{Z}$ for some small $a>0$,
satisfying the following conditions:
\item{1.} $h(\nu,p)=h(-\nu,-p)$,
\item{2.} $h(\nu,p)$ is regular,
\item{3.} $h(\nu,p)\ll (1+|\nu|+|p|)^{-4-b}$ with a small $b>0$.
\par
\noindent
Then we have, for any non-zero $\omega_1,\omega_2\in\B{Z}[i]$,
$$
\leqalignno{
&\sum_V|c_V(1)|^2t_V(\omega_1)t_V(\omega_2) h(\nu_V,p_V)&(10.1)\cr
&\hskip 1cm+\sum_{p\in2\B{Z}}{1\over2\pi i}\left({\omega_1\omega_2
\over|\omega_1\omega_2|}\right)^p
\int_{(0)}{\sigma_\nu(\omega_1,-p/2)\sigma_\nu(\omega_2,-p/2)
\over|\omega_1\omega_2|^\nu|\zeta_\r{F}(1+\nu,p/2)|^2}
h(\nu,p)\,d\nu\cr
&={\delta_{\omega_1,\omega_2}+ \delta_{\omega_1,-\omega_2}\over4\pi^3
i}\sum_{p\in\B{Z}}\int_{(0)}h(\nu,p)(p^2-\nu^2)d\nu\cr
&\hskip 1cm+\sum_{c\ne0}{1\over|c|^2}S_\r{F}(\omega_1,\omega_2;c)
\r{B}h\left({2\pi\over c}\sqrt{\omega_1\omega_2}\right).
\cr
}
$$
Here $V$ runs over all Hecke invariant right-irreducible cuspidal
subspaces of $L^2(\varGamma\backslash G)$ together with the
specifications in Section 8; and
$$
\r{B}h(u)=\sum_{p\in\B{Z}}{1\over8\pi i}
\int_{(0)}\e{K}_{\nu,p}(u)h(\nu,p)(p^2-\nu^2)d\nu\leqno(10.2)
$$
with $\e{K}_{\nu,p}$ as in $(7.21)$. Convergence of these expressions
is absolute throughout. }
\smallskip
\noindent
Proof. We denote the left and right sides of $(10.1)$ by
$L_{\omega_1,\omega_2}h$ and $R_{\omega_1,\omega_2}h$, respectively.
We may regard $L_{\omega_1,\omega_2}$, $R_{\omega_1,\omega_2}$ as
linear functionals on the space of functions defined on
$i\B{R}\times\B{Z}$. The eigenvalues of Hecke operators $t_V(\omega)$
are real, and so are the quantities $\bigl( {\omega/|\omega|}
\bigr)^p \sigma_\nu(\omega,-p/2)|\omega|^{-\nu}$. Thus
$L_{\omega,\omega}$ is positive definite for any non-zero
$\omega\in\B{Z}[i]$. We put
$$
\eta_0(\nu,p)=(1-\nu^2)^2(4-\nu^2)^{-2}
(4-\nu^2+p^2)^{-2-b/2}.\leqno(10.3)
$$
In $(9.30)$ we may set $\eta(\nu,p)=\lambda_l(\nu,p)^{-1}\eta_0(\nu,p)
e^{\delta\nu^2}$, and $\theta(\nu,p)=e^{\delta\nu^2}$ with
$\delta>0$. Thus we have
$$
L_{\omega_1,\omega_2}\eta_\delta^l=R_{\omega_1,\omega_2}
\eta_\delta^l,\leqno(10.4)
$$
where $\eta^l_\delta(\nu,p)=\eta_0(\nu,p)e^{2\delta\nu^2}$, if $|p|\le
l$, and $=0$ otherwise. Using this we are going to show that
$$
\lim_{l\to\infty}\,\lim_{\delta\to0+}L_{\omega_1,\omega_2}
\eta_\delta^l=L_{\omega_1,\omega_2} \eta_0.\leqno(10.5)
$$
Since $|L_{\omega_1,\omega_2}\eta_\delta^l
|\le\left(L_{\omega_1,\omega_1}\eta_\delta^l\cdot
L_{\omega_2,\omega_2}\eta_\delta^l\right)^{1/2}$, and $\eta_\delta^l$
is increasing on $i\B{R}\times\B{Z}$ as $\delta\downarrow0$,
$l\uparrow\infty$, it is obviously sufficient to show that
$L_{\omega,\omega}\eta_\delta^l$ is uniformly bounded for $l\ge0$,
$\delta\ge0$. To this end we shall prove that uniformly for $l\ge0$,
$\delta\ge0$
$$
\r{B}\eta^l_\delta(u)\ll |u|^{1+\varepsilon}\leqno(10.6)
$$
as $|u|\downarrow0$. Here $\varepsilon>0$ is small and may depend on
$b$. Then, by the bound $(8.14)$, one may confirm our claim by
showing that $\lim_{l\to\infty}\lim_{\delta\to0+}
R_{\omega,\omega}\eta^l_\delta=R_{\omega,\omega}\eta_0$.
\par
By definition we have
$$
\leqalignno{\qquad
&\r{B}\eta^l_\delta(u)&(10.7)\cr
&=\sum_{|p|\le l}{1\over8\pi i}
\int_{(0)}\e{K}_{\nu,p}(u)\eta_0(\nu,p)
e^{2\delta\nu^2}(p^2-\nu^2)d\nu\cr
&=-\sum_{|p|\le l}{1\over4\pi i}
\int_{(\alpha)}{\e{J}_{\nu,p}(u)\over\sin\pi\nu}\eta_0(\nu,p)
e^{2\delta\nu^2}(p^2-\nu^2)d\nu+{1\over4\pi}\sum_{|p|\le l}p^2
\e{J}_{0,p}(u)\eta_0(0,p), }
$$
where ${1\over2}<\alpha<1$. Since we have
$$
\e{J}_{0,p}(u)=(-1)^p|J_{|p|}(u)|^2\ll (|u|/2)^{2|p|}/
(|p|!)^2,\leqno(10.8)
$$
the last sum of $(10.7)$ is negligible compared with $(10.6)$. Also,
we have, as $|u|\downarrow0$, $\Re\nu=\alpha$,
$$
\leqalignno{ {\e{J}_{\nu,p}(u)\over\sin\pi\nu}&\ll
|u|^{2\alpha}{|\Gamma(|p|-\nu)|\over|\Gamma(|p|+1+\nu)|}&(10.9)\cr
&=|u|^{2\alpha}{|\Gamma(-\nu)|\over|(|p|+\nu)\Gamma(\nu)|}
\prod_{j=0}^{|p|-1}{|j-\nu|\over |j+\nu|}\ll
{|u|^{2\alpha}\over||p|+\nu||\nu|^{2\alpha}}.\cr
}
$$
Inserting this into $(10.7)$ we indeed get $(10.6)$.
\par
Next, we put
$$
h_\delta(\nu,p)=-ie^{\delta\nu^2}\sqrt{\delta\over\pi}
\int_{(0)}{h(\xi,p)\over\eta_0(\xi,p)}e^{\delta(\xi-\nu)^2}
d\xi.\leqno(10.10)
$$
We repeat the above discussion with $\theta=h_\delta$ and the same
$\eta$. We have first
$L_{\omega_1,\omega_2}f^l_\delta=R_{\omega_1,\omega_2} f^l_\delta$,
where $f^l_\delta(\nu,p)=\eta_0(\nu,p)h_\delta(\nu,p)
e^{\delta\nu^2}$ for $|p|\le l$, and $=0$ otherwise. Then we note
that $f^l_\delta(\nu,p)\ll \eta_0(\nu,p)$ uniformly for $\delta\ge0$,
$l\ge0$, and that $f_\delta^l\to h$ on $i\B{R}\times\B{Z}$ as
$\delta\downarrow0$, $l\uparrow\infty$. Thus, by $(10.5)$, we get
$$
\lim_{l\to\infty}\lim_{\delta\to0+}L_{\omega_1,\omega_2}
f_\delta^l=L_{\omega_1,\omega_2}h.\leqno(10.11)
$$
Corresponding to $(10.6)$ we have to estimate $\r{B}f^l_\delta(u)$. In
$(10.7)$ we replace $\eta_\delta^l$ by $f_\delta^l$ and set
$\alpha={1\over2}+a$ with $a>0$ given in the condition {\it 2} above.
Accordingly, we shift the contour in $(10.10)$ to $(\alpha)$, and see
that $f_\delta^l(\nu,p)\ll |\eta_0(\nu,p)|$ uniformly for $l\ge0$,
$\delta\ge0$ with $\Re\nu=\alpha$ and arbitrary integer $p$. Hence we
have the counterpart of $(10.6)$ for $f_\delta^l$. This gives
$$
\lim_{l\to\infty}\lim_{\delta\to0+}R_{\omega_1,\omega_2}
f_\delta^l=R_{\omega_1,\omega_2}h,\leqno(10.12)
$$
which ends the proof.
\smallskip
As the first application of the sum formula $(10.1)$ we shall prove
\smallskip
\noindent
{\bf Corollary 10.1.}\quad{\it There exist infinitely many $V$'s, and
we have, uniformly for $N,\,P\ge1$ and non-zero $\omega\in\B{Z}[i]$,
$$
\sum_{|\nu_V|\le N,\,|p_V|\le P}|c_V(1)|^2t_V(\omega)^2 \ll
(NP+|\omega|^{1+\varepsilon})(N^2+P^2)\leqno(10.13)
$$
with any fixed $\varepsilon>0$. In particular, we have the bound
$(8.15)$. }
\smallskip
\noindent
Proof. The deduction of the last assertion from $(10.13)$ is analogous
to the case of $\r{PSL}_2(\B{Z})$; see the proof of Lemma 3.3 of
[30]. To prove the first and the second assertions we put in $(10.1)$
$\omega_1=\omega_2=\omega$ and
$$
h(\nu,p)=h(\nu,p;N,P) =\exp((\nu/N)^2-(p/P)^2).\leqno(10.14)
$$
We note that $\zeta_\r{F}(1+\nu,p)\gg \log^{-1}(|\nu|+|p|+2)$ for
$\Re\nu=0$, with the implied constant being absolute. This can be
proved as in Sections 3.10--3.11 of [36]. The necessary uniform upper
bound for $\zeta_\r{F}(s,p)$ in the critical strip follows from the
functional equation $(5.34)$ and the convexity argument of Phragm\'en
and Lindel\"of. Thus we have
$$
\leqalignno{
&\sum_V|c_V(1)|^2t_V(\omega)^2h(\nu_V,p_V;N,P)+O(NP\sigma_0(\omega)^2
\log^2(NP+2))&(10.15)\cr
&={1\over8\pi^{5/2}}\sum_{p\in\B{Z}}(2Np^2+N^3)e^{-(p/P)^2}
+\sum_{c\ne0}{1\over|c|^2}S_\r{F}(\omega,\omega;c)\r{B}
h(2\pi\omega/c). }
$$
We are going to show
$$
\r{B}h(2\pi\omega/c)\ll\min(1,|\omega/c|^2)(N^2+P^2). \leqno(10.16)
$$
This and the bound $(8.14)$ give $(10.13)$, as well as an asymptotic
formula for the first sum in $(10.15)$, which gives the first
assertion.
\par
By $(12.1)$ below we have
$$
\r{B}h(2\pi\omega/c)={N^3\over2\pi^{3/2}}\int_1^\infty
C_{\omega/c}(y;N,P)\exp(-(N\log y)^2){dy\over y},\leqno(10.17)
$$
where
$$
\leqalignno{ C_{\omega/c}(y;N,P)=& \sum_{p\in\B{Z}}(-1)^p
((p/N)^2-(N\log y)^2+\txt{1\over2})&(10.18) \cr
&\times J_{2p}\left(2\pi|\omega/c||x|\right) \exp(-(p/P)^2+2pi\psi) }
$$
with $x=|x|e^{i\psi}=ye^{i\vartheta}+(ye^{i\vartheta})^{-1}$,
$\vartheta=\arg(\omega/c)$. Using an integral representation for
$J_{2p}$ we have also
$$
\leqalignno{\qquad
C_{\omega/c}(y;N,P)&={1\over2\pi}\int_{-\pi}^\pi\exp(2\pi
i|\omega/c||x|\cos(\theta+\psi))&(10.19)
\cr
&\times\sum_{p\in\B{Z}}((p/N)^2-(N\log y)^2+\txt{1\over2})
\exp(-(p/P)^2-2pi\theta)d\theta\cr
&={P\over2\sqrt{\pi}}\int_{-\pi}^\pi\exp(2\pi
i|\omega/c||x|\cos(\theta+\psi))\cr
&\hskip-2cm\times\sum_{q\in\B{Z}}\left((P/N)^2(\txt{1\over2}
-(P(\theta+q\pi))^2)-(N\log y)^2+\txt{1\over2}\right)
\exp(-(P(\theta+q\pi))^2)d\theta, }
$$
where the last line is due to Poisson's sum formula. This gives
$$
C_{\omega/c}(y;N,P)\ll(P/N)^2+(N\log y)^2+1\leqno(10.20)
$$
uniformly for all parameters involved. Now, the case $|\omega/c|\ge1$
is settled by inserting $(10.20)$ into $(10.17)$. If $|\omega/c|<1$
then we divide the integral in $(10.17)$ at $y=\sqrt{|c/\omega|}$.
For the infinite integral thus obtained we use again $(10.20)$, and
see that the contribution is negligible compared with $(10.16)$. To
estimate the remaining part we use $(10.18)$ together with
$J_{2p}(a)= (1+O(a^2))(a/2)^{2p}/(2p)!$ for small $a>0$. We have, for
$1\le y\le \sqrt{|c/\omega|}$,
$$
C_{\omega/c}(y;N,P)=-(N\log y)^2+{1\over2}+
O\left(|\omega/c|^2x^2((N\log y)^2+1)\right).\leqno(10.21)
$$
The contribution of this error term to $(10.17)$ is
$\ll|\omega/c|^2N^2$. As to the main term, we note that
$$
\int_1^\infty \left(-(N\log y)^2+\txt{1\over2}\right)\exp(-(N\log
y)^2) {dy\over y}=0.\leqno(10.22)
$$
Thus the relevant contribution to $(10.17)$ is easily seen to be
negligible. This ends the proof.
\smallskip
\noindent
{\csc Remark.} The class of test functions in Theorem 10.1 is as large
as possible. The strip on which the test functions are required to be
defined is narrow, due to the Weil bound $(8.14)$ (cf.\ $(3.6.24)$ of
[30]). The use of general seed functions in Lemma 9.1 leads to an
extension step with a functional analytic flavour. The proof is
similar to those in [1] (for $\r{PSL}_2(\B{R}))$, and [3] (for
$\r{PSL}_2$ over the product of the archimedean completions of a
number field). The proof in [32] for the $K$-trivial case is an
extension of Kuznetsov's original treatment [22] of the rational
case.
\par
The corollary is a counterpart of Kuznetsov's estimate for the
spectral mean square of the Fourier coefficients of Maass forms over
the modular group; see Lemma 2.4 of [30]. In Lemma 11 of [32] the
$K$-trivial case of the corollary is given. The bound $(10.13)$ is
essentially the best possible. We could prove an asymptotic result in
which the main term is a constant multiple of $NP(N^2+P^2)$. Note
that the proof of the corollary requires the rather deep integral
representation of $\e{K}_{\nu,p}$ in Theorem 12.1, whereas the series
expansion defining $\e{J}_{\nu,p}$ suffices for the theorem. One may
also consider the spectral large sieve estimate for
$$
\sum_{|\nu_V|\le N,\,|p_V|\le P}|c_V(1)|^2\Big|\sum_{\omega\ne0}
a(\omega)t_V(\omega)\Big|^2\leqno(10.23)
$$
with an arbitrary finite vector $\{a(\omega)\}$. To this we shall
return elsewhere, entailing a fuller treatment of the corollary. For
the rational case see [16] and Section 3.5 of [30].
\par
Generalization of the spectral sum formula and the bound $(10.13)$ to
other imaginary quadratic number fields and congruence subgroups
seems possible, as long as we have a counterpart of $(8.14)$. Without
such a bound, one has to be content with test functions which are
holomorphic for $|\Re\nu|\le1+a$ with an $a>0$, and have prescribed
zeros at $\nu=\pm 1$ for most values of $p$ (see $(9.18)$ and Lemma
9.1). In contrast to what we have seen in the above there might be,
in general, irreducible subspaces $V$ of complementary series type as
well, corresponding to exceptional eigenvalues. We note also that the
assumption that the spaces $V$ are Hecke invariant is not essential
for the sum formula; that is, the sums over $V$ in $(10.1)$ and
$(13.1)$ could be formulated in terms of Fourier coefficients in
place of Hecke eigenvalues.
\medskip
\centerline{\bf 11. A Bessel inversion}
\smallskip
\noindent
The aim of this section is to demonstrate a one-sided inversion of the
transform $\r{B}$ defined by $(10.2)$. Results of the present and the
next sections will play basic r\^oles in the proof of the second
version of our sum formula for $S_\r{F}$, which is to be developed in
Section 13.
\smallskip
\noindent
{\bf Theorem 11.1.}\quad{\it We put
$$
\r{K}f(\nu,p)=\int_{\B{C}^\times}\e{K}_{\nu,p}(u)f(u)d^\times\!u.
\leqno(11.1)
$$
Then, for any $f$ that is even, smooth and compactly supported on
$\B{C}^\times$, we have}
$$
2\pi\r{B}\r{K}f=f.\leqno(11.2)
$$
Proof. We shall prove, instead, the Parseval identity
$$
\int_{\B{C}^\times}f(u)g(u)d^\times\!u=\sum_{p\in\B{Z}}{1\over4i}
\int_{(0)}\r{K}f(\nu,p) \r{K}g(\nu,p)(p^2-\nu^2)d\nu,\leqno(11.3)
$$
where $f$, $g$ are to satisfy the condition given in the theorem. This
implies $(11.2)$, since a simple manipulation shows that the right
side is equal to
$$
\int_{\B{C}^\times}\r{B}\r{K}f(u)\cdot g(u)d^\times\!u.\leqno(11.4)
$$
Here the necessary absolute convergence follows from the estimate
$$
\r{K}f(\nu,p)\ll (1+|\nu|+|p|)^{-A},\leqno(11.5)
$$
where $A>0$ is arbitrary, and the implied constant depends on
$\Re\nu$, $A$, and the support of $f$. To show this we put
$$
\r{J}f(\nu,p)=\int_{\B{C}^\times}f(u)\e{J}_{\nu,p}(u)d^\times\!u,
\leqno(11.6)
$$
so that
$$
\r{K}f(\nu,p)={1\over\sin\pi\nu}\left\{\r{J}f(-\nu,-p)-
\r{J}f(\nu,p)\right\}.\leqno(11.7)
$$
By definition we have
$$
\r{J}f(\nu,p)=2\pi\sum_{m,\,n\ge0}{(-1)^{m+n}2^{-2(\nu+m+n)}
\r{M}f(\nu+m+n,p-m+n)\over m!n!\Gamma(\nu-p+m+1)\Gamma(\nu+p+n+1)}
,\leqno(11.8)
$$
where
$$
\r{M}f(\nu,p)=\int_0^\infty f_p(r)r^{2\nu-1}dr,\quad
f_p(r)={1\over2\pi}\int_0^{2\pi}f(re^{i\theta})e^{-2pi\theta}
d\theta. \leqno(11.9)
$$
A multiple application of partial integration gives, for any integers
$A,\,B\ge0$,
$$
f_p(r)\ll {r_f^A\over(1+|p|)^A},\quad \r{M}f(\nu,p)\ll
{1\over(1+|p|)^A}\left|{\Gamma(2\nu)\over\Gamma(2\nu+B)}\right|
r_f^{2|\Re\nu|+A+B},\leqno(11.10)
$$
where the constant $r_f$ depends on the compact support of $f$, and
the implied constants only on $A,\,B$ and $f$. Collecting these, we
get $(11.5)$. Thus we see also that the right side of $(11.3)$ is
equal to
$$
\leqalignno{ \lim_{T\to\infty}&\sum_{|p|\le P}{1\over4 i}
\int_{-iT}^{iT}\r{K}f(\nu,p) \r{K}g(\nu,p)(p^2-\nu^2)d\nu&(11.11)\cr
=&\lim_{T\to\infty} \left\{\int_{|u|<|v|}+\int_{|u|>|v|}\right\}
f(u)g(v)R_{P,T}(u,v)d^\times\!u d^\times\!v,\cr
}
$$
where $P=P(T)\in\B{Z}$ is to be chosen later, and
$$
R_{P,T}(u,v)=\sum_{|p|\le P}{1\over 4i}
\int_{-iT}^{iT}\e{K}_{\nu,p}(u) \e{K}_{\nu,p}(v)(p^2-\nu^2)d\nu.
\leqno(11.12)
$$
\par
We shall consider the case $|u|<|v|$. We indent the contour
$[-iT,\,iT]$ with the right half of a small circle centered at
$\nu=0$. Denoting the new contour by $L_T$, we have, by $(11.7)$,
$$
R_{P,T}(u,v)=\sum_{|p|\le P}{1\over 4i}\left\{\int_{L_T^-}-\int_{L_T}
\right\}\e{J}_{\nu,p}(u)
\e{K}_{\nu,p}(v){p^2-\nu^2\over\sin\pi\nu}d\nu,\leqno(11.13)
$$
where $L^-_T=\{\nu:-\nu\in L_T\}$. This implies that
$$
\leqalignno{ R_{P,T}(u,v)&=\sum_{|p|\le P}{i\over 2}
\int_{C_{P,T}}\e{J}_{\nu,p}(u)
\e{K}_{\nu,p}(v){p^2-\nu^2\over\sin\pi\nu}d\nu&(11.14)\cr
&+{1\over2}\sum_{|p|\le P}\sum_{0\le q\le P}(-1)^qa_q
\e{J}_{q,p}(u)\e{K}_{q,p}(v)(p^2-q^2), }
$$
where $a_0=1$, $a_q=2$ for $q>0$, and $C_{P,T}$ is the oriented
polygonal line connecting the points $-iT$, $P+{1\over2}-iT$,
$P+{1\over2}+iT$, $iT$ in this order. This double sum vanishes,
because of the identities, for any $p,q\in\B{Z}$,
$$
\e{J}_{q,p}=(-1)^{p+q}\e{J}_{p,q}
=(-1)^{p+q}\e{J}_{|p|,q\,\r{sign}(p)},\quad
\e{K}_{q,p}=\e{K}_{|p|,q\,\r{sign}(p)}.\leqno(11.15)
$$
In fact, the first identity is trivial, and the second follows from
the expression
$$
\e{K}_{q,p}(u)={(-1)^{q+1}\over\pi}
\left\{\b{Y}_{q-p}(u)J_{p+q}(\bar{u})
+J_{p-q}(u)\b{Y}_{p+q}(\bar{u})\right\}.\leqno(11.16)
$$
Here $\b{Y}_n$ with $n\in\B{Z}$ is the Hankel function, which
satisfies $\b{Y}_n=(-1)^n\b{Y}_{-n}$ (see p.\ 59 of [42]). Hence we
have
$$
\leqalignno{\qquad R_{P,T}(u,v)&=\sum_{|p|\le P}{i\over2}
\int_{C_{P,T}}\e{J}_{\nu,p}(u)
\e{K}_{\nu,p}(v){p^2-\nu^2\over\sin\pi\nu}d\nu&(11.17)\cr
&=\sum_{|p|\le P}{i\over2} \int_{C_{P,T}}\left[\e{J}_{\nu,p}(u)
\e{J}_{-\nu,-p}(v)-\e{J}_{\nu,p}(u)\e{J}_{\nu,p}(v)\right]
{p^2-\nu^2\over(\sin\pi\nu)^2}d\nu\cr
&=R^{-}_{P,T}(u,v)-R^{+}_{P,T}(u,v), }
$$
say, in an obvious mode of division.
\par
Now we have trivially
$$
\int_{|u|<|v|}f(u)g(v)R^\pm_{P,T}(u,v)d^\times\!u\,d^\times\!v
=\int_{r_1<r_2}(r_1r_2)^{-1}Q^\pm_{P,T}(r_1,r_2;f,g)
dr_1\,dr_2,\leqno(11.18)
$$
where
$$
Q^\pm_{P,T}(r_1,r_2;f,g)
=\int_0^{2\pi}\int_0^{2\pi}f(r_1e^{i\theta_1})g(r_2e^{i\theta_2})
R^\pm_{P,T}(r_1e^{i\theta_1},r_2e^{i\theta_2}) d\theta_1\,d\theta_2.
\leqno(11.19)
$$
The series expansion of $\e{J}_{\nu,p}$ implies that
$$
\leqalignno{ \e{J}_{\nu,p}(u)&\e{J}_{-\nu,-p}(v)
{p^2-\nu^2\over(\sin\pi\nu)^2}=-\pi^{-2} \left|{u\over
v}\right|^{2\nu}\left({u\over|u|}\right)^{-2p}
\left({v\over|v|}\right)^{2p}&(11.20)\cr
&\times\sum_{k,l,m,n\ge0}{(-1)^{k+l+m+n}(u/2)^{2k}(\bar{u}/2)^{2l}
(v/2)^{2m}(\bar{v}/2)^{2n}\over k!l!m!n!\lambda(\nu,p;k,l,m,n)}, }
$$
where $\lambda(\nu,p;k,l,m,n)=
(\nu-p+1)_k(\nu+p+1)_l(-\nu+p+1)_m(-\nu-p+1)_n$ with $(\alpha)_k$ as
in $(5.19)$. Thus we have
$$
\leqalignno{ Q^-_{P,T}(r_1,r_2;f,g)=&\sum_{|p|\le P}\sum_{k,l,m,n\ge0}
{(-1)^{k+l+m+n}(r_1/2)^{2(k+l)} (r_2/2)^{2(m+n)}\over
k!l!m!n!}&(11.21)\cr
&\times f_{p-k+l}(r_1)g_{-p-m+n}(r_2)S_{P,T}^-(r_1/r_2,p; k,l,m,n), }
$$
where $g_q$ is analogous to $f_q$, and
$$
S_{P,T}^-(\rho,p;k,l,m,n) =-2i\int_{C_{P,T}} {\rho^{2\nu}\over
\lambda(\nu,p;k,l,m,n)}d\nu. \leqno(11.22)
$$
Similarly we have
$$
\leqalignno{ Q^+_{P,T}(r_1,r_2;f,g)=&\sum_{|p|\le
P}\;\sum_{k,l,m,n\ge0} {(-1)^{k+l+m+n}(r_1/2)^{2(k+l)}
(r_2/2)^{2(m+n)}\over k!l!m!n!}&(11.23)\cr
&\times f_{p-k+l}(r_1)g_{p-m+n}(r_2)S_{P,T}^+(r_1r_2,p; k,l,m,n), }
$$
where
$$
\leqalignno{\qquad S_{P,T}^+(\rho,p;k,l,m,n)
=2\pi^2i\int_{C_{P,T}}&{(p^2-\nu^2)(\rho/4)^{2\nu}\over
(\sin\pi\nu)^2\Gamma(\nu-p+k+1)\Gamma(\nu+p+l+1)}&(11.24)\cr
&\times{1\over\Gamma(\nu-p+m+1) \Gamma(\nu+p+n+1)}d\nu. }
$$
\par
Assuming that $2P\le T$, we shall estimate $Q^\pm_{P,T}(r_1,r_2;f,g)$;
implicit constants may depend only on the parameter $A$ and the
supports of $f$ and $g$. To this end we stress that the first
estimate in $(11.10)$ implies readily that for any $A>0$
$$
\sum_{k,l,m,n\ge0} {(r_1/2)^{2(k+l)} (r_2/2)^{2(m+n)}\over
k!l!m!n!}|f_{p-k+l}(r_1)g_{\pm p-m+n}(r_2)|\ll{1\over(1+|p|)^A}.
\leqno(11.25)
$$
Arguing as in $(10.9)$ we see that the integrand in $(11.24)$ is
$O((\rho/4)^{2\Re\nu}|\Gamma (-\nu)/\Gamma(\nu)|^2)$. Hence we get
immediately
$$
Q^+_{P,T}(r_1,r_2;f,g)\ll{1\over\log T}+{1\over P!},\leqno(11.26)
$$
where the terms on the right come from the horizontal and the vertical
parts of $C_{P,T}$, respectively. As to $Q^-_{P,T}(r_1,r_2;f,g)$ we
note first that
$$
S_{P,T}^-(\rho,p;0,0,0,0) =-2i\int_{-iT}^{iT}\rho^{2\nu}d\nu.
\leqno(11.27)
$$
If $(k,l,m,n)\ne(0,0,0,0)$, then on the horizontal part of $C_{P,T}$
we have $|\lambda (\nu,p;k,l,m,n)|\allowbreak\gg T$, and the
corresponding contribution in $(11.22)$ is $\ll\min(1,
(T|\log\rho|)^{-1})$, provided $\rho\le1$. We restrict ourselves to
the vertical part of $C_{P,T}$. We may assume naturally $p\ge0$.
Then, if either $l>0$ or $n>0$, we have $|\lambda (\nu,p;k,l,m,n)|\gg
P+|\nu|$. By partial integration we see that the corresponding
contribution to $(11.22)$ is
$$
\ll(k+l+m+n)\log T\min\left(1,{1\over P|\log\rho|}\right),
\leqno(11.28)
$$
provided $\rho\le1$. If $l=n=0$, and $0\le p\le P/2$, then obviously
we get the same conclusion. Otherwise the contribution in question is
$\ll\log T$, provided $\rho\le1$. Collecting these, we have, for
$r_1<r_2$, $2P\le T$,
$$
\leqalignno{ Q^-_{P,T}(r_1,r_2;&f,g)=-2i\sum_{|p|\le
P}f_p(r_1)g_{-p}(r_2)\int_{-iT}^{iT}(r_1/r_2)^{2\nu}d\nu&(11.29)\cr
&+O\left(\log T\min\left(1,{1\over
P|\log(r_1/r_2)|}\right)\right)+O\left({\log T\over P^A}\right), }
$$
which ends the discussion of the case $|u|<|v|$.
\par
The case $|u|>|v|$, i.e., $r_1>r_2$, can be treated in just the same
way. We return to $(11.12)$, and this time we shift the relevant
contours to the left, getting the same assertions as $(11.26)$ and
$(11.29)$. In this way we now have
$$
\leqalignno{
&\int_{\B{C}^\times\times\B{C}^\times}
f(u)g(v)R_{P,T}(u,v)d_*u\,d_*v&(11.30)
\cr
&=-2i\sum_{|p|\le P}\int_{-iT}^{iT} \r{M}f(\nu,p)\r{M}g(-\nu,-p)d\nu
+O\left({1\over\log T}+{\log^2T\over P}\right). }
$$
Hence we set $P=[\log^3T]$. We find that the right side of $(11.3)$ is
equal to
$$
-2i\sum_{p\in\B{Z}}\int_{(0)}
\r{M}f(\nu,p)\r{M}g(-\nu,-p)d\nu.\leqno(11.31)
$$
With this and the Parseval formulas for Mellin transform and Fourier
series expansion, we finish the proof.
\smallskip
We shall also need the following property of the transform $\r{K}$:
\smallskip
\noindent
{\bf Lemma 11.1.}\quad{\it Let $f$ be an even smooth function on
$\B{C}^\times$ with a compact support. Then we have}
$$
\sum_{p\in\B{Z}}\int_{(0)}\r{K}f(\nu,p)(p^2-\nu^2)d\nu=0.
\leqno(11.32)
$$
Proof. Let $\nu\in\B{C}$ be such that $\Re\nu\ge0$ and $|\nu-\nobreak
k|\ge {1\over2}$ for all $k\in\B{Z}$. Then $(11.8)$ and $(11.10)$
give
$$
\leqalignno{
\r{J}f(\nu,p)\ll\sum_{m,\,n\ge0}&{(r_f/2)^{2\Re\nu+2(m+n)} \over
m!n!|\Gamma(\nu-p+m+1)||\Gamma(\nu+p+n+1)|}&(11.33) \cr
&\times{r_f^{A+B}\over(1+|p-m+n|)^A|\nu+m+n|^B}. }
$$
If $|p-m+n|\le{1\over2}|p|$, then we have $n+m\ge{1\over2}|p|$, and
consequently $|\nu+m+n|\gg|\nu|+|p|$; otherwise we have
$|p-m+n||\nu+m+n|\gg |p\nu|$. Thus we have, for any fixed large
$C>0$,
$$
\r{J}f(\nu,p) \ll{(r_f/2)^{2\Re\nu} \over|\Gamma(\nu-|p|+1)
\Gamma(\nu+|p|+1)|}(|\nu|+|p|)^{-C}.\leqno(11.34)
$$
This estimate allows us to carry out the same procedure as in
$(11.13)$--$(11.14)$: We have, for any positive integer $P$,
$$
\sum_{p=-P}^P\int_{(0)}\r{K}f(\nu,p)(p^2-\nu^2)d\nu=
2\sum_{p=-P}^P\int_{(P+1/2)}\r{J}f(\nu,p)
{p^2-\nu^2\over\sin\pi\nu}d\nu. \leqno(11.35)
$$
As before, the sum of residues arising from this shift of contour
vanishes because of the first relation in $(11.15)$. The last
integral is, in view of $(10.9)$ and $(11.34)$,
$$
\ll\left({r_f\over2P}\right)^{2P}\int_{-\infty}^\infty
(P+|t|)^{2-C}dt\leqno(11.36)
$$
uniformly for $|p|\le P$. This obviously ends the proof.
\smallskip
\noindent
{\csc Remark.} The inversion formula $(11.2)$ could be formulated as a
discontinuous integral of new type in the theory of Bessel functions.
The idea of the proof is to view the transformation $\r{K}$ in (11.1)
as a perturbation of the Mellin--Fourier transformation on
$\B{C}^\times$. Insert the power series expansion of
$\e{J}_{\pm\nu,\pm p}$ into the integrals hidden in $\r{K}f(\nu,p)
\r{K}g(\nu,p)$ on the right of (11.3). Two of the four lowest order
terms describe the Mellin--Fourier transformation on $\B{C}^\times$
in polar coordinates, as is well indicated by $(11.20)$. The proof of
the inversion consists of showing that all other terms do not
contribute. The key to achieve this is the vanishing of the double
sum in $(11.14)$. That is, a certain rearrangement of products of
$J$-Bessel functions of various orders is taking place behind our
argument, which further points to a relation with the Neumann
expansion (see Chapter XVI of [42]). The basic idea is present in
Section 2.5 of Kuznetsov's preprint [21], which deals with the Bessel
inversion for the modular case, and is indeed the first instance of
such investigations (see also Section 2.4 of [30]).
\par
Our proof is, however, admittedly technical, and one may wish to find
a more structural proof that takes into account the way through which
the functions $\e{J}_{\nu,p}$ and $\e{K}_{\nu,p}$ come into our
discussion. They correspond to functions on the big cell in the
Bruhat decomposition of $G$, transforming on the left and the right
according to non-trivial characters of the subgroup $N$, and turn out
to be a basis of the solutions of $\Omega_\pm f = {1\over8}( (\nu\mp
p)^2-1)f$. These and the adjoint formulation $(11.3)$ suggest that
the inversion should be a part of the spectral theory on the big
cell. A proof along such a line might work for other Lie groups of
rank one as well. To this we hope to return elsewhere.
\par
As to Lemma 11.1, we remark that there are test functions $h$ such as
the one introduced in $(10.14)$, for which
$$
\sum_{p\in\B{Z}}\int_{(0)}h(\nu,p)(p^2-\nu^2)d\nu\ne0.\leqno(11.37)
$$
Thus Lemma 11.1 shows that $(11.2)$ gives only a one-sided inversion
of the transformation $\r{B}$. Also see the remark at the end of
Section 13 for an alternative argument .
\medskip
\centerline{\bf 12. The Bessel kernel $\e{K}_{\nu,p}$}
\smallskip
\noindent
The main feature of Theorem 10.1 rests precisely in the integral
transform $\r{B}$ defined by $(10.2)$, and thus in the kernel
$\e{K}_{\nu,p}$. In this section we shall prove an integral formula
for $\e{K}_{\nu,p}$, which has a practical value for our purpose, and
an interest of its own.
\smallskip
\noindent
{\bf Theorem 12.1.}\quad{\it Let $|\Re\nu|<{1\over4}$. Then we have,
for any $p\in\B{Z}$ and non-zero $u\in\B{C}$,
$$
\leqalignno{ \e{K}_{\nu,p}(u)=(-1)^p{2\over\pi}\int_0^\infty
&y^{2\nu-1}\left({ye^{i\vartheta}+(ye^{i\vartheta})^{-1}
\over|ye^{i\vartheta}+(ye^{i\vartheta})^{-1}|}\right)^{2p}&(12.1)
\cr
&\times J_{2p}\left(|u||ye^{i\vartheta}+(ye^{i\vartheta})^{-1}|
\right)dy,\qquad }
$$
where $u=|u|e^{i\vartheta}$. }
\smallskip
\noindent
Proof. Since $\overline{\e{K}_{\bar{\nu},p}(u)}=\e{K}_{\nu,-p}(u)$, we
may assume that $p$ is non-negative. We shall first show that we
have, for $\Re\nu>-p$,
$$
\leqalignno{ \e{J}_{\nu,p}(2\pi u)&=(-1)^p|u/2|^{2p}\sum_{|m|\le p}
{2p\choose p+m}(iu/|u|)^{-2m}&(12.2)
\cr
&\times{1\over2\pi i}\int_{(1)}\exp(2\pi\xi-\pi\Re(u^2)/\xi)
I_{\nu-m}(\pi|u|^2/\xi)\xi^{-2p-1}d\xi.
\cr
}
$$
The basis for this formula is $(6.20)$, where the function
$\e{J}_{\nu,p}$ enters into our investigations. We set there
$\omega_1=1$, $\omega_2=u^2$, $l=p$, and $q=p$. On noting $(6.18)$,
we have, for $\Re\nu>0$, $r>0$, and $\r{k}\in K$,
$$
\leqalignno{ \int_\B{C}e^{-2\pi i\Re
z}&\e{B}_1\varphi_{p,p}(\nu,p)(\r{h}[u]\r{w}
\r{n}[z]\r{a}[r]\r{k})d_+z&(12.3)\cr
&=\pi^{-2\nu}|u|^2\e{J}_{\nu,p}(2\pi u)
\e{A}_1\varphi_{p,p}(\nu,p)(\r{a}[r]\r{k}),
\cr
}
$$
where, by $(5.9)$,
$$
\leqalignno{
&\r{h}[u]\r{w}\r{n}[z]\r{a}[r]&(12.4)\cr
&= \r{n}\left[{-u^2\bar{z}\over r^2+|z|^2}\right]\r{a}\left[
{|u|^2r\over r^2+|z|^2}\right]\r{h}\left[{u\over|u|}\right]
\r{k}\left[{\bar{z}\over\sqrt{r^2+|z|^2}},
{-r\over\sqrt{r^2+|z|^2}}\right]. }
$$
Using $(3.24)$, $(5.26)$--$(5.27)$, and $(6.13)$--$(6.14)$, we equate
the coefficients of $\Phi_{0,p}^p(\r{k})$ on both sides of $(12.3)$,
getting
$$
K_\nu(2\pi r)\e{J}_{\nu,p}(2\pi u)={|u|^{2p}r^{-2p}\over2\pi}
\sum_{|m|\le p} i^m(u/|u|)^{-2m}U_{\nu,m}(u,r).\leqno(12.5)
$$
Here
$$
\leqalignno{ U_{\nu,m}(u,r)=&\int_\B{C} \exp\Big(-2\pi ir\Re z-2\pi
i{\Re u^2\bar{z}\over r(1+|z|^2)} \Big)&(12.6)\cr
&\times I_{\nu-m}\Big({2\pi|u|^2\over r(
1+|z|^2)}\Big)\Phi_{m,0}^{p,*}\left(\left[{\bar{z}\atop1}\; {-1\atop
z}\right]\right){d_+z\over (1+|z|^2)^{2p+1}}, }
$$
where the asterisk denotes that we have extended the definition
$(3.18)$ in an obvious way. We are going to compute $U_{\nu,p}(u,r)$
asymptotically when $r$ tends to infinity, so that the result yields
the cancellation of the factor $K_\nu(2\pi r)\sim {1\over2}
r^{-1/2}e^{-2\pi r}$ on the left side of $(12.5)$. We put $z=x+iy$
with $x,y\in\B{R}$ in $(12.6)$, and then regard the integral as a
double complex integral. Studying partial derivatives of the argument
of the exponentiated factor, we see that a saddle point exists at
$x=-i+c_1/r$, $y=c_2/\sqrt{r}$, where $c_1$, $c_2$ are asymptotically
constant as $r\uparrow\infty$. Because of this we make the change of
variables $(x,\,y)\mapsto(-i(1-1/r)+x/r,\, y/\sqrt{r})$. Here the
factor $1-1/r$ is to avoid the singularity at $(-i,0)$. We have
$$
\leqalignno{
&U_{\nu,m}(u,r)=r^{2p-{1\over2}}e^{-2\pi (r-1)}&(12.7)\cr
&\times\int_{-\infty}^\infty\int_{-\infty}^\infty \exp\Big(-2\pi
ix+2\pi i{ai-a(x+i)/r-by/\sqrt{r}\over 2+y^2-2ix+(x+i)^2/r}\Big)\cr
&\times\Phi_{m,0}^{p,*}\Big( \Big[{-i+(x+i)/r-iy/\sqrt{r}\atop
1}{-1\atop -i+(x+i)/r+iy/\sqrt{r}} \Big]\Big) \cr
&\times I_{\nu-m}\Big({2\pi|u|^2\over
2+y^2-2ix+(x+i)^2/r}\Big){dxdy\over (2+y^2-2ix+(x+i)^2/r)^{2p+1}}, }
$$
where $a=\Re u^2$, $b=\Im u^2$. It is easy to check that this double
integral converges absolutely and uniformly as $r\uparrow\infty$,
provided $\Re\nu>0$. Thus we have
$$
\leqalignno{ \lim_{r\to\infty}&{U_{\nu,p}(u,r) \over r^{2p}K_\nu(2\pi
r)}=2\int_{-\infty}^\infty\int_{-\infty}^\infty
\exp\Big(2\pi(1-ix)-{2\pi a\over2+y^2-2i x}\Big)&(12.8)\cr
&\times \Phi_{m,0}^{p,*}\Big( \Big[{-i\atop 1}\;{-1\atop -i}\Big]\Big)
I_{\nu-m}\Big({2\pi|u|^2\over2+y^2-2ix}\Big){dxdy\over
(2+y^2-2ix)^{2p+1}}. }
$$
We shift the contour of the $x$-integral to $\Im x=-y^2/2$. We then
find that
$$
\leqalignno{
&\lim_{r\to\infty}{U_{\nu,p}(u,r) \over r^{2p}K_\nu(2\pi r)}&(12.9)\cr
&=-i2^{-2p}\int_{(1)} \exp\Big(2\pi\xi-{\pi a\over
\xi}\Big)\Phi_{m,0}^{p,*}\Big( \Big[{-i\atop 1}\;{-1\atop
-i}\Big]\Big) I_{\nu-m}\Big({\pi|u|^2\over\xi}\Big){d\xi\over
\xi^{2p+1}}. }
$$
On noting that $\Phi_{m,0}^{p,*}\left( \left[{-i\atop 1}\;{-1\atop
-i}\right]\right) =(-1)^pi^m{2p\choose p-m}$, we get, from $(12.5)$
and $(12.9)$, the representation $(12.2)$ at least for $\Re\nu>0$,
and then we use analytic continuation with respect to $\nu$.
\par
We move to the proof of $(12.1)$. We shall treat first the case $p>0$.
Since $\e{J}_{-\nu,-p}(u)=\overline{\e{J}_{-\bar{\nu},p}(u)}$, we see
that $(12.2)$ gives, for $|\Re\nu|<p$,
$$
\leqalignno{ \e{K}_{\nu,p}(2\pi u)&=(-1)^p{2\over\pi}
\left|{u\over2}\right|^{2p}\sum_{|m|\le p} {2p\choose
p+m}\left({u\over|u|}\right)^{-2m}&(12.10)\cr
&\times{1\over2\pi
i}\int_{(1)}\exp\left(2\pi\xi-{\pi\Re(u^2)\over\xi}\right)
K_{\nu-m}\left({\pi|u|^2\over\xi}\right)\xi^{-2p-1}d\xi.
\cr
}
$$
Observing that $\Re\xi^{-1}>0$, we have
$$
\leqalignno{
&\sum_{|m|\le p} {2p\choose p+m}\left({u\over|u|}\right)^{-2m}
K_{\nu-m}\left({\pi|u|^2\over\xi}\right)&(12.11)
\cr
&={1\over2}\int_0^\infty \left(\sqrt{y}e^{i\vartheta}+
{1\over\sqrt{y}e^{i\vartheta}}\right)^{2p}
\exp\left(-{\pi|u|^2\over2\xi} \left(y+{1\over
y}\right)\right)y^{\nu-1}dy\cr
&=I_{\nu,p}(u,\xi)+I_{-\nu,p}(\bar{u},\xi), }
$$
say, where $\vartheta$ is as in $(12.1)$, and $I_{\nu,p}(u,\xi)$ is
the part corresponding to $y\ge1$. We put
$$
\e{K}^*_{\nu,p}(2\pi u)=(-1)^p{1\over\pi^2i}
\left|{u\over2}\right|^{2p}\int_{(1)}
\exp\left(2\pi\xi-{\pi\Re(u^2)\over\xi}\right)I_{\nu,p}(u,\xi)
\xi^{-2p-1}d\xi.\leqno(12.12)
$$
We insert into this the integral representation of $I_{\nu,p}(u,\xi)$.
The resulting double integral is absolutely convergent for
$\Re\nu<0$, and $\e{K}^*_{\nu,p}(2\pi u)$ is regular there. To get
analytic continuation to $\Re\nu\ge0$, we turn the line of
integration of $I_{\nu,p}(u,\xi)$ around the point $y=1$ through a
small angle which has the same sign as $\Im\xi$. We get immediately
the bound $I_{\nu,p}(u,\xi)\ll|\xi|^{\Re\nu+p}$, which means that
$\e{K}_{\nu,p}^*(2\pi u)$ is regular for $\Re\nu<p$. Thus we have the
decomposition
$$
\e{K}_{\nu,p}(2\pi u)=\e{K}^*_{\nu,p}(2\pi
u)+\e{K}^*_{-\nu,p}(2\pi\bar{u}), \leqno(12.13)
$$
provided $|\Re\nu|<p$. We now assume that $-p<\Re\nu<0$. Then, because
of the absolute convergence mentioned above, we may exchange the
order of integration in $(12.12)$. We get
$$
\leqalignno{\qquad \e{K}^*_{\nu,p}(2\pi u)&=(-1)^p{1\over2\pi^2i}
\left|{u\over2}\right|^{2p}\int_1^\infty y^{\nu-1}
\left(\sqrt{y}e^{i\vartheta}+
{1\over\sqrt{y}e^{i\vartheta}}\right)^{2p}&(12.14)
\cr
&\hskip 1cm\times\int_{(1)} \exp\left(2\pi\xi-{\pi|u|^2\over2\xi}
\left|\sqrt{y}e^{i\vartheta}+
{1\over\sqrt{y}e^{i\vartheta}}\right|^2\right) \xi^{-2p-1}d\xi\,dy\cr
&={(-1)^p\over\pi}\int_1^\infty y^{\nu-1}\left({\sqrt{y}e^{i\vartheta}
+(\sqrt{y}e^{i\vartheta})^{-1} \over|\sqrt{y}e^{i\vartheta}
+(\sqrt{y}e^{i\vartheta})^{-1}|}\right)^{2p}\cr
&\hskip 1cm\times J_{2p}\left(2\pi|u|\left|\sqrt{y}e^{i\vartheta}+
{1\over\sqrt{y}e^{i\vartheta}}\right|\right)dy. }
$$
By the asymptotic property of the $J$-Bessel function, the last
integral converges absolutely for $\Re\nu<{1\over4}$. This fact and
the identity $(12.13)$ gives rise to $(12.1)$, if $p\ne0$.
\par
We shall next consider the case $p=0$. We are unable to use the
formula $(12.2)$. Nonetheless, the right side of $(12.1)$ converges
absolutely for $p=0$ and $|\Re\nu|<{1\over4}$. By Neumann's addition
theorem for $J_0$ we see that it is equal to
$$
{2\over\pi}\sum_{m\in\B{Z}}(-1)^me^{2mi\vartheta}\int_0^\infty
y^{2\nu-1}J_m(|u|y)J_m(|u|/y)dy,\leqno(12.15)
$$
where the necessary absolute convergence is easy to check. By Lemma 6
of [33] this is transformed into
$$
\leqalignno{
{8\over\pi^2}\cos\pi\nu&\sum_{m\in\B{Z}}(-1)^me^{2mi\vartheta}
\int_0^{\pi/2}J_{2m}(|u|\sin\tau)K_{2\nu}(2|u|\cos\tau)d\tau
&(12.16)\cr
&={8\over\pi^2}\cos\pi\nu\int_0^{\pi/2}\cos(2|u|\cos\vartheta
\sin\tau)K_{2\nu}(2|u|\cos\tau)d\tau, }
$$
where the second line depends on the definition of the $J$-Bessel
function of an integral order. According to Lemma 8 of [32], the last
integral is equal to ${1\over8}\pi^2\e{K}_{\nu,0}(u)/\cos\pi\nu$.
This ends the proof.
\smallskip
Next, we turn to the mean of $\e{K}_{\nu,p}$ over $\B{C}^\times$. We
shall later require that it is not too large.
\smallskip
\noindent
{\bf Lemma 12.1.}\quad{\it Let
$$
2|\Re\nu|<\rho<{1\over2}\,.\leqno(12.17)
$$
Then we have
$$
\int_{\B{C}^\times}|\e{K}_{\nu,p}(u)|^2|u|^{2\rho}d^\times\!u\ll
(1+|p|)^{2\rho-1}, \leqno(12.18)
$$
where the implied constant depends only on $\rho$ and $\Re\nu$. }
\smallskip
\noindent
Proof. The formula $(12.1)$ and the Cauchy--Schwartz inequality give
$$
\leqalignno{
|\e{K}_{\nu,p}(re^{i\vartheta})|^2\le{4\over\pi^2}\int_1^\infty&
(y^{2\Re\nu}+y^{-2\Re\nu})^2y^{-\eta-1}dy&(12.19)\cr
&\times\int_1^\infty y^{\eta-1}J_{2|p|}
(r|ye^{i\vartheta}+(ye^{i\vartheta})^{-1}|)^2dy, }
$$
where we assume $4|\Re\nu|<\eta<2\rho$. We multiply both sides by
$r^{2\rho-1}$ and integrate with respect to $r$ and $\vartheta$,
getting
$$
\leqalignno{
&\int_{\B{C}^\times}|\e{K}_{\nu,p}(u)|^2|u|^{2\rho}d^\times\!u
&(12.20)\cr
&\ll\int_1^\infty \int_0^{2\pi}y^{\eta-1}
|ye^{i\vartheta}+(ye^{i\vartheta})^{-1}|^{-2\rho}d\vartheta\,dy
\int_0^\infty J_{2|p|}(r)^2r^{2\rho-1}dr. }
$$
The right side converges. Then we invoke
$$
\int_0^\infty J_{2|p|}(r)^2r^{2\rho-1}dr=2^{2\rho-1}
{\Gamma(1-2\rho)\Gamma(2|p|+\rho)\over
\Gamma(1-\rho)^2\Gamma(2|p|+1-\rho)}\leqno(12.21)
$$
(see p.\ 403 of [42]). This ends the proof.
\smallskip
\noindent
{\csc Remark.} The integral formula $(12.1)$ appears to be new;
despite its classical outlook we have not been able to find its
tabulation. An alternative proof of $(12.1)$ is indicated in Section
15.
\par
It seems worth remarking that $\e{K}_{\nu,p}$ appears in the context
of Section 2 as well: We consider, more generally than $(1.1)$, the
mean value
$$
\int_{-\infty}^\infty |\zeta_\r{F}(\txt{1\over2}+it)
\zeta_\r{F}(\txt{1\over2}+it, \txt{1\over2}p)|^2g(t)dt,\leqno(12.22)
$$
where $p\in2\B{Z}$. Then we need to treat
$$
\sum_{n(n+m)\ne0}\sigma_\alpha(n,\txt{1\over2}p)
\sigma_\beta(n+m,-\txt{1\over2}p) g^*(n/m;\gamma,\delta)\leqno(12.23)
$$
with $m\ne0$. Corresponding to $(2.31)$ we have the sum of Kloosterman
sums
$$
\sum_{c\ne0} {1\over {|c|^2}}S_\r{F}(m,n\,;c)\,[g]_p\Big({2\pi\over c}
\sqrt{mn};\alpha,\beta,\gamma,\delta\Big),\leqno(12.24)
$$
where
$$
\leqalignno{\qquad
&[g]_p(u;\alpha,\beta,\gamma,\delta) =(|u|/2)^{-2(1+\alpha+\beta)}
\sum_{q\in\B{Z}} (-1)^{\max(|p|,|q|)}(iu/|u|)^{2(p+q)}&(12.25)\cr
&\times\int_{(\eta)}{\Gamma(1-s+{1\over2}|p+q|)
\Gamma(1+\alpha-s+{1\over2}|p-q|)\over \Gamma(s+{1\over2}|p+q|)
\Gamma(s-\alpha+{1\over2}|p-q|)}
\tilde{g}_{p+q}(s;\gamma,\delta)(|u|/2)^{4s}ds.
\cr
}
$$
with $\eta<1+\min(0,\Re\alpha)$. We compare this with $(14.13)$, and
are led to the expression
$$
\leqalignno{\qquad
&[g]_p(u;\alpha,\beta,\gamma,\delta)&(12.26)\cr
&=\pi i(u/|u|)^{2p}
(|u|/2)^{2(1-\beta)}\int_{\B{C}^\times}g^*(v;\gamma,\delta)
\e{K}_{-\alpha,p}(iu\sqrt{v})(v/|v|)^p |v|^{2+\alpha}d^\times\!v. }
$$
This belongs to the family of Vorono{\"\i} transforms in the theory of
lattice points.
\medskip
\centerline{\bf 13. Sum formula. II}
\smallskip
\noindent
We are now ready to invert the sum formula $(10.1)$. The result is one
of the main assertions of the present article, and is embodied in
\smallskip
\noindent
{\bf Theorem 13.1. (Kloosterman--Spectral sum formula)}\quad {\it Let
$f$ be an even function on $\B{C}^\times$. Let us suppose that there
exist constants $\rho$ and $\sigma$ such that
$0<\rho<{1\over2}<\sigma$, and
\item{1.} $f(u)=O(|u|^{2\sigma})$ as $|u|\downarrow0$,
\item{2.} $f$ is six times continuously differentiable, and for
$a+b\le6$
$$
\int_{\B{C}^\times}|(u\partial_u)^af(u)|^2
|u|^{b-2\rho}d^\times\!u<\infty.
$$
\par
\noindent
Then we have, for any non-zero $\omega_1,\omega_2\in\B{Z}[i]$,
$$
\leqalignno{ \sum_{c\ne0}&{1\over|c|^2}S_\r{F}(\omega_1,\omega_2;c)
f\left({2\pi\over c}\sqrt{\omega_1\omega_2}\right)&(13.1)
\cr
&= 2\pi\sum_{V}|c_V(1)|^2t_V(\omega_1)t_V(\omega_2)\r{K}
f(\nu_V,p_V)\cr
&-i\sum_{p\in2\B{Z}}\left({\omega_1\omega_2
\over|\omega_1\omega_2|}\right)^p
\int_{(0)}{\sigma_\nu(\omega_1,-p/2)\sigma_\nu(\omega_2,-p/2)
\over|\omega_1\omega_2|^\nu|\zeta_\r{F}(1+\nu,p/2)|^2}
\r{K}f(\nu,p)\,d\nu, }
$$
where the transformation $\r{K}$ is defined in $(11.1)$, and $V$ runs
over all Hecke invariant right-irreducible cuspidal subspaces of
$L^2(\varGamma\backslash G)$ together with the specifications in
Section 8. The contour is the imaginary axis, and the convergence is
absolute throughout.}
\smallskip
\noindent
Proof. We denote the left and the right sides of $(13.1)$ by
$L^*_{\omega_1,\omega_2}(f)$ and $R^*_{\omega_1,\omega_2}(f)$,
respectively. Let $X>0$ be large, and let $\chi(r)$ be a smooth
function which is equal to $1$ for $1/X\le r\le X$, and $0$ for
$r\le1/(2X)$ and $r\ge2X$, and monotonic otherwise. Also let $\phi$
be a smooth function on $\B{C}$ such that
$$
\int_{\B{C}^\times}\phi\, d^\times\!u=1,\quad \hbox{$\phi(u)\ge0$, and
$\phi(u)=0\,$ if $\,|u-1|\ge X^{-1}$}.\leqno(13.2)
$$
With these we put
$$
f_X(u)=\int_{\B{C}^\times}\phi(u/v) \chi(|v|)f(v)\,
d^\times\!v.\leqno(13.3)
$$
This function is smooth and compactly supported on $\B{C}^\times$, and
converges pointwise to $f$ as $X\uparrow\infty$. By virtue of
Theorems 10.1 and 11.1 coupled with Lemma 11.1, we have
$$
L^*_{\omega_1,\omega_2}(f_X)
=R^*_{\omega_1,\omega_2}(f_X).\leqno(13.4)
$$
By the condition {\it 1} and by the bound $(8.14)$ we have readily
$$
\lim_{X\to\infty} L^*_{\omega_1,\omega_2}(f_X)=
L^*_{\omega_1,\omega_2}(f).\leqno(13.5)
$$
To deal with the right side of $(13.4)$, we observe that
$$
u\partial_uf_X(u)
=\int_{\B{C}^\times}\phi(u/v)
v\partial_v[\chi(|v|)f(v)]\,d^\times\!v,\leqno(13.6)
$$
and that
$$
\r{K}[\r{b}_uf_X](\nu,p)= (\nu-p)^2\r{K}f_X(\nu,p)\leqno(13.7)
$$
with $\r{b}_u=(u\partial_u)^2+u^2$.
The latter is due to the fact that $\e{J}_{\nu,p}$ and thus
$\e{K}_{\nu,p}$ are eigenfunctions of Bessel's differential operator
$\r{b}_u$ with eigenvalue $(\nu-\nobreak p)^2$. Invoking Lemma 12.1,
we have
$$
\r{K}f_X(\nu,p) \ll \min\left(\Vert f_X\Vert_\rho,\, \Vert\r{b}_u^3f_X
\Vert_\rho |\nu-p|^{-6}\right) \leqno(13.8)
$$
for any $\nu\in i\B{R}$, $p\in\B{Z}$, where $\Vert\cdot\Vert_\rho$ is
the norm of the Hilbert space
$L^2(\B{C}^\times,|u|^{-2\rho}d^\times\!u)$. By the definition we
have $\Vert f_X\Vert_\rho\ll\Vert f \Vert_\rho$. A multiple use of
$(13.6)$ gives that
$$
\Vert u^{a}(u\partial_u)^bf_X\Vert_\rho\ll
\Vert u^{a}(u\partial_u)^b\chi f\Vert_\rho\ll \sum_{j=0}^b\Vert
u^{a}(u\partial_u)^jf\Vert_\rho,
\leqno(13.9)
$$
since $(u\partial_u)^j\chi$ is bounded. This and the condition {\it 2}
imply that $\Vert\r{b}_u^3f_X \Vert_\rho$ is uniformly bounded. That
is, we have, from $(13.8)$,
$$
\r{K}f_X(\nu,p) \ll(1+|\nu|+|p|)^{-6} \leqno(13.10)
$$
uniformly in all involved parameters. Following the argument leading
to $(10.5)$, we find that
$$
\lim_{X\to\infty} R^*_{\omega_1,\omega_2}(f_X)
=R^*_{\omega_1,\omega_2}(f).\leqno(13.11)
$$
This ends the proof of the theorem.
\smallskip
\noindent
{\csc Remark.} Note the similarity between Theorem 13.1 and Theorem
2.3 of [30]. Comparing $(11.1)$ with $(2.4.8)$ of [30], the
impression will be enhanced. One may improve Theorem 13.1 by relaxing
the second condition. It appears, however, that our assertion is
sufficiently precise for practical purposes. For a possible
alternative approach to Theorem 13.1, see the final section.
\par
The identity $(13.4)$ does not contain the delta-term corresponding to
that in $(10.1)$. This is of course due to Lemma 11.1. Then, it might
be worth remarking that the assertion $(11.32)$ is a consequence of
the spectral sum formula as well. Namely, one may prove it
alternatively in the following way: We write $(13.4)$ with the $f$ in
Lemma 11.1, in place of $f_X$, but without using the lemma, so that
the delta-term remains in the identity. Specialize it by setting
$\omega_1=1$ and $\omega_2=n$, multiply both sides by
$\zeta_\r{F}(1-\beta)\sigma_{-\alpha}(n)|n|^{\alpha+\beta-1}$ with
$|\Re\alpha|+\Re\beta<-2$, and sum over all non-zero $n\in\B{Z}[i]$.
Applying some rearrangement partly depending on the arguments in
Sections 2 and 14, one is lead to the conclusion that
$\zeta_\r{F}(1-\beta)$ times the left side of $(11.32)$ is regular at
$\beta=0$, which is naturally equivalent to $(11.32)$. This is
definitely far more complicated than the above proof, but seems to
have certain interest of its own.
\par
It seems possible to extend Theorem 13.1 to other discrete subgroups
of $\r{PSL}_2(\B{C})$, provided we have a non-trivial bound for the
sums corresponding to our $S_\r{F}$. Otherwise, the test functions in
Theorem 10.1 have to live on a wider strip $|\Re\nu|\leq
1+\varepsilon$ with an $\varepsilon>0$, and have prescribed zeros at
$\nu=\pm1$. Then the problem is that the transforms $\r{K}f$, for $f$
even, smooth, compactly supported, need not posses those zeros.
\medskip
\centerline{\bf 14. An explicit formula}
\smallskip
\noindent
In this final section we shall apply Theorem 13.1 to the sum $(2.31)$
and establish a spectral decomposition of $\e{Z}_2(g,\r{F})$ defined
in $(1.1)$. The underlying principle is the same as in the rational
case but the procedure is naturally more involved (cf.\ Sections
4.4--4.7 of [30]).
\smallskip
We have first to examine if the two conditions in Theorem 13.1 are
satisfied by the function $[g](u;\alpha,\beta,\gamma,\delta)$ while
$(2.27)$ is assumed. The first condition is easy to check. As to the
second we observe that $[g]$ is smooth and, together with its
derivatives, of rapid decay as $|u|\uparrow\infty$; indeed, by virtue
of Lemma 2.1, it is enough to move the contour in $(2.32)$ far to the
left. Thus we may restrict ourselves to the vicinity of the point
$u=0$. Then we note that $u\partial_u=
{1\over2}(r\partial_r-i\partial_\theta)$ with $u=re^{i\theta}$. This
implies that when applied to $[g]$ the operator $u\partial_u$
does not change essentially the asymptotic behaviour of the function
as $|u|\downarrow0$. Note that again Lemma 2.1 plays a r\^ole. Hence
we have to check only the case $a=b=0$ in the condition {\it 2\/} of
Theorem 13.1. The confirmation is then immediate.
\par
Now, let us put
$$
\Phi_p(\nu)=\Phi_p(\nu;\alpha,\beta,\gamma,\delta;g)
=\int_{\B{C}^\times}\e{K}_{\nu,p}
(u)\,[g](u;\alpha,\beta,\gamma,\delta)d^\times\!u.\leqno(14.1)
$$
Theorem 13.1 gives, on $(2.27)$,
$$
S_{m,n}(\alpha,\beta,\gamma,\delta;g)= \{S^{(1)}_{m,n}+S^{(2)}_{m,n}\}
(\alpha,\beta,\gamma,\delta;g),\leqno(14.2)
$$
where
$$
\leqalignno{\qquad S^{(1)}_{m,n}(\alpha,\beta,\gamma,\delta;g)&=
2\pi\sum_{V}|c_V(1)|^2t_V(m)t_V(n) \Phi_{p_V}(\nu_V),&(14.3)\cr
S^{(2)}_{m,n}(\alpha,\beta,\gamma,\delta;g)
&=-i\sum_{p\in2\B{Z}}\left({mn \over|mn|}\right)^p
\int_{(0)}{\sigma_\nu(m,-p/2)\sigma_\nu(n,-p/2)
\over|mn|^\nu|\zeta_\r{F}(1+\nu,p/2)|^2} \Phi_p(\nu)\,d\nu.
&(14.4) }
$$
To consider the function $\Phi_p$, let $K_{\nu,p}(r,q)$ and $G(r,q)$
be the $2q$-th Fourier coefficients, in $\arg u$, of the functions
$\e{K}_{\nu,p}(u)$ and $[g](u;\alpha,\beta,\gamma,\delta)$,
respectively. We have, from $(2.32)$,
$$
\leqalignno{ G(r,q)&=(r/2)^{-2(1+\alpha+\beta)}&(14.5)\cr
&\times\int_{(\eta)}{\Gamma(1-s+{1\over2}|q|)
\Gamma(1+\alpha-s+{1\over2}|q|)\over \Gamma(s+{1\over2}|q|)
\Gamma(s-\alpha+{1\over2}|q|)}
\tilde{g}_q(s;\gamma,\delta)(r/2)^{4s}ds\cr
}
$$
with $\eta<1+\min(0,\Re\alpha)$. Note that on $(12.17)$ the assertion
$(12.18)$ induces
$$
\sum_{q\in\B{Z}}\int_0^\infty|K_{\nu,p}(r,q)|^2
r^{2\rho-1}dr<\infty,\leqno(14.6)
$$
and that by the above discussion we have, given $(2.27)$ and
$0<\rho<{1\over2}$,
$$
\sum_{q\in\B{Z}}\int_0^\infty|G(r,q)|^2
r^{-2\rho-1}dr<\infty.\leqno(14.7)
$$
These imply that
$$
\Phi_p(\nu)=2\pi\sum_{q\in\B{Z}} \int_0^\infty
K_{\nu,p}(r,q)G(r,q){dr\over r}.\leqno(14.8)
$$
Further, let $\tilde{K}_{\nu,p}(s,q)$ and $\tilde{G}(s,q)$ be the
Mellin transforms of $K_{\nu,p}(r,q)$ and $G(r,q)$, respectively, as
functions of $r$. Then the last expression is transformed into
$$
\Phi_p(\nu)=-i\sum_{q\in\B{Z}}\int_{(\rho)}
\tilde{K}_{\nu,p}(s,q)\tilde{G}(-s,q)ds.\leqno(14.9)
$$
This is the result of an application of the Parseval formula for
Mellin transforms to each integral in $(14.8)$ (see Theorem 72 of
[37]). The formula $(14.5)$ is obviously equivalent to
$$
\leqalignno{ \tilde{G}(s,q)&=\pi i
\tilde{g}_q(\txt{1\over2}(1+\alpha+\beta)-\txt{1\over4}s
;\gamma,\delta)&(14.10)\cr
&\times2^{s-1}{\Gamma({1\over4}s+{1\over2}(1-\alpha-\beta+|q|))
\Gamma({1\over4}s+{1\over2}(1+\alpha-\beta+|q|))\over
\Gamma({1\over2}(1+\alpha+\beta+|q|)-{1\over4}s)
\Gamma({1\over2}(1-\alpha+\beta+|q|)-{1\over4}s)},
\cr
}
$$
provided $(2.27)$ and $\Re s>2(\Re\beta+|\Re\alpha|-1)$. To find
$K_{\nu,p}(r,q)$ we combine Theorem 12.1 with Graf's addition theorem
(formula $(1)$ on p.\ 359 of [42]), getting
$$
K_{\nu,p}(r,q)=(-1)^{\max(|p|,|q|)}{2\over\pi}\int_0^\infty
y^{2\nu-1}J_{|p+q|}(ry) J_{|p-q|}(r/y)dy\leqno(14.11)
$$
for $|\Re\nu|<{1\over4}$ (cf.\ $(12.15)$). Via the formula
$$
\int_0^\infty J_{\xi}(y)y^{s-1}dy=2^{s-1}{\Gamma({1\over2}(\xi+s))
\over\Gamma({1\over2}(\xi-s)+1)}\leqno(14.12)
$$
with $-\Re\xi<\Re s<{1\over2}$, and the Parseval formula for Mellin
transforms, one may express $K_{\nu,p}(r,q)$ by an inverse Mellin
transform. Then we get
$$
\leqalignno{
&\tilde{K}_{\nu,p}(s,q)
=(-1)^{\max(|p|,|q|)}{2^{s-2}\over\pi}&(14.13)\cr
&\times {\Gamma({1\over4}s+{1\over2}(\nu+|p+q|))
\Gamma({1\over4}s-{1\over2}(\nu-|p-q|))
\over\Gamma(1-{1\over4}s-{1\over2}(\nu-|p+q|))
\Gamma(1-{1\over4}s+{1\over2}(\nu+|p-q|))}, }
$$
provided $2|\Re \nu|<\Re s<1-2|\Re\nu|$.
\par
The combination of $(14.9)$, $(14.10)$, and $(14.13)$ yields
\smallskip
\noindent
{\bf Lemma 14.1.}\quad{\it The function $\Phi_p$ continues
meromorphically to $\B{C}^5$. We have the representation
$$
\leqalignno{ \Phi_p(\nu;\alpha,\beta&,\gamma,\delta;g)
={1\over2}\sum_{q\in\B{Z}}(-1)^{\max(|p|,|q|)}&(14.14)\cr
&\times\int_{-i\infty}^{i\infty}
\tilde{g}_q(s;\gamma,\delta)\Gamma_q(s;\alpha)
\Gamma_{p,q}(s,\nu;\alpha,\beta) ds,\cr
}
$$
where
$$
\leqalignno{ \Gamma_q(s;\alpha)&={\Gamma(1-s+{1\over2}|q|)
\Gamma(1-s+\alpha+{1\over2}|q|)\over \Gamma(s+{1\over2}|q|)
\Gamma(s-\alpha+{1\over2}|q|)},&(14.15)\cr
\Gamma_{p,q}(s,\nu;\alpha,\beta)&=
{\Gamma(s-{1\over2}(\alpha+\beta-\nu+1)
+{1\over2}|p+q|)\over\Gamma(1-s+{1\over2} (\alpha+\beta-\nu+1)
+{1\over2}|p+q|)}&(14.16)\cr
&\times{\Gamma(s-{1\over2}(\alpha+\beta+\nu+1) +{1\over2}|p-q|)\over
\Gamma(1-s+{1\over2}(\alpha+\beta+\nu+1) +{1\over2}|p-q|)}.\cr
}
$$
In $(14.14)$ it is supposed that the poles of
$\tilde{g}_q(s;\gamma,\delta)\Gamma_q(s;\alpha)$ and those of
$\Gamma_{p,q}(s,\nu;\alpha,\beta)$ are separated by the contour to
the right and the left, respectively; and the parameters are such
that the contour can be drawn. It follows in particular that if
$\Re\nu$, $\alpha$, $\beta$, $\gamma$, and $\delta$ are bounded, then
we have, for any fixed $A>0$,
$$
\Phi_p(\nu;\alpha,\beta,\gamma,\delta;g)\ll (1+|\nu|+|p|)^{-A}
\leqno(14.17)
$$
as $|\nu|+|p|$ tends to infinity.}
\smallskip
\noindent
Proof. We assume first $(2.27)$ and $(12.17)$. Then we may insert
$(14.10)$ and $(14.13)$ into $(14.9)$. After the change of variable
$s\mapsto 4s-2(1+\alpha+\beta)$, we get the expression $(14.14)$ with
the contour $(0)$. The expression for the general situation follows
by analytic continuation. The meromorphy of $\Phi_p$ is an immediate
consequence of $(14.14)$. As to the bound $(14.17)$, we need only to
shift the contour in $(14.14)$ far to the left. The resulting
integral and the residues are estimated by Stirling's formula and
$(2.12)$. This ends the proof.
\smallskip
We assume $(2.27)$ again, and collect $(2.30)$--$(2.31)$ and
$(14.2)$--$(14.4)$. We obtain
$$
B_m^{(1)}(\alpha,\beta;g^*(\cdot\,;\gamma,\delta))=
\left\{B_m^{(1,1)}+B_m^{(1,2)}\right\}
(\alpha,\beta;g^*(\cdot\,;\gamma,\delta)).\leqno(14.18)
$$
Here we have
$$
\leqalignno{ B_m^{(1,1)}&(\alpha,\beta;g^*(\cdot\,;\gamma,\delta))
=-2i\pi^{2\beta}|m|^{\alpha+\beta+1} \sum_V|c_V(1)|^2t_V(m)&(14.19)\cr
&\qquad\times H_V(\txt{1\over2}(1-\alpha-\beta))
H_V(\txt{1\over2}(1+\alpha-\beta))\Phi_{p_V}(\nu_V),\cr
}
$$
and
$$
\leqalignno{ B_m^{(1,2)}(\alpha,&\beta;g^*(\cdot\,;\gamma,\delta))
=-\pi^{2\beta-1}|m|^{\alpha+\beta+1}&(14.20)\cr
&\times\sum_{p\in\B{Z}} \left({m\over|m|}\right)^{4p}
\int_{(0)}{\sigma_\nu(m,-2p)
\over|m|^\nu|\zeta_\r{F}(1+\nu,2p)|^2}\cr
&\times\zeta_\r{F}(\txt{1\over2}(1-\alpha-\beta+\nu),p)
\zeta_\r{F}(\txt{1\over2}(1-\alpha-\beta-\nu),-p)\cr
&\times\zeta_\r{F}(\txt{1\over2}(1+\alpha-\beta+\nu),p)
\zeta_\r{F}(\txt{1\over2}(1+\alpha-\beta-\nu),-p)
\Phi_{4p}(\nu)d\nu.\cr
}
$$
The expression $(14.19)$ depends on $(8.15)$, $(8.21)$ and $(10.13)$.
On the other hand, $(14.20)$ depends on the following formula: For
any $a,b,c\in\B{Z}$, $\tau,\xi\in\B{C}$ we have, in the region of
absolute convergence,
$$
\leqalignno{
&{1\over4}\sum_{n\ne0}(n/|n|)^{4a}\sigma_\tau(n,b)\sigma_\xi(n,c)
|n|^{-2s}&(14.21)\cr
&={\zeta_\r{F}(s,a)\zeta_\r{F}(s-\tau,a+b) \zeta_\r{F}(s-\xi,a+c)
\zeta_\r{F}(s-\tau-\xi,a+b+c) \over\zeta_\r{F}(2s-\tau-\xi,2a+b+c)}.
\cr}
$$
This is an extension of Ramanujan's well-known identity for the
product of four values of the Riemann zeta-function, and the proof is
similar (see $(1.3.3)$ of [36]).
\smallskip
Returning to $(2.6)$, we specialize $(2.29)$ and $(14.18)$ with
$\alpha=z_1-z_2$, $\beta=z_3-z_4$, $\gamma=z_1$, and $\delta=z_3$. We
shall assume temporarily that
$$
1<\Re z_1<\Re z_2<\Re z_1+1,\quad 1<\Re z_3<\Re z_4-3.\leqno(14.22)
$$
The condition $(2.27)$ is satisfied. Hence, via $(2.28)$--$(2.29)$ and
$(14.18)$--$(14.20)$, the formula $(2.6)$ is transformed into
$$
\e{I}(z_1,z_2,z_3,z_4;g)=\left\{\e{I}^{(r)}+\e{I}^{(c)}+\e{I}^{(e)}
\right\}(z_1,z_2,z_3,z_4;g)\leqno(14.23)
$$
in the domain $(14.22)$. Here we have
$$
\leqalignno{\hskip 1cm \e{I}^{(r)}(z_1&,z_2,z_3,z_4;g)
={{\zeta_\r{F}(z_1+z_3)\zeta_\r{F}(z_1+z_4)\zeta_\r{F}(z_2+z_3)
\zeta_\r{F}(z_2+z_4)}\over4\zeta_\r{F}(z_1+z_2+z_3+z_4)}
\hat{g}(0)&(14.24)\cr
&+\pi{{\zeta_\r{F}(z_1+z_3-1)\zeta_\r{F}(z_2+z_4)
\zeta_\r{F}(1+z_2-z_1)\zeta_\r{F}(1+z_4-z_3)}
\over2\zeta_\r{F}(2+z_2+z_4-z_1-z_3)} \tilde{g}_0(1;z_1,z_3)\cr
&+\pi{{\zeta_\r{F}(z_2+z_3-1)\zeta_\r{F}(z_1+z_4)
\zeta_\r{F}(1+z_1-z_2)\zeta_\r{F}(1+z_4-z_3)}
\over2\zeta_\r{F}(2+z_1+z_4-z_2-z_3)}\cr
&\hskip 1cm\times\tilde{g}_0(1+z_1-z_2;z_1,z_3), }
$$
$$
\leqalignno{\qquad
\e{I}^{(c)}(z_1,&z_2,z_3,z_4;g)={\pi^{2(z_3-z_4)}\over2i}\sum_V
|c_V(1)|^2H_V(\txt{1\over2}(z_1+z_2+z_3+z_4-1))&(14.25)\cr
&\times H_V(\txt{1\over2}(z_2+z_4-z_1-z_3+1))
H_V(\txt{1\over2}(z_1+z_4-z_2-z_3+1))\cr
&\hskip 1cm\times\Phi_{p_V}(\nu_V;z_1-z_2,z_3-z_4,z_1,z_3;g), }
$$
and
$$
\e{I}^{(e)}(z_1,z_2,z_3,z_4;g)=\sum_{p\in\B{Z}}\e{I}^{(e)}_p
(z_1,z_2,z_3,z_4;g),\leqno(14.26)
$$
where
$$
\leqalignno{\qquad
\e{I}^{(e)}_p&(z_1,z_2,z_3,z_4;g)=-{\pi^{2(z_3-z_4)-1}\over4}
\int_{(0)}
\zeta_\r{F}(\txt{1\over2}(z_1+z_2+z_3+z_4-1+\nu),p)&(14.27)\cr
&\times\zeta_\r{F}(\txt{1\over2}(z_1+z_2+z_3+z_4-1-\nu),-p)
\zeta_\r{F}(\txt{1\over2}(z_2+z_4-z_1-z_3+1+\nu),p)\cr
&\times\zeta_\r{F}(\txt{1\over2}(z_2+z_4-z_1-z_3+1-\nu),-p)
\zeta_\r{F}(\txt{1\over2}(z_1+z_4-z_2-z_3+1+\nu),p)\cr
&\times\zeta_\r{F}(\txt{1\over2}(z_1+z_4-z_2-z_3+1-\nu),-p)
{\Phi_{4p}(\nu;z_1-z_2,z_3-z_4,z_1,z_3;g)
\over\zeta_\r{F}(1+\nu,2p)\zeta_\r{F}(1-\nu,-2p)}d\nu. }
$$
We then observe that by Lemma 8.1 the function $H_V(s)$ is entire, and
of polynomial order in $s$, $p_V$ and $\nu_V$ if $\Re s$ is bounded.
Thus $\e{I}^{(c)}$ is meromorphic over $\B{C}^4$ by virtue of
Corollary 10.1 and Lemma 14.1. To see the situation at the point
$\r{p}_{1\over2}$, we note that if $(z_1,z_2,z_3,z_4)$ is
sufficiently close to $\r{p}_{1\over2}$, and $\nu\in i\B{R}$ (see
$(8.4)$), then we may take $({3\over4})$ as a contour in $(14.14)$.
This implies readily that $\e{I}^{(c)}$ is regular at
$\r{p}_{1\over2}$, and we have
$$
\e{I}^{(c)}(\r{p}_{1\over2};g)={1\over2i}
\sum_V|c_V(1)|^2H_V(\txt{1\over2})^3
\Phi_{p_V}(\nu_V;0,0,\txt{1\over2},\txt{1\over2};g).\leqno(14.28)
$$
As to the sum of $\e{I}^{(e)}_p$ over $p\ne0$, it is analogous to
$\e{I}^{(c)}$, and we have
$$
\sum_{\scr{p\in\B{Z}}\atop\scr{p\ne0}}
\e{I}^{(e)}_p(\r{p}_{1\over2};g)=-{1\over4\pi}
\sum_{\scr{p\in\B{Z}}\atop\scr{p\ne0}}\int_{(0)}
{|\zeta_\r{F}({1\over2}(1+\nu),p)|^6\over|\zeta_\r{F}(1+\nu,2p)|^2}
\Phi_{4p}(\nu;0,0,\txt{1\over2},\txt{1\over2};g)d\nu.\leqno(14.29)
$$
It remains to consider the function $\e{I}_0^{(e)}$. We note first
that it is meromorphic over $\B{C}^4$. This can be proved either
shifting the contour appropriately or simply observing that all
terms, except for $\e{I}_0^{(e)}$, in $(14.23)$ are already known to
be meromorphic over $\B{C}^4$. Thus, to see the nature of
$\e{I}_0^{(e)}$ near the point $\r{p}_{1\over2}$, we may let
$(z_1,z_2,z_3,z_4)$ approach to $\r{p}_{1\over2}$ in a specific way,
as we shall do shortly.
\par
We start from the domain defined by $(14.22)$, where we have the
representation $(14.27)$ with $p=0$. We shall move the contour,
closely following the discussion on the corresponding part of the
rational case (Section 4.7 of [30]). In the process we shall
 encounter singularities of the integrand, and the difficulty lies in
that they depend on $z_1,z_2,z_3,z_4$. To facilitate the discussion
we put
$$
\eqalign{
\nu_1=z_1+z_2+&z_3+z_4-3,\,\nu_2=z_2+z_4-z_1-z_3-1,\cr
&\nu_3=z_1+z_4-z_2-z_3-1. }\leqno(14.30)
$$
The zeta-part of the integrand in $\e{I}^{(e)}_0$ has singularities
only at the six points $\pm\nu_1$, $\pm\nu_2$, $\pm\nu_3$ and at the
zeros of $\zeta_\r{F}(1+\nu)\zeta_\r{F}(1-\nu)$. Then we make an
observation:
\smallskip
\noindent
{\bf Lemma 14.2.}\quad{\it The singularities of
$\Phi_0(\nu;z_1-z_2,z_3-z_4,z_1,z_3;g)$ as a function of $\nu$ is
contained in the set}
$$
\Big\{\pm(\nu_1+2a),\pm(\nu_2+2b),\pm(\nu_3+2c): \B{Z}\ni
a,b,c\ge1\Big\}.\leqno(14.31)
$$
Proof. The singularities can occur only when we are unable to draw the
contour in $(14.14)$, that is,
$$
\leqalignno{
&\Big\{\txt{1\over2}(\alpha+\beta\pm\nu+1) -l_\pm:\,\B{Z}\ni
l_\pm\ge0\Big\}&(14.32)\cr
&\quad\cap \Big\{1+l_1,\,1+\alpha+l_2,\,\gamma+\delta+l_3:\,\B{Z}\ni
l_1,l_2,l_3\ge0 \Big\}\not=\emptyset. }
$$
Such situations are covered by $(14.31)$ under the current
specialization. This ends the proof.
\smallskip
We now set
$$
z_1={1\over2}+t,\,z_2={1\over2}+2t,\,
z_3={1\over2}+t,\,z_4={1\over2}+6t;\quad |\Im
t|<\varepsilon_0\leqno(14.33)
$$
with a sufficiently small $\varepsilon_0$. If ${2\over3}<\Re t<1$,
then $(14.22)$ is satisfied, and moreover the points $\nu_1=10t-1$,
$\nu_2=6t-1$, $\nu_3=4t-1$ are not in the set $(14.31)$. Thus, the
last lemma implies that we can move the contour in $\e{I}_0^{(e)}$ so
that the points $\nu_1$, $\nu_2$, $\nu_3$ are on the left of the new
contour, but none of the points in $(14.31)$ and zeros
$\zeta_\r{F}(1-\nu)\zeta_\r{F}(1+\nu)$ are encountered in the
process. Leaving the residues at $\nu_1$, $\nu_2$ and $\nu_3$ for a
later discussion, we consider the resulting integral as a function of
$t$ as $t\to0$, while keeping $\Re t>0$ and moving the contour
stepwise. We observe, via the last lemma, that, except for the cases
$t={2\over3},\,{1\over2},\,{1\over3}$, we can draw the contour. These
exceptional points obviously make no trouble; for instance we may
assume $\Im t>0$. Thus the integral continues analytically to a small
right semicircle centered at the origin. Then, having $t$ in this
domain, we shift the contour back to the original, i.e., the
imaginary axis. This time we encounter singularities at $-\nu_1$,
$-\nu_2$, and $-\nu_3$ but none else. Note that at this stage we may
leave the specialization $(14.33)$, and suppose, instead, that
$(z_1,z_2,z_3,z_4)$ are in a small neighbourhood of
$\r{p}_{1\over2}$. The integral $\e{I}^{(e)}_{0,*}$ thus obtained is
regular at $\r{p}_{1\over2}$, and we see readily that
$$
\e{I}^{(e)}_{0,*}(\r{p}_{1\over2};g)=-{1\over4\pi}\int_{(0)}
{|\zeta_\r{F}({1\over2}(1+\nu))|^6\over
|\zeta_\r{F}(1+\nu)|^2}\Phi_0(\nu;0,0,\txt{1\over2},
\txt{1\over2};g)d\nu.\leqno(14.34)
$$
\par
Gathering these, we obtain the assertion
$$
\e{I}(\r{p}_{1\over2};g)=\Big\{M+\e{I}^{(c)}
+\e{I}^{(e)}_*\Big\}(\r{p}_{1\over2};g).\leqno(14.35)
$$
Here $\e{I}^{(e)}_*$ is the sum of the right side of $(14.27)$ over
all $p\in\B{Z}$, but with a different $(z_1,z_2,z_3,z_4)$. On the
other hand $M$ is the sum of $\e{I}^{(r)}$ and the contribution of
the poles at $\nu=\pm\nu_1$, $\pm\nu_2$, and $\pm\nu_3$ that we
encountered in the above procedure. We stress that $M$ is regular at
$\r{p}_{1\over2}$. This is because all other functions involved in
$(14.35)$ are regular at $\r{p}_{1\over2}$.
\par
Hence it remains for us to compute $M(z_1,z_2,z_3,z_4;g)$. We have, in
a small neighbourhood of $\r{p}_{1\over2}$,
$$
\leqalignno{\hskip 1cm
&\qquad M(z_1,z_2,z_3,z_4;g)=
{{\zeta_\r{F}(z_1+z_3)\zeta_\r{F}(z_1+z_4)\zeta_\r{F}(z_2+z_3)
\zeta_\r{F}(z_2+z_4)}\over4\zeta_\r{F}(z_1+z_2+z_3+z_4)}
\hat{g}(0)&(14.36)\cr
&+\pi{{\zeta_\r{F}(z_1+z_3-1)\zeta_\r{F}(z_2+z_4)
\zeta_\r{F}(1+z_2-z_1)\zeta_\r{F}(1+z_4-z_3)}
\over2\zeta_\r{F}(2+z_2+z_4-z_1-z_3)} \tilde{g}_0(1;z_1,z_3)\cr
&+\pi{{\zeta_\r{F}(z_2+z_3-1)\zeta_\r{F}(z_1+z_4)
\zeta_\r{F}(1+z_1-z_2)\zeta_\r{F}(1+z_4-z_3)}
\over2\zeta_\r{F}(2+z_1+z_4-z_2-z_3)}
\tilde{g}_0(1+z_1-z_2;z_1,z_3)\cr
&+i\pi^{2(z_3-z_4)+1} {\zeta_\r{F}(z_2+z_4-1)\zeta_\r{F}(2-z_1-z_3)
\zeta_\r{F}(z_1+z_4-1)\zeta_\r{F}(2-z_2-z_3)\over
4\zeta_\r{F}(4-z_1-z_2-z_3-z_4)}\cr
&\hskip1cm\times\Phi_0(z_1+z_2+z_3+z_4-3;
z_1-z_2,z_3-z_4,z_1,z_3;g)\cr
&+i\pi^{2(z_3-z_4)+1} {\zeta_\r{F}(z_2+z_4-1)\zeta_\r{F}(z_1+z_3)
\zeta_\r{F}(z_4-z_3)\zeta_\r{F}(z_1-z_2+1)\over
4\zeta_\r{F}(2-z_2-z_4+z_1+z_3)}\cr
&\hskip1cm\times\Phi_0(z_2+z_4-z_1-z_3-1;
z_1-z_2,z_3-z_4,z_1,z_3;g)\cr
&+i\pi^{2(z_3-z_4)+1} {\zeta_\r{F}(z_1+z_4-1)\zeta_\r{F}(z_2+z_3)
\zeta_\r{F}(z_4-z_3)\zeta_\r{F}(z_2-z_1+1)\over
4\zeta_\r{F}(2-z_1-z_4+z_2+z_3)}\cr
&\hskip1cm\times\Phi_0(z_1+z_4-z_2-z_3-1;
z_1-z_2,z_3-z_4,z_1,z_3;g).\cr
}
$$
This is obviously a linear integral transform of the weight function
$g$. The six members on the right side have singularities at
$\r{p}_{1\over2}$, but these have to cancel out each other as $M$ is
regular at the point. Thus, what matters actually are the constant
terms $m_j$ $(1\le j\le6)$, respectively, in the Laurent series
expansions of these members at $\r{p}_{1\over2}$. That is, we have
$$
M(\r{p}_{1\over2};g)=\sum_{j=1}^6m_j.\leqno(14.37)
$$
The computation of $m_j$ can be carried out in much the same way as in
the rational case (see pp.\ 176--178 of [30]). It is possible to
write $M(\r{p}_{1\over2};g)$ down explicitly in terms of $g$ and
derivatives of the $\Gamma$-function, but we stop here to restrict
ourselves to the description of the overall structure of our subject.
\smallskip
To state our final result we put
$\Lambda_{\nu,p}(g)=(2i)^{-1}\Phi_p(\nu;0,0,{1\over2}, {1\over2};g)$,
and $M_\r{F}(g)=M(\r{p}_{1\over2},g)\allowbreak + a_0g({1\over2}i)
+b_0g(-{1\over2}i)+ a_1g'({1\over2}i)+b_1g'(-{1\over2}i)$ with $a_0$,
$a_1$, $b_0$, $b_1$ as in $(2.2)$. We thus have established
\smallskip
\noindent
{\bf Theorem 14.1.}\quad {\it Let $g$ be as in $(1.1)$. Then, with the
transformations $M_\r{F}$ and $\Lambda_{\nu,p}$ defined above, we
have the identity
$$
\leqalignno{ {\eusm
Z}_2(g,\r{F})=M_\r{F}(g)&+\sum_V|c_V(1)|^2H_V(\txt{1\over2})^3
\Lambda_{\nu_V,p_V}(g)&(14.38)\cr&+{1\over2\pi
i}\sum_{p\in\B{Z}}\int_{(0)} {|\zeta_\r{F}({1\over2}(1+\nu),p)|^6\over
|\zeta_\r{F}(1+\nu,2p)|^2}\Lambda_{\nu,4p}(g)d\nu\,,\cr }
$$
where $V$ runs over all Hecke invariant right-irreducible cuspidal
subspaces of $L^2(\varGamma\backslash G)$ together with the
specifications in Section 8. The contour is the imaginary axis, and
the convergence is absolute throughout.}
\medskip
\noindent
{\csc Remark.} Note that despite the special nature of our dissection
in $(2.6)$ the arithmetic ingredients in $(14.38)$, i.e., the
functions $H_V(s)$ and $\zeta_\r{F}(s,p)$ are in fact defined over
integral ideals of $\B{Q}(i)$.
\medskip
\centerline{\bf 15. Concluding remarks}
\smallskip
\noindent
In this final section we shall develop a further discussion of
elements involved in our main result $(14.38)$, in the light of
recent developments made to understand the explicit formula for
$\e{Z}_2(g,\B{Q})$, the fourth moment of the Riemann zeta-function
$\zeta=\zeta_\B{Q}$. We shall also ponder on an intriguing nature of
$(14.38)$ that is briefly remarked in the introduction.
\smallskip
Thus, Theorem 4.2 of [30] is now translated into
$$
\e{Z}_2(g,\B{Q}) =M_\B{Q}(g)+\sum_\r{V}|\r{c}_\r{V}(1)|^2
\r{H}_\r{V}(\txt{1\over2})^3 \Lambda_{\nu_\r{V}}(g)+{1\over2\pi
i}\int_{(0)} {|\zeta({1\over2}(1+\nu))|^6\over
|\zeta(1+\nu)|^2}\Lambda_{\nu}(g)d\nu. \leqno(15.1)
$$
Here $M_\B{Q}$ has a construction similar to $M_\r{F}$, $\r{V}$ runs
over all Hecke invariant right-irreducible cuspidal subspaces of
$L^2(\r{PSL}_2(\B{Z})\backslash\r{PSL}_2(\B{R}))$, and
$\r{c}_\r{V}(n)$ $(\B{Z}\ni n\ne0$) are the Fourier coefficients of
$\r{V}$, to which the Hecke series $\r{H}_\r{V}$ is associated. The
$\nu_\r{V}$ is the spectral parameter of $\r{V}$; that is, being
restricted to $\r{V}$, the Casimir operator over $\r{PSL}_2(\B{R})$
becomes the constant multiplication $(\nu_\r{V}^2-{1\over4})\cdot1$.
The functional $\Lambda_\nu$ is to be made precise shortly.
\par
The similarity between the formulas $(14.38)$ and $(15.1)$ appears to
the authors to suggest the existence of a geometric structure yet to
be discovered. In particular, these results are expected to extend to
a wide family of automorphic $L$-functions (cf.\ [19]). To enhance
this observation, we quote, from [34] with minor changes of notation,
the integral representation
$$
\Lambda_\nu(g)=\int_0^\infty{\hat{g}(\log(1+1/r))\over
(r(r+1))^{1\over2}}\Xi_\nu(r)dr,\leqno(15.2)
$$
where $\hat{g}$ is as above, and
$$
\Xi_\nu(r)={1\over2}\int_{\B{R}^\times}
 j_\nu(u/r)j_0(-u){d^\times\!u\over\sqrt{|u|}},\quad
\B{R}^\times=\B{R}\backslash\{0\},\quad
d^\times\!u=du/|u|.\leqno(15.3)
$$
Here
$$
j_\nu(u)=\pi{\sqrt{|u|}\over\sin\pi\nu}
\left\{J_{-\nu}^{\r{sgn}(u)}(4\pi\sqrt{|u|})-
J_\nu^{\r{sgn}(u)}(4\pi\sqrt{|u|})\right\},\leqno(15.4)
$$
with $J_\nu^+=J_\nu$, $J_\nu^-=I_\nu$. Correspondingly, we have, for
$\r{F}=\B{Q}(i)$,
$$
\Lambda_{\nu,p}(g)=\int_\B{C}{\hat{g}(2\log|1+1/u|)
\over|u(u+1)|}\Xi_{\nu,p}(u)d_+\!u,\leqno(15.5)
$$
where
$$
\Xi_{\nu,p}(u)={1\over16\pi}\int_{\B{C}^\times}
 j_{\nu,p}\left(\sqrt{v/u}\right)
j_{0,0}\left(\sqrt{-v}\right){d^\times\!v\over|v|}, \leqno(15.6)
$$
with
$$
\leqalignno{ j_{\nu,p}(u)&=2\pi^2|u|^2\e{K}_{\nu,p}(2\pi u) &(15.7)
\cr
&=2\pi^2{|u|^2\over\sin\pi\nu} \left(\e{J}_{-\nu,-p}(2\pi
u)-\e{J}_{\nu,p}(2\pi u)\right). }
$$
The last three formulas follow readily from $(12.26)$, with $p=0$, and
$(14.1)$. As a matter of fact, the Bessel kernel $j_\nu$ that
originates in Kuznetsov's works [21], [22] can be identified as the
Bessel function of irreducible representations of the group
$\r{PSL}_2(\B{R})$, and so is the $j_{\nu,p}$ in its relation with
$\r{PSL}_2(\B{C})$. For these facts see [8], whose announcement is in
[7], [34]. This interpretation of $j_\nu$ based on the harmonic
analysis over $\r{PSL}_2(\B{R})$ is originally due to Cogdell and
Piatetski-Shapiro [9]; and [8] contains an alternative and rigorous
approach to it via the concept of local functional equations of
Jacquet and Langlands [18], including its extension
$(15.5)$--$(15.7)$ as well. Hence, the resemblance between $(14.38)$
and $(15.1)$ in fact reaches deeper than the sheer outlook suggests.
At any events, the last expressions show how tightly the mean values
of zeta-functions are related to the structure of function spaces
over linear Lie groups.
\par
This is naturally the same with sum formulas of Kloosterman sums. The
work [9] in fact indicates a way to directly connect Kuznetsov's sum
formula with $j_\nu$, without the inversion procedure of the
Spectral-Kloosterman sum formula for $\r{PSL}_2(\B{R})$ as Kuznetsov
did. Since the principal means on which [9] is based have been
extended to the present situation, as remarked above, one may argue
that we could prove our Theorem 13.1 without first establishing
Theorem 10.1. That appears to be the case, but in the present work we
have chosen the way to extend the argument of [26] to include all
$K$-aspects. This is because the combination of the Jacquet and the
Goodman--Wallach operators provides us with a flexibility, perhaps
greater than the extension of [9] could. Moreover, the present
version of the sum formula for $\r{PSL}_2(\B{C})$ is more suitable
than existing ones for applications in the study of Kloosterman sums
and in the investigation of the distribution of automorphic spectral
data. It should, however, be remarked that Theorem 12.1, the above
proof of which depends on the Goodman--Wallach operator, could be
derived also from the interpretation of $j_{\nu,p}$ as Bessel
functions of representation of $\r{PSL}_2(\B{C})$. For this see [8],
[34].
\smallskip
We now turn to a comparison between $(14.38)$ and $(15.1)$, in their
asymptotic aspects, which one may assert is more important than to
discuss their structural similarities. Here exists a remarkable
difference between these formulas. That concerns the nature of the
term $M_\r{F}(g)$ in $(14.38)$. In $(15.1)$ the $M_\B{Q}(g)$ is
indeed the main term in the sense that with a specialization of $g$
it gives rise to the main term in the asymptotic formula
$$
\int_{-T}^T|\zeta(\txt{1\over2}+it)|^4dt=TP_4(\log
T)+O(T^{2\over3}\log^8T)\leqno(15.8)
$$
as $T\uparrow\infty$, where $P_4$ is a polynomial of order 4 (Theorem
5.2 of [30]). The function that is the counterpart of $M$ for $\zeta$
has an expression similar to $(14.36)$, and at the point
$\r{p}_{1\over2}$ the terms involved in it have singularities of
order 4 at most; see Section 4.7 of [30]. Analogously the functions
on the right of $(14.36)$ have singularities of order 4 at
$\r{p}_{1\over2}$, and the structure of $M_\r{F}(g)$ is similar to
that of $M_\B{Q}(g)$. We have, however, the lower bound
$$
\int_{-T}^T|\zeta_\r{F}(\txt{1\over2}+it)|^4dt\gg
T\log^8T,\leqno(15.9)
$$
which can be proved by the argument in Section 7.19 of [36]. Thus, it
is hard to regard $M_\r{F}(g)$ as the main term in $(14.38)$, in the
present context. This appears to raise a basic question about the
fourth moment of $\zeta_\r{F}$, and remotely the same about the
eighth moment of the Riemann zeta-function. Most probably it would be
expedient to study the Mellin transform
$$
Z_2(s,\r{F})=\int_1^\infty|\zeta_\r{F}(\txt{1\over2}+it)|^4t^{-s}dt
\leqno(15.10)
$$
as in the rational case, i.e., [29] and section 5.3 of [30], in order
to see finer analytic aspects of $(14.38)$. This and the related
issues such as the spectral mean values of $H_V({1\over2})$ are left
for future works. That is, the explicit formula $(14.38)$ is to be
regarded as just the beginning of a new story.
\par
The above argument is, however, not the unique way to investigate the
fourth moment of $\zeta_\r{F}$. This is already mentioned in [33],
and concerns a possible extension of Vinogradov--Takhtadjan--Jutila's
functional treatment of the fourth moment of the Riemann
zeta-function ([39], [19]). In the above we have obtained all
necessary means for this purpose: especially, the spectral sum
formula $(10.1)$, the Fourier expansion of Eisenstein series
$(5.32)$, and results on Bessel transforms.
\par
Finally, we remark that the mean square of the Dedekind zeta-function
of any imaginary or real quadratic number field can be treated with
the sum formula over the relevant congruence subgroup of the modular
group; and the result is analogous to that for $\e{Z}_2(g,\B{Q})$
(see [31]). Also, repeating a remark in the introduction, the present
authors have established, in [6], an explicit spectral decomposition
for the fourth power moment of the Dedekind zeta-function of any real
quadratic number field with class number one. The necessary spectral
theory is developed over the Hilbert modular group attached to the
respective field. \vskip 1.5cm \centerline{REFERENCES}
\smallskip
\item{[1]} R.W. Bruggeman. {\it Fourier Coefficients of Automorphic
Forms}. Lecture Notes in Math., {\bf865}, Springer-Verlag, Berlin,
1981.
\item{[2]} R.W. Bruggeman and R.J. Miatello. Estimates of Kloosterman
sums for groups of real rank one. {\it Duke Math.\ J}., {\bf 80}
(1995), 105--137.
\item{[3]} ---. Sum formula for $\r{SL}_2$ over a number field and
Selberg type estimate for exceptional eigenvalues. {\it Geom.\
Funct.\ Anal.}, {\bf8} (1998), 627--655.
\item{[4]} R.W. Bruggeman, R.J. Miatello and I. Pacharoni. Estimates
for Kloosterman sums for totally real number fields. {\it J. reine
angew.\ Math.}, {\bf 535} (2001), 103-164.
\item{[5]} R.W. Bruggeman and Y. Motohashi. A note on the mean value
of the zeta and $L$-functions.\ X. {\it Proc.\ Japan Acad.}, {\bf
77}(A) (2001), 111-114.
\item{[6]} ---. Fourth power moment of Dedekind zeta-functions of real
quadratic number fields with class number one. {\it Functiones et
Approximatio}, {\bf 29} (2001), 41--79.
\item{[7]} ---. A note on the mean value of the zeta and
$L$-functions.\ XIII. preprint.
\item{[8]} ---. Projections Poincar\'e series into irreducible
subspaces. preprint.
\item{[9]} J.W. Cogdell and I.I. Piatetski-Shapiro. {\it The
Arithmetic and Spectral Analysis of Poincar\'e Series}. Academic
Press, San Diego, 1990.
\item{[10]} J. Elstrodt, F. Grunewald and J. Mennicke. {\it Groups
Acting on Hyperbolic Space}. Springer -Verlag, Berlin, 1998.
\item{[11]} J.D. Fay. Fourier coefficients of the resolvent for a
Fuchsian group. {\it J. Reine Angew.\ Math.}, {\bf 293}/{\bf294}
(1977), 143--203.
\item{[12]} R. Goodman and N.R. Wallach. Whittaker vectors and conical
vectors. {\it J.\ Funct.\ Anal.}, {\bf39} (1980), 199--279.
\item{[13]} Harish Chandra. {\it Automorphic Forms on Semisimple Lie
Groups}. Lecture Notes in Math., {\bf 62}, Springer-Verlag, Berlin,
1968.
\item{[14]} E. Hecke. Eine neue Art von Zetafunktionen und ihre
 Beziehungen zur Verteilung der Primzahlen. {\it Math.\ Z.}, {\bf 6}
(1920), 11--51.
\item{[15]} G. Humbert. Sur la mesure des Classes d'Hermite de
discriminant donn\'e dans un corps quadratique imaginaire, et sur
certains volumes non euclidiens. {\it C.R. Acad.\ Sci.\ Paris},
{\bf169} (1919), 448--454.
\item{[16]} H. Iwaniec. Non-holomorphic modular forms and their
applications. in: R.A. Rankin (Ed.), {\it Modular Forms}, Ellis
Horwood, Chichester, 1984, pp.\ 157--196.
\item{[17]} H. Jacquet. Fonctions de Whittaker associ\'ees aux groupes
de Chevalley. {\it Bull.\ Soc.\ Math.\ France}, {\bf95} (1967),
243--309.
\item{[18]} H. Jacquet and R.P. Langlands. {\it Automorphic Forms on
$\r{GL}(2)$}. Lecture Notes in Math., {\bf 114}, Springer-Verlag,
Berlin, 1970.
\item{[19]} M. Jutila. Mean values of Dirichlet series via Laplace
transforms. in: Y. Motohashi (Ed.), {\it Analytic Number Theory},
Cambridge Univ.\ Press, Cambridge, 1997, pp.\ 169--207.
\item{[20]} A.W. Knapp. {\it Representation Theory of Semisimple Lie
Groups, An Overview based on Examples}. Princeton University Press,
Princeton, 1986.
\item{[21]} N.V. Kuznetsov. The Petersson conjecture for forms of
weight zero and the conjecture of Linnik. A mimeographed preprint,
Khabarovsk, 1977. (in Russian)
\item{[22]} ---. The Petersson conjecture for parabolic forms of
weight zero and the conjecture of Linnik. Sums of Kloosterman sums.
{\it Mat.\ Sb.}, {\bf111} (1980), 334--383. (in Russian)
\item{[23]} R.P. Langlands. {\it On the Functional Equations Satisfied
by Eisenstein Series}. Lecture Notes in Math., {\bf544},
Springer-Verlag, Berlin, 1976.
\item{[24]} H. Matumoto. Whittaker vectors and the Goodman--Wallach
operators. Acta Math., {\bf 161} (1988), 183--241.
\item{[25]} R. Miatello and N.R. Wallach. Automorphic forms
constructed from Whittaker vectors. {\it J.\ Funct.\ Anal.}, {\bf86}
(1989), 411--487.
\item{[26]} ---. Kuznetsov formulas for real rank one groups. {\it J.\
 Funct.\ Anal.}, {\bf93} (1990), 171--206.
\item{[27]} Y. Motohashi. An explicit formula for the fourth power
mean of the Riemann zeta-function. {\it Acta Math.}, {\bf 170}
(1993), 181--220.
\item{[28]} ---. The binary additive divisor problem. {\it Ann.\ Sci.\
\'Ecole Norm.\ Sup.}, (4) {\bf27} (1994), 529--572.
\item{[29]} ---. A relation between the Riemann zeta-function and the
hyperbolic Laplacian. {\it Ann.\ Scuola Norm.\ Sup.\ Pisa Cl.\ Sci.},
(4) {\bf 22} (1995), 299--313.
\item{[30]} ---. {\it Spectral Theory of the Riemann Zeta-Function}.
 Cambridge Univ.\ Press, Cambridge, 1997.
\item{[31]} ---. The mean square of Dedekind zeta-functions of
quadratic number fields. in: G.R.H. Greaves, G. Harman and M.N.
Huxley (Eds.), {\it Sieve Methods, Exponential Sums, and their
Applications in Number Theory}, Cambridge Univ.\ Press, Cambridge,
1997, pp.\ 309--324.
\item{[32]} ---. Trace formula over the hyperbolic upper half space,
in: Y. Motohashi (Ed.), {\it Analytic Number Theory}, Cambridge
Univ.\ Press, Cambridge, 1997, pp.\ 265--286.
\item{[33]} ---. New analytic problems over imaginary quadratic number
fields. in: M. Jutila M and T. Mets\"ankyl\"a (Eds.), {\it Number
Theory}, de Gruyter, Berlin, 2001, pp.\ 255--279.
\item{[34]} ---. A note on the mean value of the zeta and
$L$-functions.\ XII. Proc.\ Japan Acad., {\bf78}(A) (2002), 36--41.
\item{[35]} M.E. Picard. Sur un group de transformations des points de
l'espace situ\'es du m\^eme c\^ot\'e d'un plan. {\it Bull.\ Soc.\
Math. France}, {\bf12} (1884), 43--47.
\item{[36]} E.C. Titchmarsh. {\it The Theory of the Riemann
Zeta-Function}. Clarendon Press, Oxford, 1951.
\item{[37]} ---. {\it Introduction to the Theory of Fourier
Integrals}, Clarendon Press, Oxford, 1967.
\item{[38]} Ja.N. Vilenkin and A.U. Klimyk. {\it Representations of
Lie Groups and Special Functions}. Mathematics and Its Applications
(Soviet Series), Vol.\ 1, Kluwer, Amsterdam, 1991.
\item{[39]} A.I. Vinogradov and L.A. Takhtadjan. The zeta-function of
the additive divisor problem and the spectral decomposition of the
automorphic Laplacian. {\it Zap.\ Nauchn.\ Sem.\ LOMI}, {\bf 134}
(1984), 84--116. (Russian)
\item{[40]} N.R. Wallach. {\it Real Reductive Groups I}. Pure and
 Applied Mathematics, {\bf132}, Academic Press, New York, 1988.
\item{[41]} ---. {\it Real Reductive Groups II}. ibid., {\bf 132},
1992.
\item{[42]} G.N. Watson. {\it A Treatise on the Theory of Bessel
 Functions}. Cambridge Univ.\ Press, Cambridge, 1944. \vskip 1cm
\noindent
{\small Roelof W. Bruggeman
\par\noindent
Department of Mathematics, Utrecht University,
\par\noindent
P.O.Box 80.010, TA 3508 Utrecht, the Netherlands
\par\noindent
Email: bruggeman@math.uu.nl}
\medskip
\noindent
{\small Yoichi Motohashi
\par\noindent
Honkomagome 5-67-1-901, Tokyo 113-0021, Japan
\par\noindent
Email: am8y-mths@asahi-net.or.jp } \bye